\documentclass[smallextended,envcountsect]{svjour3}       
\smartqed  

\usepackage{appendix}
      \usepackage{amssymb}
      \usepackage{graphicx}

 \usepackage[usenames,dvipsnames]{pstricks}
 \usepackage{epsfig}
 \usepackage{pst-grad} 
\usepackage{pst-plot} 
 \usepackage[space]{grffile} 
 \usepackage{etoolbox} 
\usepackage{color}
\usepackage{dsfont}
\usepackage{amsmath}
\RequirePackage[OT1]{fontenc}
\RequirePackage{amsmath}
\RequirePackage[numbers]{natbib}
\RequirePackage[colorlinks,citecolor=blue,urlcolor=blue]{hyperref}

\spnewtheorem{thm}{Theorem}[section]{\bfseries}{\itshape}
\spnewtheorem{cor}[theorem]{Corollary}{\bfseries}{\itshape}
\spnewtheorem{lem}[theorem]{Lemma}{\bfseries}{\itshape}
\spnewtheorem{ntn}{Notations}[section]{\bfseries}{\itshape}
\spnewtheorem{pro}{Proposition}[section]{\bfseries}{\itshape}
\spnewtheorem{dfn}{Definition}[section]{\bfseries}{\itshape}
\spnewtheorem{as}{Assumption}[section]{\bfseries}{\itshape}
\spnewtheorem{rem}{Remark}[section]{\bfseries}{\itshape}
\spnewtheorem{ob}{Observation}[section]{\bfseries}{\itshape}

\numberwithin{equation}{section}

\begin{document}


\title{Polynomial ballisticity conditions and invariance principle for random walks in strong mixing environments
\thanks{E. Guerra was supported by CAPES PNPD20130824, Nucleus Millennium SMCDS NC130062 and CONICYT FONDECYT Postdoctorado 3180255. G. Valle was partially supported by CNPq grants 308006/2018-6 and 421383/2016-0. M. E. Vares was partially supported by  CNPq grant 305075/2016-0 and FAPERJ grant E-26/203.048/2016.}}


\titlerunning{Polynomial ballisticity conditions in mixing environments}        

\author{Enrique Guerra \and
        Glauco Valle \and
				Maria Eulália Vares
}

\authorrunning{E. Guerra, G. Valle and M. E. Vares} 

\institute{Enrique Guerra \at
              The Hebrew University of Jerusalem \\
              \email{eaguerra@mat.puc.cl}
           \and
           Glauco Valle \at
              Instituto de Matem\'atica - Universidade Federal do Rio de Janeiro \\
					    \email{glauco.valle@im.ufrj.br}
					\and
           Maria Eulália Vares\at
              Instituto de Matem\'atica - Universidade Federal do Rio de Janeiro \\
					    \email{eulalia@im.ufrj.br}  					
}

\date{Received: date / Accepted: date}

\maketitle

\begin{abstract}
We study ballistic conditions for $d$-dimensional random walks in strong mixing environments (RWRE), with underlying dimension $d\ge 2$. Specifically, we introduce an effective polynomial condition similar to that given by Berger, Drewitz and Ram\'{\i}rez (\emph{Comm. Pure Appl. Math.} {\bf 77} (2014) 1947--1973). In a mixing setup we prove this condition to imply the corresponding stretched exponential decay, and obtain an annealed functional central limit theorem for the  random walk process centered at velocity.

This paper complements previous work of Guerra (\emph{Ann. Probab.} {\bf 47} (2019) 3003--3054) and completes the answer about the meaning of condition $(T')|\ell$ in a mixing setting, an open question posed by Comets and Zeitouni  (\emph{Ann. Probab.} {\bf 32} (2004) 880--914).

\end{abstract}





{
  \hypersetup{linkcolor=black}
  \tableofcontents
}

\section{Introduction}
\label{secIntro}

Random walk in random environment (RWRE) is a widely studied stochastic model. However, its asymptotic behaviour is still poorly understood when the underlying dimension $d$ of the walk is larger than one, especially for the non-i.i.d. random environment case. The article \cite{Gue17}
opened a path to study ballistic properties of the random walk process under weaker assumptions than the usual Kalikow's condition. This condition intrinsically implies the ballistic nature of the walk (see \cite[Lemma 2.2]{SZ99}) and was the standing assumption in previous works such as \cite{CZ01}, \cite{CZ02} and \cite{RA03}, among others. The i.i.d. random environment case has already shown that Kalikow's condition is not the appropriate criterion in order to detect ballistic behaviour, indeed it is stronger (see \cite{Sz03}).

\vspace{0.5ex}
The main issue to establish a precise relation connecting local conditions on the environment law and a given asymptotic behaviour (which might be considered the main objective in this field) comes from the fact that under the annealed measure the random walk process is not Markovian. The non-Markovian character, already present when the environment is endowed with an i.i.d. structure, makes it harder to use ergodic devices to handle its asymptotic laws in terms of conditions on the random environment distribution. In turn, the loss of the Markov property is associated to the very definition of the annealed law as the semidirect product between quenched and environment laws.
As a result of this increasing complexity for more general environmental structures, the random walk process loses further properties. For instance, one can show that the multidimensional i.i.d. renewal structure of \cite{SZ99} fails to be a regeneration in mixing environments.

\vspace{0.5ex}
Our main objectives in this article are firstly to provide an effective condition under which ballistic behaviour is ensured. Secondly, to derive a proof of the equivalence between all the ballisticity conditions $(T^\gamma)|\ell$, for $\gamma \in (0,1)$ and $\ell\in\mathbb S^{d-1}$, where as a matter of definition
$$
\mathbb S^{d-1}:= \Big\{x\in\mathbb R^d:\, \sum_{i=1}^d x_i^2=1\Big\},
$$
therefore extending the original proof of Berger, Drewitz and Ram\'{\i}rez in \cite{BDR14} to mixing random environments.
These objectives are intertwined. Assuming the construction of an effective criterion equivalent to $(T^\gamma)|\ell$, the former aim will be essentially accomplished using diffusive controls on the fluctuations orthogonal to the \textit{quasi-asymptotic direction} of the walk process; these controls are provided under stretched exponential upper bounds on the probability of exiting large boxes through the unlikely faces. Renormalization procedures under those diffusive controls yield sharp inequalities for the environmental probability of unlikely quenched events. Following Sznitman's construction \cite{Sz00} in a mixing setting (see \cite [Section 5] {Gue17}), these yield appropriate estimates for renewal time tails. The proof will be concluded in these terms after an application of the functional CLT of \cite{CZ02}, which in turn, once the integrability conditions are fulfilled, depends only on the theory of Markov chains with infinite connections \cite{IG80}.

\vspace{0.5ex}
For the second aim we will follow the argument given in \cite{BDR14}, where the authors proved the so-called effective criterion starting from a class of polynomial conditions. Indeed, a first target for us shall be to construct an \textit{effective criterion}, denoted by $(EC)|\ell$, which will imply condition $(T')|\ell$. Afterwards, an exhaustive use of renormalization schemes shall produce the decay required by $(EC)|\ell$ under an \textit{effective polynomial condition} introduced in Section \ref{secpoly} and denoted by $(P_J)|\ell$. Alongside, the big picture behind the proof strategy in \cite{BDR14} works for our purposes. Nevertheless, we modify their proof argument due to minor problems related to split quasi independent environmental events, which indeed seem to be present even for
their i.i.d. setting. In Section \ref{secpoly}, we present a proof within the present more general framework of mixing random environments. We start with formal developments in the next paragraph.

\vspace{0.5ex}
A $d$-dimensional random walk in a random environment (RWRE) evolves on the lattice $\mathbb Z^d$, where the environment is specified by random transitions to nearest neighbors at each lattice site. More precisely, let $d\geq2$ (the case $d=1$ is well understood; see the survey \cite{Ze04}) be the underlying dimension,
$\kappa\in (0, 1/(4d))$ and consider the $(2d-1)$-dimensional simplex $\mathcal P_\kappa$ defined by
\begin{equation}\label{simplex}
\mathcal P_\kappa:=\left\{x\in\mathbb R^{2d}: \, \sum_{i=1}^{2d}x_i=1,\, x_i\geq 2\kappa, \,\,\mbox{  for  }i\in[1,2d] \right\}.
\end{equation}
Consider the product space $\Omega:=\mathcal P_\kappa^{\mathbb Z^d}$ which denotes the set of environment configurations, endowed with the canonical product $\sigma$-algebra $\mathfrak{F}_{\mathbb Z^d}$, where we fix an ergodic probability measure $\mathbb P$. We use the notation $|\,\cdot\,|_1$ and $|\,\cdot\,|_2$ to denote the $\ell_1$ and $\ell_2$ distances on $\mathbb R^d$, respectively; and furthermore, for $A, B\subset \mathbb Z^d$, $i\in \{1,2\}$, the notation $d_i(A,B)$ stands for the canonical $\ell_i$-distance between sets $A,\, B$, i.e., $d_i(A,B):=\inf\{|x-y|_i:\, x\in A, y\in B\}$. Also set $\Lambda=\{e \in \mathbb{Z}^d \colon |e|_1=1\}$. For a given environment $\omega:=(\omega(y, e))_{y\in \mathbb Z^d, e\in \mathbb Z^d, e \in \Lambda}=(\omega_y)_{y\in\mathbb Z^d} \in \Omega$, and  $x\in \mathbb Z^d$, the \textit{quenched law} $P_{x,\omega}$ is defined as the law of the canonical Markov chain $(X_n)_{n\geq 0}$ with state space $\mathbb Z^d$, starting from $x$, and with stationary transition probabilities given by the environment $\omega$, i.e.,
\begin{gather*}
  P_{x,\omega}[X_0=x]=1 \mbox{  and}\\
  P_{x,\omega}[X_{n+1}=X_n+e|X_n]=\omega(X_n,e),\,\, \mbox{  for  } e \in \Lambda.
\end{gather*}
The \textit{annealed law} $P_x$ is defined as the semidirect product probability measure $\mathbb P\otimes P_{x,\omega}$ on the product space $\Omega\times (\mathbb Z^d)^{\mathbb N}$. With some abuse of notation, we also denote by $P_x$ the marginal law of the random walk process $(X_n)_{n\geq0}$ under this measure.

  For an arbitrary set $A$ in a universal set $U$, let $A^c$ or $U\setminus A$ denote its complement.
Before introducing the class of mixing assumptions that will be considered, we recall the definition of $r$-Markovian fields:

\begin{definition} \label{rmarkov}
For $r\geq 1$ and $V \subset \mathbb{Z}^d$, let $\partial^r V=\{z\in V^c: d_1(z, V)\leq r\}$ be the (exterior) $r$-boundary of the set $V$. To simplify notation we will also write $\partial^1 V = \partial V$ for sets $V \subset \mathbb{Z}^d$. A random environment $(\mathbb P, \mathfrak{F}_{\mathbb Z^d})$ is called $r$-Markovian if for any finite $V\subset \mathbb Z^d$,
\begin{equation*}
  \mathbb P[(\omega_{x})_{x\in V}\in \cdot|\mathfrak{F}_{V^c}]=\mathbb P[(\omega_x)_{x\in V}\in \cdot|\mathfrak{F}_{\partial ^r V}], \,\, \mathbb P-a.s.,
\end{equation*}
where $\mathfrak{F}_{A}=\sigma(\omega_x, \, x\in A)$.
\end{definition}
We can now introduce a mixing condition on the environment which will be used throughout this article as standing assumptions:

\begin{definition} \label{def:smg}
Let $C$ and $g$ be positive real numbers. We say that an $r$-Markovian field $(\mathbb P, \mathfrak{F}_{\mathbb{Z}^d} )$ satisfies Guo's strong mixing condition \textbf{(SMG)}$_{C,g}$ if for all finite subsets $\Delta\subset V \subset \mathbb Z^d$ with $d_1(\Delta, V^c)\geq r$, and $A\subset V^c$,
\begin{equation}
\label{smg}
\frac{d\mathbb P[(\omega_x )_{x\in \Delta}\in \cdot | \eta]}{d \mathbb P[(\omega_x )_{x\in \Delta}\in \cdot | \eta']}\leq \exp\left( C \sum_{x\in \Delta, y \in  A}e^{-g |x-y|_1}\right) \,
\end{equation}
for $\mathbb P$-a.s. all pairs of configurations $\eta, \,\eta'\in \mathcal P_{\kappa}^{V^c} $ which agree over the set $V^c \backslash A$. Here we have used the notation
\begin{equation*}
\mathbb P[(\omega_x )_{x\in \Delta}\in \cdot | \eta]=\mathbb P[(\omega_x )_{x\in \Delta}\in \cdot |\mathfrak{F}_{V^c}]|_{(\omega_x)_{x\in V^c}=\eta}.
\end{equation*}
\end{definition}


In the RWRE context, condition  \textbf{(SMG)}$_{C,g}$ was used in the work of Guo \cite{Guo14}. It is in the spirit of what was considered by Dobrushin and Shlosman \cite{DS85}. Indeed, Lemma 9 in \cite{RA03} shows that it holds for strong-mixing Gibbs measures as in \cite{DS85}. This provide a class of interesting examples.

It is straightforward to see that the random walk process $(X_n)_{n\geq 0}$ under the annealed law $P_0$ is not Markovian. As an example, the first return time to the initial position is not a stopping time. In this mixing framework, worse things may happen  under
\textbf{(SMG)}$_{C,g}$. For instance, it can be checked that the renewal structure introduced in \cite{SZ99} fails to be a regeneration for any $g$. Indeed, it is unknown whether any regeneration structure exists at all.

\vspace{0.5ex}
For a subset $A\subset\mathbb Z^d$
we
introduce the random variable $T_A$ as
\begin{gather*}
T_A:=\inf\{n\geq 0:\, X_n\notin A\}
= \inf\{n\geq 0:\, X_n\in \partial A\}.
\end{gather*}
Following the original formulations in \cite{Sz02} for i.i.d. random environments and as a result of Lemma \ref{lemmaTgamma} in Section \ref{equivTgamma} for $\gamma\in (0,1]$ and  $\ell\in \mathbb S^{d-1}$, we define the so-called ballisticity conditions as follows.

\begin{definition}
\label{deftgammaandtprime}
We say that condition $(T^\gamma)|\ell$ holds if there exists some neighborhood $U\subset \mathbb S^{d-1}$ of $\ell$, such that
\begin{equation*}
\limsup_{\substack{L\rightarrow\infty}} L^{-\gamma}\ln\left(P_0\left[\widetilde{T}_{-bL}^{\ell'}<T_{L}^{\ell'}\right]\right)<0
\end{equation*}
for each $\ell'\in U$ and each $b>0$, where we have used the standard notation in the RWRE literature: for $c\in \mathbb R$ and $u\in \mathbb S^{d-1}$,
\begin{gather}\label{extimslab}
T_{c}^u:=\inf\{n\geq 0: X_n\cdot u\geq c\},\,\, \mbox{along with}   \\
\widetilde{T}_{c}^{u}:=\inf\{n\geq 0: X_n\cdot u\leq c \}.
\end{gather}
One then defines $(T')|\ell$ as the requirement that $(T^\gamma)|\ell$ holds for all $\gamma\in (0,1)$.
\end{definition}

Conditions requiring weaker decays than the above ones will be used in this paper.
\begin{definition}
\label{defpolasympandec}
For $J>0$ we say that condition $(P^*_J)|\ell$ holds if there exists some neighborhood $U\subset \mathbb S^{d-1}$ of $\ell$ such that
\begin{equation}
\lim_{\substack{L\rightarrow\infty}}L^{J}P_0\left[\widetilde{T}_{-bL}^{\ell'}<T_{L}^{\ell'}\right]=0,
\end{equation}
for all $\ell'\in U$ and all $b>0$. Furthermore, for a given rotation $R$ of $\mathbb R^d$ with $R(e_1)=\ell$ and real numbers $L, L', \widetilde{L}>3\sqrt{d}$, we define the \textit{box specification} $\mathcal B:= (R, L, L', \widetilde L)$ to which we associate the box
\begin{equation*}
B_{\mathcal B}:=R\left((-L,L')\times (-\widetilde L, \widetilde L)^{d-1}\right)
\end{equation*}
along with its \textit{frontal boundary part}
\begin{equation}\label{fb+}
\partial^+B_{\mathcal B}:=\partial B_{\mathcal B}\cap\{z: z\cdot \ell\geq L', |z\cdot R(e_j)|<\widetilde L\, \forall j\in [2,d]\},
\end{equation}
where $\partial = \partial^1$ as before. (With a little abuse of notation, we sometimes use $[a,b]$ to denote $[a,b] \cap \mathbb{Z}$, being always clear from the context.) We further introduce and attach to $\mathcal B$ the following random variables: 
\begin{gather}
\nonumber
  q_{\mathcal B}:=P_{0,\omega}\left[X_{T_{B_{\mathcal B}}}\notin \partial^+B_{\mathcal B}\right]=1-p_{\mathcal B}\,\, \mbox{ and} \\
  \label{enrique}
  \rho_{\mathcal B}:=\frac{q_{\mathcal B}}{p_{\mathcal B}}.
\end{gather}
\end{definition}

\smallskip

\begin{definition}
\label{effectiveee}
We say that the \textit{effective criterion} in direction $\ell$ holds, and denote it by $(EC)|\ell$, if
\begin{equation*}
  \inf_{\substack{\mathcal B, a}}C'\widetilde L^{d-1}L^{4(d-1)+1}\mathbb E\left[\rho_{\mathcal B}^{a}\right]<1,
\end{equation*}
where the infimum runs over $a\in (0,1]$, all the box specifications $(R,L-2,L+2, \widetilde L )$ with $R(e_1)=\ell$, $L^4>\widetilde L>L$,
$L>C''$, and where in turn $C',C''>3\sqrt d$ are certain prescribed dimensional constants.
\end{definition}

As it was anticipated, one of the advantages of this criterion is its effective character, i.e., it can be checked by inspection on a finite box.

\begin{remark} \leavevmode
\begin{itemize}
\item[(1)] Notice that the condition involved in Definition \ref{defpolasympandec} 
is not an effective condition. It is only instrumental, and will imply our effective polynomial condition to be given later (see Definition \ref{defpoly} in Section \ref{secpoly}). We start by this asymptotic condition, instead of the effective one, because it is simpler and displays better the connection with the conditions $(T^\gamma)|\ell$.
\item[(2)] $(EC)|\ell$ is similar to the effective criterion for i.i.d. environments introduced in \cite[Statement of Theorem 2.4]{Sz02}. The difference is in the exponent of $L$, since we allow the orthogonal dimensions to be of order $L^4$ instead of $L^3$.
That adjustment in the criterion helps to simplify the proof of Theorem \ref{mainth1} and it stems from necessary adjustments to make the criterion work for mixing environments.
\end{itemize}
\end{remark}

Our first result is a proof of the RWRE conjecture of Sznitman \cite{Sz02} in the present mixing setting.
\begin{theorem}
\label{mainth1}
Let $J>9d$, $\gamma\in (0,1)$, $C,g>0$ and $\ell\in \mathbb S^{d-1}$. Then the following statements are equivalent under 
\textbf{(SMG)}$_{C,g}$:
\begin{center}
\begin{itemize}
  \item $(T^\gamma)|\ell$ holds.
  \item $(T')|\ell$ holds.
  \item $(P^*_J)|\ell$ holds.
  \item $(EC)|\ell$ holds.
\end{itemize}
\end{center}
\end{theorem}

\begin{remark} For i.i.d. environments, the equivalence between any of the above ballisticity conditions and  condition $(T^1)|\ell$ (usually denoted as $(T)|\ell$) was proved in \cite{GR18}. This remains open in a mixing environment framework.
\end{remark}

We now state the result regarding the annealed invariance principle.

\begin{theorem}
\label{mainth2}
Let $C>0$, $J>9d$, $\ell\in \mathbb S^{d-1}$ with $g>2\ln(1/\kappa)$, and assume that the RWRE satisfies $(P^*_J)|\ell$ and 
\textbf{(SMG)}$_{C,g}$. Then there exist a deterministic non-degenerate covariance matrix $\mathcal{R}$ and a deterministic vector $v$ with $v\cdot \ell>0$ such that, setting
$$
S_n(t):=\frac{X_{[nt]}-vt}{\sqrt{n}},
$$
 the process $S_n(\cdot)$ with values in the space of right continuous functions possessing left limits endowed with the supremum norm, converges in law under $P_0$ to a standard Brownian motion with covariance matrix $\mathcal{R}$.
\end{theorem}

Notice that Theorem \ref{mainth1} requires an assumption on the mixing strength: $g>2\ln(1/\kappa)$. Of course $g$ can be taken arbitrarily large under the assumption that the environment is i.i.d..  The proof in Section \ref{secpoly} will explain our requirement on $g$. 
On the other hand, a detailed analysis of the proof will show that it actually requires just a little more than $J=9d-1$.

\vspace{0.5ex}
Let us explain the structure of the paper. Section \ref{secPrel} introduces primary properties about conditions $(T^\gamma)|\ell$,  including the statement that a suitable \textit{box version} is equivalent to the slab version given in Definition \ref{deftgammaandtprime}. The proofs are totally similar to those in \cite{Gue17}. For completeness they are given in an appendix. 
We also recall the approximate renewal structure of Comets and Zeitouni \cite{CZ01} and
we finish Section \ref{secPrel} stating that stretched exponential moments are finite under condition $(T^\gamma)|\ell$ for the maximun displacement of the random walk between regeneration times.
In Section \ref{sectionce} we show that $(EC)|\ell$ implies $(T')|\ell$.
Section \ref{secpoly} has the purpose of introducing a polynomial effective condition $(P_J)|\ell$ (see Definition \ref{defpoly}) under which $(EC)|\ell$ holds. A proof for Theorem \ref{mainth1} is provided at the end of that section. Finally in Section \ref{secproofmainth} we develop an appropriate machinery to get estimates for the tails of the approximate renewal time, and we prove Theorem \ref{mainth2}.

\smallskip
\noindent {\emph {Overview of the proofs}}

This is a brief description of the general strategy to prove Theorems \ref{mainth1} and \ref{mainth2}. It is quite informal, and we refer to corresponding sections for mathematical developments, after a quick description of some auxiliary results which might help the reader understand the general idea behind the proofs and avoid focusing on cumbersome computations. The rough plan is easy to explain:
\begin{enumerate}
  \item \label{1} Using condition $(P^*_J)|\ell$ we prove $(EC)|\ell$ in mixing random environments. This is done by first obtaining an effective polynomial condition $(P_J)|\ell$ from $(P^*_J)|\ell$, see Section \ref{secpoly}.
  \item \label{2} $(EC)|\ell$ implies condition $(T')|\ell$, proved in Section \ref{sectionce}. On the other hand, the converse implication  $(T')|\ell$ implies $(P^*_J)|\ell$ is straightforward. This will prove Theorem \ref{mainth1}.
  \item \label{3}Using condition $(T')|\ell$ we obtain several controls that at the end will provide suitable estimates for tails of the approximate regeneration time of Subsection \ref{section2}. The result is proven after an application of the CLT given in \cite{CZ02}.
\end{enumerate}
We explain quickly each point in the outline above.

\noindent
\begin{itemize}

\item The effective criterion $(EC)|\ell$ will imply $(T')|\ell$. The proof is carried out along a one-step-renormalization inequality, which is iterated to obtain stretched exponential decays. We shall develop the formal procedure in Section \ref{sectionce}.

\item The strategy to prove that $(P_J)|\ell$ implies  $(EC)|\ell$ is split into proving  Lemma \ref{lemmakpolybad} and Proposition \ref{propquenpol}. The combination of these two results shows that  $\mathbb E[\rho_{\mathcal B}^{a}]$ (recall notation in (\ref{enrique})) decays faster than any polynomial function on $L$, with a box specification $(R, L-2,L+2,4L^3)$ and $a\in (0,1)$ depends on $L$. This is the content of Section \ref{secpoly}.

\item Under $(T')|\ell$, we prove Proposition \ref{propatyquenest}, which combined with Lemma \ref{decomtail1} proves finite moments of all orders for the approximate regeneration time $\tau_1^{(L)}$, which follows from an estimate on the tail of the regeneration time given in Theorem \ref{thboundtails}. Based on arguments in \cite{CZ02} this implies the CLT. The purpose of Lemma \ref{decomtail1} is to show that under assumption $(T')|\ell$, the tails of $\tau_1$ could be estimated by controlling in turn, the probability of an atypical event, which is suitably bounded using Proposition \ref{propatyquenest}. These topics will be the main concern in Section \ref{secproofmainth}.
\end{itemize}

\smallskip

\section{Preliminaries} \label{secPrel}

In this section, we discuss some equivalent formulations for condition $(T^\gamma)$, and we also introduce an approximate renewal structure for random walks in mixing environments already used in \cite{Gue17}.

\subsection{On equivalent formulations for $(T^\gamma)$} \label{equivTgamma}

The goal of this section is to provide different formulations for conditions $(T^\gamma)|\ell$.
We begin by introducing the term \textit{directed systems of slabs}.

\begin{definition}
We say that $l_0,l_1, \ldots ,l_k\, \in \mathbb S^{d-1}$, $a_0=1, a_1, a_2, \ldots, a_k\,>0$, $b_0,b_1, \ldots,b_k>0 $ generate an $l_0$-directed system of slabs of order $\gamma$, when
\begin{itemize}
  \item $l_0, l_1,\ldots, l_k$ span $\mathbb R^d$,
  \item $\mathcal D=\{x\in \mathbb R^d: x\cdot l_0\in [-b_0,1], l_i\cdot x\geq -b_i, i\in[1,k]\}\,\subset \{x\in \mathbb R^d: l_i\cdot x<a_i, i\in[1,d]\}$,
  \item $\limsup_{\substack{L\rightarrow \infty}}L^{-\gamma}\,\ln P_0\left[\widetilde{T}_{-b_iL}^{l_i}<T_{a_iL}^{l_i}\right]<0$,
  for $i\in [0,k]$.
\end{itemize}
\end{definition}

For positive numbers $L$ and $L'$, we define the box
\begin{equation}
\label{generalboxes}
B_{L,L',l}(x):=
x+R\left(\left(-L,L\right)\times\left(-L',L'\right)^{d-1}\right)
\cap\mathbb{Z}^d,
\end{equation}
where $R$ is a rotation on $\mathbb R^d$, $x\in \mathbb Z^d$ and $l\in\mathbb S^{d-1}$ such that $R$ satisfies
$R(e_1)=l$ (the specific form of such a rotation is immaterial for our purposes).

We can prove:

\begin{lemma}\label{lemmaTgamma} For any $\gamma\in(0,1)$,
the following assertions are equivalent to condition $(T^\gamma)|_{l_0}$ in Definition \ref{deftgammaandtprime}:
\begin{itemize}
\item [(i)]  Together with $l_0$, there exists a choice of
$$
l_1, \ldots ,l_k\, \in \mathbb S^{d-1}, \ a_0=1, \ a_1, a_2, \ldots, a_k\,>0,\ b_0,b_1, \ldots,b_k>0
$$
that generate an $l_0$-directed system of slabs of order $\gamma$.
\item [(ii)]For some positive constants $b$, $\hat{r}$, and large $L$, there are finite subsets $\Delta_L\subset \mathbb Z^d$, with $0\in \Delta_L\subset \{x\in \mathbb Z^d: x\cdot l_0\geq -bL\}\cap \{x\in \mathbb R^d: |x|_2\leq \hat{r}L\}$ such that
    $$
    \limsup_{\substack{L\rightarrow \infty}}L^{-\gamma}\, \ln P_0\left[X_{T_{\Delta_L}}\notin \partial^+\Delta_L\right]<0,
    $$
where $\partial^+\Delta_L=\partial \Delta_L \cap \{x\in \mathbb R^d: x\cdot l_0 \geq L \}$.
\item [(iii)]For some $c>0$, one has
\begin{equation}\label{Tgammasquare}
\limsup_{\substack{L\rightarrow\infty}}\, L^{-\gamma}\ln P_0\left[X_{T_{B_{L, cL, l_0}(0)}}\notin \partial^+ B_{L, cL, l_0}(0)\right]<0.
\end{equation}
\end{itemize}
\end{lemma}

The proof follows a similar procedure as that of  \cite[Lemma 2.2]{Gue17}, where the case $\gamma=1$ is considered. For the sake of completeness we provide the proof in Appendix \ref{apx}.

\subsection{Approximate renewal structure}
\label{Apre}
We recall the approximate renewal structure of Comets and Zeitouni (see \cite{CZ01}). The results are stated without proofs which can be found in \cite[Section 2.2]{Gue17}.

Throughout this section we assume that $(T^\gamma)|\ell$ holds for some $\gamma\in (0,1)$ and a fixed $\ell\in \mathbb S^{d-1}$. Observe that one can assume the direction $\ell$ is such that there exists $h\in (0, \infty)$ with
\begin{equation}\label{rul}
l:=h\ell\in \mathbb Z^d.
\end{equation}
 As it was argued in \cite{Gue17}, this is not a further restriction. Recall that $\Lambda=\{e \in \mathbb{Z}^d \colon |e|_1=1\}$ and set
 $\mathcal{W}=\Lambda \cup \{0\}$,  let the probability measure $\overline P_0$ be defined as
\begin{equation*}
\overline P_0:=\mathbb P\otimes Q\otimes P_{0,\omega, \varepsilon}\,\,\,\, \mbox{on}\,\,\, \Omega\times (\mathcal W)^{\mathbb N} \times (\mathbb Z^d)^{\mathbb N},
\end{equation*}
where $Q$ is a product
probability measure on $(\mathcal W)^{\mathbb N}$ such that writing
$\varepsilon=(\varepsilon_1, \varepsilon_2, \,\ldots)\in (\mathcal
W)^{\mathbb N}$, we have for all $n \ge 1$,  $Q[\varepsilon_n=
e]=\kappa$ for all $e \in \Lambda$, and $Q[\varepsilon_n=0]=1-2d \kappa$.
Then for fixed random elements $\varepsilon\in (\mathcal W)^{\mathbb N}$ and
$\omega \in \Omega$, we define $P_{0,\omega, \varepsilon}$ as the
law of the Markov chain $(X_n)_{n\ge 0}$ with state space
$\mathbb Z^d$, starting from $0\in \mathbb Z^d$ and transition
probabilities
\begin{equation*}
  P_{0,\omega, \varepsilon}[X_{n+1}=X_n+e| X_n]=\mathds{1}_{\{\varepsilon_{n+1}=e\}}+\frac{\mathds{1}_{\{\varepsilon_{n+1}=0\}}}{1-2d\kappa}\left(\omega(X_n,e)-\kappa\right),
\end{equation*}
where $e \in \Lambda$.
An important property of this auxiliary probability space stems from the fact that the law of $(X_n)_{n\geq0}$
under $Q\otimes P_{0,\omega, \varepsilon}$ coincides with the law
under $P_{0,\omega}$, while its law under $\overline P_0$ coincides with that under $P_0$, as one easily verifies.


\vspace{1ex}
\noindent
Define now the sequence $\bar{\varepsilon}$ of length $|l|_1\in \mathbb N$ in
the following form:
$\bar{\varepsilon}_1=\bar{\varepsilon}_2=\ldots=\bar{\varepsilon}_{|l_1|}=\mbox{sign}(l_1)e_1$,
$\bar{\varepsilon}_{|l_1|+1}=\bar{\varepsilon}_{|l_1|+2}=\ldots=\bar{\varepsilon}_{|l_1|+|l_2|}=\mbox{sign}(l_2)e_2$,
$\ldots , \,
\bar{\varepsilon}_{|l|_1-|l_d|+1}=\ldots=\bar{\varepsilon}_{|l|_1}=\mbox{sign}(l_d)e_d$.
Notice that $l_i:=\ell\cdot e_i$ for $1\le i \le d$.

\noindent
For $\zeta>0$ small, $x\in \mathbb Z^d$, the cone
$C(x,l,\zeta)\subset \mathbb Z^d$ will be given by
\begin{equation}\label{cone}
C(x,l,\zeta):=\{y\in \mathbb Z^d: (y-x)\cdot l \geq \zeta |l|_2
|y-x|_2 \}.
\end{equation}
We assume that $\zeta$ is small enough in order to satisfy the
following requirement:
\begin{equation*}
  \bar{\varepsilon}_1,\bar{\varepsilon}_1+\bar{\varepsilon}_2,\, \ldots\, , \bar{\varepsilon}_1+\bar{\varepsilon}_2+\ldots+\bar{\varepsilon}_{|l|_1}\in C(0, l, \zeta).
\end{equation*}
For $L\in |l|_1\mathbb N$ we will denote by
$\bar{\varepsilon}^{(L)}$ the vector
\begin{equation*}
\bar{\varepsilon}^{(L)}=\overbrace{(\bar{\varepsilon}, \bar{\varepsilon},\, \ldots\, ,\bar{\varepsilon},\bar{\varepsilon})}^{L/|l|_1-\mbox{times}}
\end{equation*}
of length equal to $L$. Setting
\begin{equation*}
  D':=\inf\{n\geq0:\, X_n\notin C(X_0, l , \zeta)\},
\end{equation*}
we have:

\begin{lemma}\label{Dundert}
Assume $(T^\gamma)|\ell$ holds for some $\gamma \in (0,1)$, and fix a constant $\mathfrak{r}$ and a rotation $R$ as in item $(iii)$ of
Lemma \ref{lemmaTgamma}. Then there exists $c_1>0$ such that if $\zeta
<\min\left\{\frac{1}{9d}, \frac{1}{3d\mathfrak{r}}\right\}$, then
$$
P_0[D'=\infty]\geq c_1.
$$
\end{lemma}

The proof is exactly as that of \cite[Lemma 2.3]{Gue17}, which uses \cite[Proposition 5.1]{GR17}.

\vspace{0.5ex}
Let $\zeta>0$ satisfy the hypotheses of Lemma \ref{Dundert}. For each $L\in |l|_1\mathbb N$,
set $S_0=0$, and denoting by $\theta$ the canonical time shift, we set
\begin{gather*}
S_1=\inf\{n\geq L:\, X_{n-L}\cdot l >\max_{0\leq j<n-L}\{X_j\cdot
l\}, \,
(\varepsilon_{n-L+1},\ldots,\varepsilon_{n})=\bar{\varepsilon}^{(L)}\},\\
R_1=D'\circ \theta_{S_1}+S_1,
\end{gather*}
and for $n\ge 2$,
\begin{gather*}
S_n=\inf\{n>R_{n-1}:\, X_{n-L}\cdot l >\max_{0\leq j<n-L}\{X_j\cdot
l\}, \,
(\varepsilon_{n-L+1},\ldots,\varepsilon_{n})
=\bar{\varepsilon}^{(L)}\},\\
R_n=D'\circ \theta_{S_n}+S_n.
\end{gather*}
Given $L$ as above, these random variables are stopping times for the canonical filtration of the pair $(X_n,
\varepsilon_n)_{n\geq 0}$. Moreover
\begin{equation*}
  S_0=0<S_1\leq R_1\leq \ldots \leq S_n\leq R_n \ldots \leq \infty ,
\end{equation*}
with strict inequality whenever the left term is finite. Setting
\begin{equation*}
  K:=\inf\{n\geq 1: S_n<\infty, R_n=\infty\},
\end{equation*}
one defines the first time of
asymptotic regeneration $\tau_1:=\tau_1^{(L)}=S_K \leq \infty$. We suppress the dependence  of $\tau_1$ on $L$ when there is no risk of confusion. Notice that $\tau_1$ also depends on the direction $\ell$.


\vspace{0.5ex}
The next lemma shows that the construction makes sense.
\begin{lemma}\label{tranunderttau}
Assume $(T^{\gamma})|\ell$ holds for some $\gamma\in (0,1)$. Then $P_0$-a.s. (recall (\ref{rul}))
\begin{equation}\label{TTRANSIENT}
\lim_{\substack n\rightarrow \infty}\, X_n\cdot l=\infty.
\end{equation}
and there exists $L_0>0$ such that
for each $L\geq L_0$, with $L\in \,|l|_1\mathbb N$, one has $\overline P_0$-a.s.
\begin{equation}\label{tau1finite}
\tau_1^{(L)}<\infty.
\end{equation}
\end{lemma}

The claim (\ref{tau1finite}) is a straightforward application of
Lemma \ref{Dundert}. See \cite[Lemma 6.2]{GR17} for details. As
for the proof of claim (\ref{TTRANSIENT}), see
\cite[Lemma 6.1]{GR17} and \cite[page 517]{Sz02}. Choosing $L$ and $\zeta$ as prescribed by Lemmas \ref{Dundert}-\ref{tranunderttau}, one has that under $(T^\gamma)|\ell$, $\overline P_0$-a.s. $\{R_k<\infty\}\,=\,\{S_{k+1}<\infty\}$ and $S_1<\infty$ by (\ref{TTRANSIENT}).

We may also define the iterates:
\begin{equation*}
  \tau_n=\tau_1\circ \theta_{\tau_{n-1}}+\tau_{n-1}
\end{equation*}
for $n>1$ (by convention $\tau_0=0$). It is routine to verify that, under $(T^\gamma)|\ell$, $\overline P_0$-a.s. $\tau_k<\infty$ for any $k\in \mathbb N$.

\smallskip
\subsection{On the almost renewal structure for random walks in strong mixing environments}
\label{section2}
Our mixing assumptions provide an approximate
renewal structure when one considers the increments of the $\tau_k$
iterates. More precisely, given $j \in \mathbb{N}$, we define the $\sigma$-algebra $\mathcal
G_j$ by
\begin{equation}
\sigma\left(\omega(y, \cdot): y\cdot l< X_{\tau_j}\cdot l- (L|l|_2)/(|l|_1),\, (\varepsilon_i)_{0\leq i\leq \tau_j},\,(X_i)_{0 \leq i\leq \tau_j}\right).
\end{equation}
An important technical fact is the content of the next proposition which is \cite[Corollary 3.3]{Gue17}. 


\begin{proposition}
\label{propare}
Assume 
\textbf{(SMG)}$_{C, g}$.
For each $t\in (0,1)$ there exists $L_0:=L_0(C, g, \kappa, l,
d,r)$ such that $\overline P_0$-a.s.
\begin{align}
\label{approxre1}
&\exp\left(-e^{-g\, tL}\right)\overline{P}_0[(X_{n}-X_{0})_{n\geq0}\in \cdot\,|\, D'=\infty] \\
\nonumber
&\leq \overline P_0[(X_{\tau_j+n}-X_{\tau_j})_{n\geq0}\in \cdot \,|\,\mathcal G_j]\\
\nonumber
&\leq \exp\left(e^{-g\, tL}\right)\overline{P}_0[(X_{n}-X_{0})_{n\geq0}\in \cdot\,|\, D'=\infty]
\end{align}
holds for all $L\geq L_0$ with $L\in|l|_1 \mathbb N$.
\end{proposition}





\subsection{Further estimates under $(T^\gamma)$ in mixing environments}\label{secbasicpropTgamma}

Throughout this section, we let $\gamma \in (0,1)$, $\ell\in \mathbb S^{d-1}$ and choose $h\in(0,\infty)$ such that (\ref{rul}) is satisfied, i.e. $l=h\ell\in \mathbb Z^d $. The goal of this section is to establish that under $(T^\gamma)|\ell$ one obtains the existence of some finite stretched exponential moments for the regeneration position.

\vspace{0.5ex}
Thus, we consider for the vector $l$, the regeneration times $\tau_1^{(L)}$ constructed along $l$, which depend on $L\in|l|_1\mathbb N$ and we further assume that $L\geq L_0$, where $L_0$ is as in Proposition \ref{propare}. We recall that the construction of $\tau_1^{(L)}$ requires an underlying cone $C(x,l,\zeta)$, where $x\in\mathbb Z^d$ (see \eqref{cone}), and the cone angle $\zeta>0$ will be conveniently taken as a positive number satisfying:

\begin{equation}
\label{zeta}
\zeta<\min\left\{\frac{1}{9d},\, \frac{1}{3d\mathfrak{r}},\,\cos\left(\frac{\pi}{2}-\arctan(3\mathfrak{r})\right)\right\}.
\end{equation}

In turn, we have chosen the constant $\mathfrak{r}>0$ so that (\ref{Tgammasquare}) in item $(iii)$ of Lemma \ref{lemmaTgamma} is fulfilled. Let us remark that the form of requirement (\ref{zeta}) is used in the proof of Lemma \ref{expintM}.

\begin{proposition}
\label{expmpr}
Assume that $(T^\gamma)|\ell$ and 
\textbf{(SMG)}$_{C,g}$ hold. Then
there exist positive constants $c_2$, $c_3$ and $L_0$ such that
for all $L\geq L_0$, with $L\in |l|_1 \mathbb N$, we have that
\begin{equation}
\label{expmompos}
\overline{E}_0 \Big[\exp\left(c_2\left(\kappa^L X_{\tau_1} \cdot l\right)^\gamma\right) \Big] < c_3 .
\end{equation}
\end{proposition}


\noindent \emph{Proof.}
The proof is analogous to that of \cite[Proposition 4.1]{Gue17}. It is just a matter of following the proof using the subadditivity of $h(u) = u^\gamma$ at each step where the exponent is decomposed. \hfill $\square$



\medskip
\begin{lemma}
\label{expintM}
There exist constants $c_4$ and $c_5>0$ such that
\begin{equation}
\label{expbM}
E_0\left[\exp \left(c_4\left(\overline{M}\right)^\gamma\right),\, D'<\infty\right]<c_5.
\end{equation}
\end{lemma}
\noindent \emph{Proof.}
The proof is basically the same as that of \cite[Lemma 4.2]{Gue17}. The reader just has to carry the exponent $\gamma$ over the variables $\overline{M}$ and $M'$, where
\begin{equation*}
M':=\sup_{\substack{0\leq n \leq D'}}\{(X_n-X_0)\cdot \ell\} \, ,
\end{equation*}
for the entire proof. We simply remark that the assumption (\ref{zeta}) on the cone angle was conveniently used therein to guarantee that the renormalization procedure works. The same argument works here. \hfill $\square$

\medskip
The next corollary represents a reinforcement of Proposition \ref{expmpr}, to be used  in Section \ref{secproofmainth} (proof of Proposition \ref{protransfluctuation}). Its proof is completely analogous to that of \cite[Corollary 4.3]{Gue17}, and we omit it here.

\begin{corollary}\label{corexp}
Assume $(T^\gamma)|\ell$ and 
\textbf{(SMG)}$_{C,g}$. Let
\begin{equation}
\label{supt1}
X^*:=\sup_{\substack{0\leq n \leq \tau_1}}|X_n|_2.
\end{equation}
Then there exist positive constants $c_{6}$, $c_{7}$ and $L_0$ such that
$$
\sup_{L\geq L_0 \atop L\in |l|_1 \mathbb N} \overline{E}_0[e^{c_{6}\left(\kappa^L X^*\right)^\gamma}]\leq c_{7} \, .
$$
\end{corollary}


\section{An effective criterion for mixing environments}
\label{sectionce}

In this section we show that $(EC)|\ell$ implies $(T')|\ell$, thus extending the effective criterion of \cite{Sz02} to strongly mixing environments.
Indeed $(EC)|\ell$ will be equivalent to $(T')|\ell$ as stated in  Theorem \ref{mainth1} (see also the outline of its proof in Section \ref{secIntro}). We say that $(EC)|\ell$ is effective because it can be checked by inspection at the environment on a finite box, as opposed to the asymptotic character of the definition of $(T')|\ell$. Our approach is close to the one presented by Sznitman in \cite[Section 2]{Sz02}.

Let $\ell\in \mathbb S^{d-1}$ be a fixed direction and a rotation $R$ with $R(e_1)=\ell$. Recall the notation introduced in Definition \ref{defpolasympandec} with respect to $(\ell,R)$.
Fix positive numbers $L_0, L_1,\widetilde L_0$ and $\widetilde L_1$ satisfying
\begin{equation}\label{lolonenum}
3\sqrt{d}<L_0<L_1,\hspace{3ex}3\sqrt{d}<\widetilde L_0<\widetilde L_1,
\end{equation}
and the boxes $B_0$ and $B_1$ with specifications respectively $\mathcal B_0 = (R,L_0-1,L_0+1, \widetilde L_0 )$ and $\mathcal B_1 = (R,L_1-1,L_1+1,\widetilde L_1)$.
We now add a basic remark which will be used throughout this and the next sections.

\begin{remark}
Recalling the standard inequality $|x|_1 \leq \sqrt{d} |x|_2$ for $x \in \mathbb{R}^d$ and setting $\nu_1:= \sqrt{d}$, we have that for any pair of points $x, y\in \mathbb Z^d$ such that $|x-y|_2\leq \eta$, there exists a self avoiding and nearest neighbour path $[[x,y]]$ of length less or equal to $\nu_1\eta$. 
\end{remark}

\noindent
The next proposition is an extension of \cite [Proposition 2.1] {Sz02} to our setup. It provides an appropriate mixing estimate for moments of the random variable $\rho_1 = \rho_{\mathcal B_1}$ in terms of moments for the random variable $\rho_0 = \rho_{\mathcal B_0}$. 

\begin{proposition}\label{proprecurce}
Assume \textbf{(SMG)}$_{C,g}$.
There exist $c'_2>3\sqrt{d}$, $c'_3(d)$, $c_4'(d)$, $c_5'(d)$ and $c_6'(d)>1$ such that when $N_0:=L_1/L_0\geq 3$, $L_0\geq c_2'$, $\widetilde L_1\geq 48 N_0 \widetilde L_0$,  for any $a\in (0,1]$ one has
\begin{gather}
\nonumber
  \mathbb E\left[\rho_1^{a/4}\right]\leq c_3' \left\{(2\kappa)^{-10 \nu_1 L_1}\left(c_4'\widetilde L_1^{d-2}\frac{L_1^3}{L_0^2}\widetilde L_0\mathcal F_{SM}\mathbb E\left[q_0\right]\right)^{\frac{\widetilde L_1}{16 N_0\widetilde L_0}}\right. \\
\label{critrenor}
  \left.+\sum_{0\leq m\leq N_0+1}\left(F_{SM}\, c_4'\widetilde L_1^{d-1}\mathbb E\left[\rho_0^a\right]\right)^{\frac{N_0+m-1}{4}}\right\},
\end{gather}
where $\mathcal F_{SM}$ and $F_{SM}$  are positive functions on $L_0$, $\widetilde L_0$, $L_1$ and $\widetilde L_1$ which satisfy:
\begin{align*}
\mathcal F_{SM} &\leq \exp\left(c_5'(\widetilde L_1)^2 (L_1)^2 (\widetilde L_1)^{2(d-2)}\exp\left(-2gN_0\widetilde L_0\right)\right) \mbox{  and}\\
  F_{SM} &\leq \exp\left( c_6' (\widetilde L_1)^{2{d-1}}(L_0)^2N_0\exp\left(-\frac{3gL_0}{2}\right)\right).
\end{align*}
\end{proposition}

\noindent \emph{Proof.} It will be convenient to introduce thin slabs transversal to $\ell$. More precisely, for $i\in \mathbb Z^d$, we define:
\begin{equation}\label{Hi}
\mathcal H_i:=\left\{y\in \mathbb Z^d: \, \exists \, x\in \mathbb Z^d, \, |x-y|_1=1,\,(x\cdot\ell -iL_0)(y\cdot\ell -iL_0)\leq 0 \right\},
\end{equation}
and let $I:\mathbb Z^d\mapsto \mathbb Z$ be given by
\begin{equation}\label{I}
I(z)=i \text { if and only if  }  z\cdot \ell \in [iL_0-L_0/2, iL_0+L_0/2).
\end{equation}
We note that $I(z)=i$ {\blue for} $z \in \mathcal H_i$, since $L_0>2$. Now define the successive visit times to the slabs $\mathcal H_i$:
\begin{gather}
\nonumber
V_0=0,\,\, V_1=\inf\{n\geq 0,\,X_n\in\mathcal H_{I(X_0)+1}\cup \mathcal H_{I(X_0)-1} \}, \,\, \mbox{ and}\\
\label{visitHi}
V_i=V_1\circ \theta_{V_{i-1}}+V_{i-1}\, \mbox{ for $i\ge 2$.}
\end{gather}
For a given environment $\omega \in \Omega$ and $x\in \mathbb Z^d$, set $\widehat q_{x,\omega}$ and $\widehat p_{x,\omega}$ as:
\begin{equation}\label{slabqp}
  \widehat q_{x,\omega}=P_{x,\omega}\left[X_{V_1}\in \mathcal H_{I(X_0)-1}\right]=1-\widehat p_{x,\omega},
\end{equation}
and
\begin{equation}\label{rhohat}
  \widehat \rho_{i,\omega}:=\sup\left\{\frac{\widehat q_{x,\omega}}{\widehat p_{x,\omega}},\,\,x\in \mathcal H_i,\,\, |x\cdot R(e_j)|<\widetilde L_1, \, j\in[2,d]\right\},
\end{equation}
for $i\in \mathbb Z$.
Let $n_0:=\lfloor N_0 \rfloor$ and consider the nonnegative function $f$ on $\{n_0+2,n_0+1,\ldots \}\times\Omega$  defined by:
\begin{gather}
\nonumber
f(n_0+2,\omega)=0, \,\, \mbox{ and} \\
\label{ffunc}
f(i,\omega)=\sum_{i\leq j \leq n_0+1}\, \prod_{j< m\leq n_0+1}(\widehat \rho_{m,\omega})^{-1},\,\, \mbox{ for $i \le n_0+1$}.
\end{gather}
We shall drop the $\omega$-dependence from these random variables when $\omega$ is fixed in the estimates. Introducing the stopping time $\widehat T$ as
\begin{equation}\label{ttilderen}
\widehat T:=\inf\{n\geq 0,\,\, |X_n\cdot R(e_j)|\geq \widetilde L_1,\, \mbox{ for some }j\in [2,d] \} \, ,
\end{equation}
one has the quenched estimate:
\begin{equation}\label{quenf}
P_{0,\omega}\left[\widetilde{T}_{-L_1+1}^\ell<\widehat T\wedge T_{L_1+1}^\ell\right]\leq \frac{f(0)}{f(-n_0+1)} \, ,
\end{equation}
which is proved in \cite[proof of Proposition 2.1 - page 524]{Sz02}.

We will also need an estimate for the annealed probability $P_0[\widehat T< \widetilde T_{-L_1+1}^\ell \wedge T_{L_1+1}^\ell]$. For this end, we introduce for $z\in \mathbb R^d$ the semi-norm 
\begin{equation} \label{nperp}
|z|_{\perp}:=\sup\{z\cdot R(e_j), j\in [2,d]\}
\end{equation}
as well as the stopping times for $j\in[2,d]$ and $u > 0$
\begin{gather}
\nonumber
\sigma_{u}^{j,+}=\inf\left\{n\geq 0, \,\, X_n\cdot R(e_j)\geq  u\right\}\,\, \mbox{ and}\\
\label{sigmatimes}
\sigma_{u}^{j,-}=\inf\left\{n\geq 0, \,\, X_n\cdot R(e_j)\leq -u\right\}.
\end{gather}
We set
\begin{equation}\label{escalejlo}
  M=\left\lfloor \frac{\widetilde N_0}{2(n_0+1)}\right\rfloor, \,\, \mbox{with }\,\, \widetilde{N}_0=\frac{\widetilde L_1}{\widetilde L_0},\,\, \mbox{ and }\,\,\overline L_0=2(n_0+1)\widetilde L_0 \, .
\end{equation}
 From the hypotheses, $M\geq \lfloor 24 N_0/(n_0+1)\rfloor \geq 18$. Since $M\overline L_0<\widetilde L_1$, on the event $\{\widehat T\leq\widetilde T_{-L_1+1}^\ell\wedge T_{L_1+1}^\ell \}$, one gets that $P_0$-a.s. there exists $n\in \mathbb N$ such that $|X_n\cdot R(e_j)|\geq M\overline L_0$ for some $j\in [2,d]$ and for all $i<n$ the walk is inside box $B_1$. Thus:
\begin{equation}\label{pttilde}
P_0[\widehat T\leq \widetilde T_{-L_1+1}^\ell\wedge T_{L_1+1}^\ell] \leq \sum_{2\leq j\leq d} \big( P_0[\sigma_{M\overline L_0}^{+,j}\leq T_{B_1}]+P_0[\sigma_{M\overline L_0}^{-, j}\leq T_{B_1}] \big) .
\end{equation}

Let us simply write $\sigma_u$ by $\sigma_{u}^{+,2}$ and prove an estimate for $P_0[\sigma _{M\overline L_0}\leq T_{B_1}]$. The other terms in (\ref{pttilde}) can be treated similarly, and we obtain analogous bounds for them. We let $k\in \mathbb Z$ and define the cylinder
\begin{equation*}
c_\perp(k)=\{w\in\mathbb R^d: w\cdot R(e_2)\in [k\overline L_0,(k+1)\overline L_0), \,\, |w\cdot R(e_i)|<\widetilde L_1,\, \forall i\in [3,d]\}
\end{equation*}
together with the discrete cylinder 
\begin{equation*}
  \overline c(k)=\{x\in \mathbb Z^d: \, \inf_{w\in c_\perp(k)}|x-w|\leq (2n_0+1)\widetilde L_0,\,x\cdot l\in (-L_1+1, L_1+1) \}.
\end{equation*}

Applying the strong Markov property on the \textit{quenched probability} $P_{0,\omega}$, we have
\begin{equation}\label{sigmdes}
P_0[\sigma_{M\overline L_0}\leq T_{B_1}]\leq E_0[\sigma_{(M-2)\overline L_0}<T_{B_1}, P_{X_{\sigma_{(M-2)\overline L_0}},\omega}[\sigma_{M\overline L_0}\leq T_{B_1}]].
\end{equation}
From  (\ref{escalejlo}) and the definitions of $c_\perp(k)$ and $\overline c(k)$, on
$\{\sigma_{(M-2)\overline L_0}<T_{B_1}\}$, one has that $P_0$-a.s the random variable $X_{\sigma_{(M-2)\overline L_0}}$ is in $\overline c(M-2)$, since $2n_0+1<2(n_0+1)$. Recall that $\theta$ denotes the canonical time shift. Denoting by $H^i$, with $i\geq 0$ the iterates of the stopping time $H^1= T_{B_1}\wedge T_{X_0+B_0}$ and defining the stopping time
\begin{gather*}
S=\inf\{k\geq 0, (X_k-X_0)\cdot \ell \leq -L_0+1 \,\,\mbox{or}\,\, |(X_k-X_0)\cdot R(e_j))|\geq \widetilde L_0, \\ \mbox{for some}\,\, j\in [2,d] \},
\end{gather*}
we note that when $y\in c_\perp(M-2)\cap B_1$,  for $\omega \in \Omega$ one has that $P_{y,\omega}$-a.s.
\begin{equation}
\label{evento}
\bigcap_{i=0}^{2n_0+1}\theta_{H^i}^{-1}\{H^1<S\}\subset\{T_{B_1}<\sigma_{M\overline L_0}\}
\end{equation}
and the part of the path $X_\cdot$ involved in the event of the left on \eqref{evento} is contained in $\overline c(M-2)$  $P_{y,\omega}$-a.s.

Therefore, using the notation attached to the box specification $\mathcal B_0$ as well as denoting by $t$ the space shift on the environment $\Omega$, for $\omega \in \Omega$, $y\in \overline c(M-2)$ one has
$P_{y,\omega}$- a.s. for any $i\in [0,2n_0+1]$ the inequality
$$
P_{X_{H^i},\omega}[H^1<S]\geq\inf_{x\in \overline c(M-2)} p_{\mathcal B_0}\circ t_x(\omega).
$$
Thus an application of the strong Markov property gives
\begin{gather*}
P_{0}\mbox{-a.s. on the event }\{\sigma_{(M-2)\overline L_0}<T_{B_1}\},\\
P_{X_{\sigma_{(M-2)\overline L_0}},\omega}[T_{B_1}<\sigma_{M\overline L_0}]\geq\left\{\inf_{x\in \overline{c}(M-2)}p_{\mathcal{B}_0}\circ t_x(\omega)\right\}^{2(n_0+1)}\stackrel{def}=\varphi(M-2,\omega).
\end{gather*}
Inserting the previous estimate into (\ref{sigmdes}), we have
\begin{equation*}
P_0[\sigma_{M\overline L_0}\leq T_{B_1}]\leq \mathbb E[P_{0,\omega}[\sigma_{(M-4)\overline L_0}<T_{B_1}](1-\varphi(M-2,\omega)).
\end{equation*}
 Contrasting with \cite{Sz02}, the choice of $M-4$ in the above argument follows from the need of a larger separation between the involved regions, allowing to use the mixing hypothesis. A crucial point stems from the fact that the random variable $P_{0,\omega}[\sigma_{(M-4)\overline L_0}<T_{B_1}]$ is measurable with respect to
$$\sigma \left(\omega(y,\cdot),\, y \cdot R(e_2)< (M-4)\overline L_0, \,y\in B_1\right),$$
while $(1-\varphi(M-2,\omega))$ is measurable with respect to the  $\sigma$-algebra generated by the $\omega(y, \cdot)$ such that $y$  satisfies:
\begin{gather*}
(M-3)\overline L_0 \leq y \cdot R(e_2)< M \overline L_0, -(L_1+L_0)+2< y\cdot \ell < (L_1+L_0)+2\\ \mbox{ and } |y \cdot R(e_i)|<\widetilde L_1+ \overline L_0 \mbox { for all }\, i\in [3,d].
\end{gather*}

Let $A, B\subset \mathbb Z^d$ be the corresponding sets of sites on which the transitions need to be known to compute $P_{0,\omega}[\sigma_{(M-4)\overline L_0}<T_{B_1}]$ and $(1-\varphi(M-2,\omega))$, respectively. When $L_0>3r$ we have (recall the definitions of constants $r$ in Definition \ref{rmarkov} and $C$ and $g$ in Definition \ref{def:smg})
\begin{gather*}
\sum_{x\in A, y\in B} C\exp(-g|x-y|_1)\leq C|A||B|\exp\left(-g\overline L_0\right)\\
\leq c(d)(M \overline L_0)^2 (L_1)^2 (\widetilde L_1)^{2(d-2)}\exp\left(-g\overline L_0\right)
\end{gather*}
for some suitable $d$-dependent constant $c(d)$.
To keep control of the mixing assumption \textbf{(SMG)}$_{C,g}$ 
let
\begin{equation} \label{calF}
\mathcal F_{SM}:=e^{\sum_{x\in A, y\in B} C\exp(-g|x-y|_1)}.
\end{equation}
Applying \textbf{(SMG)}$_{C,g}$, we get
$$
P_0[\sigma_{M \overline L_0}\leq T_{B_1}]\leq \mathbb P_0[\sigma_{(M-4)\overline L_0}<T_{B_1}] \times
\mathbb E[(1-\varphi(M-2),\omega)]\times \mathcal F_{SM}.
$$
On the other hand, it can be seen that the next upper bound:
$$
\mathbb E[(1-\varphi(M-2,\omega))]\leq \widetilde c(d)\widetilde L_1^{(d-2)}\frac{L_1^3}{L_0^2}\widetilde L_0\mathbb E[q_0]
$$
holds (see \cite[page 526]{Sz02}) for a suitable constant $\widetilde c(d)>0$. Therefore, going back to \eqref{sigmdes} we get
\begin{gather*}
P_0[\sigma_{M \overline L_0}< T_{B_1}]\leq P_0[\sigma_{(M-4)\overline L_0}<T_{B_1}] \Big(
\mathcal F_{SM}\, \widetilde c(d)\widetilde L_1^{(d-2)}\frac{L_1^3}{L_0^2}\widetilde L_0\mathbb E[q_0] \Big).
\end{gather*}
Thus, a recursive argument lead us to
\begin{gather*}
P_0[\sigma_{M \overline L_0}< T_{B_1}]\leq
\left(\mathcal F_{SM}\widetilde c(d)\widetilde L_1^{(d-2)}\frac{L_1^3}{L_0^2}\widetilde L_0\mathbb E[q_0]\right)^{m},
\end{gather*}
for any $0\leq m \leq \lfloor M/4 \rfloor$. Since $\lfloor M/4 \rfloor \geq \widetilde L_1/(16 N_0\widetilde L_0)$, using \eqref{pttilde}, we find that
\begin{equation} \label{finesttransv}
P_0[\widehat T<\widetilde T_{-L_1+1}^\ell \wedge T_{L_1+1}^\ell]\leq 2(d-1)
\left(\mathcal F_{SM}\widetilde c(d)\widetilde L_1^{(d-2)}\frac{L_1^3}{L_0^2}\widetilde L_0\mathbb E[q_0]\right)^{\frac{\widetilde L_1}{16 N_0\widetilde L_0}}.
\end{equation}
From the choice for $\kappa$ in (\ref{simplex}), we have $\rho_1\leq (2\kappa)^{-2\nu_1L_1}$. Denote by $W$ the random variable
$$
P_{0,\omega}[\widetilde T_{-L_1+1}^\ell<\widehat T\wedge T_{L_1+1}^\ell]+P_{0,\omega}[\widehat T<\widetilde T_{-L_1+1}^\ell \wedge T_{L_1+1}^\ell]
$$
and consider the event
$$
\mathcal G=\{\omega\in \Omega:\, P_{0,\omega}[\widehat T<\widetilde T_{-L_1+1}^\ell \wedge T_{L_1+1}^\ell]\leq (2\kappa)^{9\nu_1 L_1} \}.
$$
Note that on $\mathcal G$ we have
\begin{equation*}
\rho_1 =  \rho_{\mathcal B_1} \leq \frac{W}{(1-W)_+}\stackrel{(\ref{quenf})}\leq \frac{f(0)+f(-n_0+1)(2\kappa)^{9\nu_1L_1}}{\left(f(-n_0+1)-f(0)-f(-n_0+1)(2\kappa)^{9\nu_1L_1}\right)_+}.
\end{equation*}
Replacing the constants $\kappa$ by $2\kappa$, the same argument as the one given in \cite[page 527]{Sz02} allows us conclude that on the event $\mathcal G$
\begin{equation}\label{quenestonG}
\rho_1(\omega)\leq 2\sum _{0\leq m\leq n_0+1}\prod_{-n_0+1<j\leq m}\widehat \rho_{j,\omega},
\end{equation}
provided $L_0\geq c_2'$, where $c_2'$ is a suitable dimensional dependent positive constant. Notice that when $i\in \mathbb Z$, we have for any $x\in \mathcal H_i$ the inequality:
$$
\frac{\widehat q_{x,\omega}}{\widehat p_{x,\omega}}\leq \rho_0\circ t_x.
$$
As a result, for $i\in \mathbb Z^d$ recalling \eqref{rhohat} and writing $\mathcal Y_i$ as the set   $\{x: |x\cdot R(e_j)|<\widetilde L_1 \ \forall \, j\in [2,d]\}\cap \mathcal H_i$, we obtain
\begin{equation}\label{rhodotilde}
\widehat \rho_{i}\leq \sup_{z\in \mathcal Y_i}\rho_0\circ t_z\stackrel{def}=\overline{\rho}_{i,\omega}.
\end{equation}
We now insert the above estimate \eqref{rhodotilde} into \eqref{quenestonG}. Then, fixing $a\in(0,1]$ and using the inequality $(w_1+w_2)^\gamma \leq w_1^\gamma+w_2^\gamma$ for $w_1,w_2\geq0$ and $\gamma\in(0,1)$, we obtain
\begin{equation}\label{estonG1}
\mathbb E\left[(\rho_1)^{\frac{a}{4}}(\omega),\mathcal G\right]\leq 2\sum _{0\leq m\leq n_0+1}\mathbb E\left[\prod_{-n_0+1<j\leq m}(\overline \rho_{j,\omega})^{\frac{a}{4}}\right].
\end{equation}
We now split each product entering at \eqref{estonG1} into four groups depending on the residues modulo $4$. For integer $j\in[0,3]$ we call $\mathcal M_j=\{i\in \mathbb Z: i=j\,\, (\mbox{mod }4)\}$. For $m\in[0,n_0+1]$ integer, applying Cauchy-Schwarz inequality twice we find that:
\begin{equation}
\label{estonG2}
\mathbb E\left[\prod_{-n_0+1<j\leq m}(\overline \rho_{j,\omega})^{\frac{a}{4}}\right] \leq
 \prod_{0\leq i\leq 3}\mathbb E^{\frac{1}{4}}\left[\prod_{\substack{-n_0+1<j\leq m\\ j\in \mathcal M_i}}(\overline \rho_{j,\omega})^a\right].
\end{equation}
Under
\textbf{(SMG)}$_{C,g}$, given $i\in [0,3]$, if $j_1, j_2\in \mathcal M_i$ with $j_1< j_2$, and if $L_0\geq 3r$, we may argue similarly to what was done before \eqref{finesttransv},  to obtain
\begin{equation}\label{decoupling}
\mathbb E\left[ \big(\prod_{-n_0+1<j\leq j_1} (\overline \rho_{j,\omega})^a \big) \, (\overline \rho_{j_2,\omega})^a\right]\leq \mathbb E\left[ \prod_{-n_0+1<j\leq j_1} (\overline \rho_{j,\omega})^a\right]F_{SM}\mathbb E\left[(\overline \rho_{j_2,\omega})^a\right].
\end{equation}
The term $F_{SM}$ in \eqref{decoupling} is defined as a mixing error as follows. We first set $j_i=\max\{k: \, k\in \mathcal M_i\cap (-n_0+1,n_0+1]\}$. Then, for each integer $j$ and $i\in[0,3]$, we define
\begin{align*}
A_{i,j}:=&\big\{z\in \mathbb Z^d: \forall k\in [2,d]\, |z\cdot R(e_k)|<\widetilde L_1+\widetilde L_0,\,\  \\
&\quad \quad  z\cdot \ell \in ((j-1)L_0, (j+1)L_0+3) \big\},
\end{align*}
together with
\begin{align}
\nonumber
B_i:=&\bigcup_{j \in \mathcal M_i\cap(-n_0+1,n_0+1]\setminus\{j_i\}}A_{i,j} .
\end{align}
In these terms, the factor $F_{SM}$ prescribed in (\ref{decoupling}), is given by:
\begin{align}
\nonumber
& \!\!\!\!\!\!\!\!\! e^{\sum_{x\in A_{i,j_i}, y \in B_i}\, C\exp\left(-g|x-y|_1\right)} \\
\label{mix2}
\leq& \exp\left(\mathfrak c (d)(\widetilde L_1)^{2{d-1}}(L_0)^2N_0\exp\left(-\frac{3gL_0}{2}\right)\right)
\end{align}
for certain dimensional constant $\mathfrak  c =\mathfrak c(d)>0$, uniformly on $i$.


Therefore for $i\in [0,3]$ and $m\in (-n_0+1,n_0+1]$, successive conditioning in each term on the right of (\ref{estonG2}) from the biggest $j$ \textit{backwards along direction} $\ell$, turns out
$$
\mathbb E^{\frac{1}{4}}\left[\prod_{\substack{-n_0+1<j\leq m\\ j\in \mathcal M_i}}(\overline \rho_{j,\omega})^a\right]\leq \prod_{\substack{-n_0+1<j\leq m\\ j\in \mathcal M_i}}\left(F_{SM}\mathbb E[(\overline \rho_{j,\omega})^{a}]\right)^{\frac{1}{4}}.
$$
Combining this last inequality in \eqref{estonG1} and (\ref{estonG2}) and the bound in (\ref{finesttransv}) with the help of Chebyshev's inequality we get
\begin{gather}
\nonumber
\mathbb E[(\rho_1)^{\frac{a}{4}}]\leq (2\kappa)^{-a\nu_1L_1}\mathbb P[\mathcal G^c]+2\sum_{0\leq m\leq n_0+1}\prod_{-n_0+1<j\leq m}\left(F_{SM}\,\mathbb E[(\overline \rho_{j,\omega})^a]\right)^{\frac{1}{4}}\\
\nonumber
\stackrel{(\ref{finesttransv})} \leq (2\kappa)^{-10\nu_1L_1}\left(\widetilde c\widetilde L_1^{d-2}\frac{L_1^3}{L_0^2}\widetilde L_0\mathcal F_{SM}\,\mathbb E[q_0]\right)^{\frac{\widetilde L_1}{16 N_0\widetilde L_0}}+\\
\label{endpropec}
2\sum_{0\leq m\leq n_0+1}\prod_{-n_0+1<j\leq m}\left(F_{SM}\mathrm c \widetilde L_1^{d-1}\mathbb E[\rho_0^{a}]\right)^{\frac{1}{4}},
\end{gather}
for some $\mathrm c>0$. The claim of the proposition follows. \hfill $\square$

\medskip

The previous proposition is instrumental and we apply it recursively along fast growing scales. For this, we let
\begin{equation}\label{u0valpha}
u_0\in (0,1),\; v=16,\; \alpha=320,
\end{equation}
and consider sequences of positive numbers $(L_k)_{k\geq 0}$ and $(\widetilde L_k)_{k\geq 0}$ satisfying:
\begin{gather}
\nonumber
L_0\geq c'_2 \, , \ \ L_0\leq \widetilde L_0 \leq L_0^4, \ \ \mbox{and} \\
\label{escalesL_k}
L_{k+1}=N_k L_k \, , \ \ \widetilde L_{k+1}=N_k^4\widetilde L_k \, , \ \ \mbox{with } N_k=\frac{\alpha \nu_1}{u_0}v^k \, , \ \ k \ge 1.
\end{gather}
Notice that from \eqref{escalesL_k}, one has for $k\geq 0$
\begin{equation}\label{expescaleL_k}
  L_k=\left(\frac{\alpha \nu_1}{u_0}\right)^k\, v^{\frac{k(k-1)}{2}} \, L_0  \quad \mbox{and} \quad \widetilde L_k=\left(\frac{L_k}{L_0}\right)^4\widetilde L_0.
\end{equation}
Then consider boxes $B_k$ with specifications $\mathcal B_k = \mathcal B(L_k-1,L_k+1, \widetilde L_k)$ for $k\geq 0$, and we set $\rho_k := \rho_{\mathcal B_k}$.

 Replacing $(L_0,L_1, \widetilde L_0, \widetilde L_1)$ by $(L_{k},L_{k+1}, \widetilde L_{k}, \widetilde L_{k+1})$ in the definitions of the events involved in the estimate \eqref{critrenor}, we obtain a sequence of constants $\mathcal F_{SM}(k)$ and $F_{SM}(k)$, $k\ge 0$, according to \eqref{calF} and \eqref{mix2}. We need appropriate upper bounds under the scales displayed in \eqref{escalesL_k} for $\mathcal F_{SM}(k)$ and $F_{SM}(k)$, provided that $L_0\geq \mathfrak c$ for a suitable constant $\mathfrak c>0$ to be determined. Observe that as a result of both expressions in \eqref{expescaleL_k}, it is clear that for some constants $c_7'(d, C,g),c_8'(d,C,g)>0$ and $\mathfrak{c}_9(d,C,g)>0$, when $L_0\geq \mathfrak{c}_9$ for all $k\geq 1$:
\begin{equation}
\label{c7c8}
c_4'\mathcal F_{SM}(k)\leq c_7',\, \hspace{0.7ex}\mbox{and }\hspace{0.7ex}c_4' F_{SM}(k)\leq c_8'.
\end{equation}
Furthermore, notice that for $k=0$ we have
\begin{equation*}
\mathcal F_{SM}(0)\leq \exp\left(c_5'\left(\frac{\alpha \nu_1 L_0}{u_0}\right)^{8d+2}e^{-2g\left(\frac{\alpha \nu_1}{u_0}\right)L_0}\right)\leq c_7',
\end{equation*}
whenever $L_0\geq \mathbf{c}_9$ for a suitable constant $\mathbf{c}_9>0$. On the other hand, an upper bound for $F_{SM}(0)$ can be obtained as follows
\begin{gather*}
F_{SM}(0) \leq \exp\left(c_6'\left(\frac{\alpha \nu_1}{u_0}\right)^{8d-3}L_0^{8d-2}e^{-\frac{3gL_0}{2}}\right)
\leq \exp\left(u_0^{-(8d-3)}e^{-gL_0}\right)\leq c_8',
\end{gather*}
provided that $L_0\geq \widetilde c(d)$ and $u_0\in [e^{-\frac{gL_0}{8d-3}},1]$. As a result, under the choice of scales given in (\ref{escalesL_k}), Proposition \ref{proprecurce} can be reformulated to get rid of the mixing terms: $\mathcal F_{SM}$  and $F_{SM}$. More precisely, defining
\begin{equation}
\mathfrak c=\max\left\{\mathfrak c_9,\, \mathbf c_9,\, \widetilde c\right\}.
\end{equation}
when $L_0\geq \mathfrak c$ and $k\geq 0$,
\begin{gather}
\nonumber
  \mathbb E\left[\rho_{k+1}^{a/4}\right]\leq c_3' \left\{(2\kappa)^{-10 \nu_1 L_{k+1}}\left(c_7'\widetilde L_{k+1}^{d-2}\frac{L_{k+1}^3}{L_k^2}\widetilde L_k\mathbb E\left[q_k\right]\right)^{\frac{\widetilde L_{k+1}}{16 N_k\widetilde L_k}}\right. \\
\label{critrenor1}
  \left.+\sum_{0\leq m\leq N_k+1}\left(c_8'\widetilde L_{k+1}^{d-1}\mathbb E\left[\rho_k^a\right]\right)^{\frac{N_k+m-1}{4}}\right\}.
\end{gather}
The next lemma provides a recursion to obtain controls of stretched exponential type on certain moments of $\rho_k, \, k\geq 0$. The main assumption will be the seed estimate, as we shall see soon that estimate is what we will call effective criterion. Keeping in mind scales $(L_k)_{k\geq 0}$, $(\widetilde L_k)_{k\geq 0}$ satisfying (\ref{escalesL_k}) we have:

\begin{lemma}
\label{lemmacriterion}
There exists a positive constant $c_9'(d)\geq\max\{\mathfrak c, c_2'\}$ such that whenever $L_0\geq c_9'$ along with $ L_0\leq \widetilde L_0\leq L_0^4$ and for some $a_0\in (0,1]$, $u_0\in [\max\{e^{-\frac{gL_0}{8d-3}},(2\kappa)^{\frac{7L_0}{4d-1}}\},1]$, the inequality
\begin{equation}\label{seedce}
\varphi_0\stackrel{\mbox{def}}=(c_7'\vee c_8')\widetilde L_1^{(d-1)}L_0\mathbb E[\rho_0^{a_0}]\leq(2\kappa)^{u_0L_0}
\end{equation}
holds, then for all $k\geq 0$ one has that
\begin{equation}\label{kestimce}
\varphi_k\stackrel{\mbox{def}}=(c_7'\vee c_8')\widetilde L_{k+1}^{(d-1)}L_k\mathbb E[\rho_k^{a_k}]\leq(2\kappa)^{u_kL_k},
\end{equation}
 where the constants $c_7'$ and $c_8'$ satisfy \eqref{c7c8} (see also \eqref {proprecurce}),
\begin{equation}\label{akuk}
a_k=a_04^{-k},\, u_k=u_0v^{-k}.
\end{equation}
\end{lemma}
\noindent \emph{Proof.}
We argue by induction following a similar procedure as in \cite[Lemma 2.2]{Sz02}. Consider the set 
$$
 E:=\{k\geq 0: \varphi_k>(2\kappa)^{u_k L_k}\}
$$
and notice that the statement of the lemma follows once we have proven that for some constant $c_9'$, $L_0\geq c_9'$ implies $E=\varnothing$. Assume that $E \neq \varnothing$, then we will show that for $L_0$ sufficiently large the hypothesis yields a contradiction. By \eqref{seedce} $0\notin E$ and thus denoting by $k+1$ the minimal natural number in $E$, we apply inequality (\ref{critrenor1}) to get
\begin{equation}\label{criterecurk}
\varphi_{k+1}\leq c_3' (c_7'\vee c_8')\widetilde L_{k+2}^{(d-1)}L_{k+1}\left\{(2\kappa)^{-10\nu_1L_{k+1}}\varphi_k^{\frac{N_k^3}{16}}+\sum_{0\leq m\leq N_k+1}\varphi_k^{\frac{[N_k+m-1]}{4}}\right\}.
\end{equation}
Since $k\notin E$ we have:
$$
(2\kappa)^{-10\nu_1L_{k+1}}\varphi_k^{\frac{N_k^3}{32}}\leq
(2\kappa)^{-10\nu_1L_{k+1}}(2\kappa)^{\frac{u_kL_{k+1}N_k^2}{32}}\stackrel{(\ref{u0valpha})-(\ref{escalesL_k})}\leq1.
$$
Using $\lfloor N_k \rfloor -1\geq N_k/2$ and $k\notin E$ once again, we have
\begin{gather}
\nonumber
\varphi_{k+1}\leq c_3'(c_7'\vee c_8')\widetilde L_{k+2}^{(d-1)}L_{k+1}\left\{\varphi_k^{\frac{N_k^3}{32}}+L_{k+1}\varphi_k^{\frac{N_k}{8}}\right\}\\
\label{almfinishlemre}
\leq 2c_3' (c_7'\vee c_8')\widetilde L_{k+2}^{(d-1)}L_{k+1}^2\varphi_k^{\frac{N_k}{16}}\,\,(2\kappa)^{u_{k+1}L_{k+1}}.
\end{gather}
The rest of the proof consists in finding a constant $c_9'>0$ such that whenever $L_0\geq c_9'$ one has
\begin{equation}\label{claim}
2c_3' (c_7'\vee c_8')\widetilde L_{k+2}^{(d-1)}L_{k+1}^2\varphi_k^{\frac{N_k}{16}}\leq 1  ,
\end{equation}
which produces the contradiction. For this end, we observe that in view of (\ref{escalesL_k}), after performing some basic estimations using $\widetilde L_0\leq L_0^4$, we get
$$
2c_3' (c_7'\vee c_8')\widetilde L_{k+2}^{(d-1)}L_{k+1}^2\varphi_k^{\frac{N_k}{16}}\leq2c_3' (c_7'\vee c_8')L_k^{4d-2}v^{4(d-1)}N_k^{4d-1}(2\kappa)^{20\nu_1 L_k} \, .
$$
Since $\kappa \leq 1/4$, it is clear by inspection at (\ref{expescaleL_k}) that one can find a constant $\widetilde c(d)$ such that whenever $L_0\geq \widetilde c$,
$$
2c_3'(c_7'\vee c_8')L_k^{4d-2}v^{4(d-1)}(2\kappa)^{\nu_1L_k}\leq 1.
$$
As a result,
\begin{gather}
\nonumber
2c_3' (c_7'\vee c_8')\widetilde L_{k+2}^{(d-1)}L_{k+1}^2\varphi_k^{\frac{N_k}{16}}\leq N_k^{4d-1}\kappa^{19\nu_1 L_k}\\
\label{ineceimp}
=\left(\frac{\alpha \nu_1}{u_0}\right)^{4d-1}v^{(4d-1)k}(2\kappa)^{19\nu_1(\frac{\alpha \nu_1}{u_0})^kv^{\frac{k(k-1)}{2}} L_0}.
\end{gather}
In order to finish the proof, we observe that one can choose $\widehat{c}(d)$ such that, when $L_0\geq \widehat c$, the right-most expression in (\ref{ineceimp}) is smaller than $1$ for all $k\geq 1$. On the other hand, notice that for $k=0$, since $u_0\in[(2\kappa)^{\frac{7L_0}{4d-1}},1]$ by choosing a further positive constant $\overline c(d)$, for $L_0\geq \overline c$ one has
\begin{equation}
\left(\frac{\alpha \nu_1}{u_0}\right)^{4d-1}(2\kappa)^{19\nu_1(\frac{\alpha \nu_1}{u_0})^k v^{\frac{k(k-1)}{2}} L_0}\leq u_0^{-4d+1}(2\kappa)^{7L_0}\leq 1.
\end{equation}
As already mentioned, the claim (\ref{claim}) contradicts the assumption that $E \neq \varnothing$ which ends the proof. \hfill $\square$

\medskip
An important point consists in proving that under an appropriate version of a seed condition as in (\ref{seedce}), condition $(T')|\ell$ holds. However we shall first introduce some further notations and remarks. We are looking for a nice expression for $\varphi_0$ of Lemma \ref{lemmacriterion}. Notice that in virtue of (\ref{escalesL_k}), $\varphi_0$ in (\ref{seedce}) equals
$$
(c_7'\vee c_8')\left(\frac{\alpha \nu_1}{u_0}\right)^{4(d-1)}\widetilde L_0^{d-1}L_0 \, \mathbb E[\rho_0^{a_0}].
$$
It will be also convenient to consider the function
$$
\lambda:[\max\{e^{-\frac{gL_0}{8d-3}},(2\kappa)^{\frac{7L_0}{4d-1}}\},1]\rightarrow[0,\infty],\,\, \lambda(u)=u^{4(d-1)}(2\kappa)^{uL_0}
$$ which has its maximum value
$$
\left(\frac{4(d-1)}{e\ln\left(\frac{1}{2\kappa}L_0\right)}\right)^{4(d-1)}
$$
at point $u_0=\frac{4(d-1)}{L_0\ln(1/(2\kappa))}$$\, \in \,$$[\max\{e^{-\frac{gL_0}{8d-3}},(2\kappa)^{\frac{7L_0}{4d-1}}\},1]$, provided that $L_0\geq \mathfrak c_{10}$ for a suitable constant $\mathfrak c_{10}(d,g,C)>0$.

We define the constant $C'':=C''(d, C,g)>0$ by
\begin{equation}
\label{const10}
C''=\max\left\{\mathfrak c_{10}, c_9'\right\}
\end{equation}
as well as the constant $C':=C'(d, C,g)>0$
\begin{equation}
\label{const10-b}
C'=2^{d-1}\left(\frac{e\ln(\frac{1}{2\kappa})}{4(d-1)}\right)^{4(d-1)}(c_7'\vee c_8')\alpha \nu_1 ,
\end{equation}
then we proceed to state the main theorem of this section:
\begin{theorem}\label{theocrite-tprime}
Assume that:
\begin{equation}\label{effectivecriterion}
\inf_{\mathcal B=(R, L-2, L+2,\widetilde L), a\in [0,1]}\left(C'\widetilde L^{d-1}L^{4(d-1)+1}\mathbb E\left[\rho_{\mathcal B}^{a}\right]\right)<1,
\end{equation}
where the infimum runs over all the box specifications $(R, L-2, L+2,\widetilde L)$ where $R$ is a rotation with $R(e_1)=\ell$, $L\geq C''$ and $L \leq \widetilde L< L^4$. Then $(T')|\ell$ is satisfied.
\end{theorem}

\noindent \emph{Proof.}
In virtue of hypothesis (\ref{effectivecriterion}), there exist some $a\in[0,1]$, $L\geq C''$ and  $\widetilde{L}\in[L,L^4)$ such that
\begin{equation}\label{efcri}
C' \widetilde{L}^{d-1}L^{4(d-1)+1}\mathbb E\left[\rho_{\mathcal B}^{a}\right]<1,
\end{equation}
for a box specification $\mathcal B=(R,L-2, L+2, \widetilde L)$. Setting $\widetilde{ L'}=(\widetilde L+1)\wedge L^4 > \widetilde L$, we have
$$
C'2^{-(d-1)}\widetilde{L'}^{d-1}L^{4(d-1)+1}\mathbb E\left[\rho_{\mathcal B}^{a}\right]<1.
$$
We take a rotation $R'$ close to $R$ such that $R'(e_1)=\ell'$ and
\begin{equation}\label{inROT}
p_{(R, L-2, L+2,\widetilde {L})}\leq p_{\mathcal B'}
\end{equation}
for the box specification $\mathcal B'= (R', L-1,L+1, \widetilde {L'})$. Hence using (\ref{inROT}) one gets
$$
C' 2^{-(d-1)}\widetilde{L'}^{d-1}L^{4(d-1)+1}\mathbb E\left[\rho_{\mathcal B'}^{a}\right]<1,
$$
which is in the spirit as the expression for $\varphi_0$ of Lemma \ref{lemmacriterion}, provided that we replace $L_0$ by $L$, $\widetilde L_0$ by $\widetilde{L'}$, $\ell$ by $R'(e_1)$ and set $u_0= 4(d-1)/\ln(L/(2\kappa))$. Therefore letting $L$ play the role of $L_0$, we can apply Lemma \ref{lemmacriterion} under scales given in (\ref{escalesL_k})-(\ref{expescaleL_k}). For this end, we proceed exactly as in the proof of \cite[Proposition 2.3]{Sz02}
to find some suitable constant $c>0$ such that for large $L$
$$
P_0[T_{bL}^{\ell'}>\widetilde{T}_{-\widetilde b L}^{\ell'}]\leq \exp\left(-\widetilde b L\, e^{-c\sqrt{\ln(\widetilde b L)}}\right) \, \ \textrm{for all } b, \ \widetilde b > 0.
$$
Using Definition \ref{deftgammaandtprime}, it is now clear that the required claim follows. \hfill $\square$

\medskip


\smallskip
\section{Polynomial condition for mixing environments and proof of Theorem \ref{mainth1}}\label{secpoly}

Here we prove that a polynomial condition, to be introduced in Definition \ref{defpoly}, implies the effective criterion previously introduced in Section \ref{sectionce}. After we prove this implication, the proof of Theorem \ref{mainth1} basically follows from what has already been established. We prove Theorem \ref{mainth1} at the end of this Section.

Concerning the polynomial condition, as mentioned in the introduction, we follow an analysis similar to the one in \cite{BDR14},  
but we need to make a modification of their Definition 3.6 and a different choice of growth for the scales. These changes are more suitable to obtain good quenched estimates in our case.

\subsection{Polynomial condition and the associated renormalization scheme}\label{secpolqe}

We start with the choice of scales, after which we shall introduce our polynomial condition for mixing environments $(P_J)|\ell$. Setting $v:=44$, for every nonnegative integer $k$ we consider positive numbers $N_k$ satisfying:
\begin{equation}
\label{scalespoly}
N_{k+1}=\left(\left\lfloor \frac{15\nu_1 N_0 \ln(1/(2\kappa))}{2\ln(N_0)}\right\rfloor + 1\right)v^{k+1} N_k, \mbox{ with } N_0\geq 3\sqrt d.
\end{equation}
Further restrictions on $N_0$ will be required later on.

We let $k$ be as above, fix $\ell\in\mathbb S^{d-1}$ and a rotation $R$ of $\mathbb R^{d}$ such that $R(e_1)=\ell$. For $z\in N_k\mathbb Z\times N_k^3\mathbb Z^{d-1}$, we define boxes $\widetilde B_{1,k}(z)$ and $B_{2,k}(z)$, via 
\begin{align}\label{boxespolynomial}
   \widetilde B_{1,k}(z):= &R\left(z+[0,N_k]\times[0,N_k^{3}]^{d-1}\right)\cap \mathbb Z^d  \mbox{  and,}\\
   \nonumber
  B_{2,k}(z) := &R\left(z+\left(-\frac{N_k}{11}, N_k+\frac{N_k}{11}\right)\times\left(-\frac{N_k^3}{10}, N_k^3+\frac{N_k^3}{10}\right)^{d-1}\right)\cap \mathbb Z^{d}.
\end{align}
with their frontal boundary parts defined as in \eqref{fb+}.
For technical reasons, we also introduce a box $\dot B_{1,k}(z)$ which is contained in $\widetilde B_{1,k}(z)$ and defined as follows:
\begin{align}
\nonumber
\dot B_{1,k}(z):=&R\left(z+(0,N_k)\times (0,N_k^3)^{d-1}\right)\cap \mathbb Z^d.
\end{align}
The choice of scales above implies the following property: whenever two boxes $B_{2,k}(z_1)$ and $B_{2,k}(z_2)$ are disjoint, they are indeed separated in $\ell^1$-distance by a factor (independent of $k$) of $N_k$. This requirement is not only needed in the mixing case, and indeed it plays a role even for i.i.d. random environments (see Remark \ref{remarkgoodpoly} for a further explanation). 


For each integer $k\geq0$, let us denote by $\mathfrak L_k$ the subset of $\mathbb Z^d$
\begin{equation}\label{setposboxpoly}
\mathfrak L_k:= N_k\mathbb Z \times N_k^3\mathbb Z^{d-1},
\end{equation}
and by $\mathfrak B_k$ the set of boxes of scale $k$:
\begin{equation}\label{setofboxesk}
  \mathfrak{B}_k:=\{B_{2,k}(z),\, z\in \mathfrak L_k\}.
\end{equation}
The upcoming renormalization procedure requires the following: for each box $B_{2,k}(x)\in \mathfrak B_k$, where $k\geq 1$ and $x\in \mathfrak L_k$, the boxes $\dot B_{1,k-1}(z)$, $z\in \mathfrak L_{k-1}$, such that $\dot B_{1,k-1}(z)\subset B_{2,k}(x)$ form a  \textit{quasi-covering} of $B_{2,k}(x)$, in the following sense:
\begin{equation}\label{unionscalespoly}
B_{2,k}(x)\subset \bigcup_{\substack{z \in \mathfrak L_{k-1}, \\
\dot B_{1,k-1}(z)\subset B_{2,k}(x) }}\widetilde B_{1,k-1}(z).
\end{equation}
This is satisfied if $N_0$ is such that for any $k\geq 1$ we have
\begin{equation}\label{requpolno}
N_{k}/N_{k-1}\in 110\mathbb N.
\end{equation}
In turn, notice that a sufficient condition to satisfy (\ref{requpolno}) is
$$
\left\lfloor \frac{15\nu_1 N_0 \ln(1/(2\kappa))}{2\ln(N_0)}\right\rfloor + 1\in 110 \mathbb N.
$$

It will be convenient to assume that $N_0$ is a fixed number satisfying the previous requirement. We now introduce our \textit{mixing polynomial condition}.

\medskip
\begin{definition}\label{defpoly}
For $J>0$ and $\ell\in \mathbb S^{d-1}$, we say that the mixing polynomial condition $(P_J)|\ell$ holds if for some $N_0$ large enough and satisfying (\ref{requpolno}) we have
\begin{equation}
\sup_{\substack x\in \widetilde B_{1,0}(0)}P_x\left[X_{T_{B_{2,0}(0)}}\notin \partial^+ B_{2,0}(0)\right]<\frac{1}{N_0^{J}}.
\end{equation}
\end{definition}


We need to introduce the notion of \text{good box} suitable into this polynomial framework.

\begin{definition}
\label{defgoodboxpoly}
Consider $z\in \mathfrak L_k$. For $k=0$ we say that the box $B_{2,0}(z)\in \mathfrak B_0$ is $N_0$-\textit{Good} if
\begin{equation}\label{goodboxpoly0}
  \inf_{\substack{x\in \widetilde B_{1,0}(z)}}P_{x,\omega} [X_{T_{B_{2,0}(z)}}\in \partial^+B_{2,0}(z)]>1-\frac{1}{N_0^5}.
\end{equation}
Otherwise we say that $B_{2,0}(z)$ is $N_0$-\textit{Bad}. Inductively, for $k\geq 1$ and $z\in \mathfrak L_k$ we say that the box $B_{2,k}(z)\in\mathfrak B_k$ is $N_k$-\textit{Good} if there exists $t\in \mathfrak L_{k-1}$ with $B_{2,k-1}(t)\in \mathfrak B_{k-1}$, $\dot B_{1,k-1}(t)\subset B_{2,k}(z)$ such that for each $y\in \mathfrak L_{k-1}$, $y\neq t$, such that $B_{2,k-1}(y)\in\mathfrak B_{k-1}$, $\dot B_{1, k-1}(y)\subset B_{2,k}(z)$ and $B_{2, k-1}(y)\cap B_{2,k-1}(t)=\varnothing$, the box $B_{2, k-1}(y)$ is $N_{k-1}$-\textit{Good}.
Otherwise, we say that box $B_{2,k}(z)$ is $N_k$-\textit{Bad}.
\end{definition}

Several remarks are needed for reference purposes later on.

\begin{remark}
\label{remarkgoodpoly}
\leavevmode
\begin{itemize}
\item[(1)]  Informally, for any $k \ge 1$ we would like to say that the box $B_{2,k}(0)$ is $N_k$-Good if there is \textit{at most one bad box} $B_{2,k-1}(t)$, $t\in \mathfrak L_{k-1}$ of scale $k-1$ with $\dot B_{1,k-1} (t) \subset B_{2,k} (0)$. This is not quite so, but the definition says that all the boxes of scale $k-1$ contained in $B_{2,k}(0)$  and disjoint from $B_{2,k-1}(t)$ must be good. That is, all the boxes of scale $k-1$ intersecting $B_{2,k-1}(y)$ might also be bad, but not more. So in total there are at most $3^d$ bad boxes of scale $k-1$ inside a box of scale $k$.
 \item[(2)] Notice that the property of being $N_k$-Good for a box $B_{2,k}(x) \in \mathfrak B_k$, with $x\in \mathfrak L_k, \, k\geq 0,$ 
 depends at most on transitions at sites of the set $\mathcal{B}_{k,x}$ defined by: 
    \begin{align}\label{attachedboxpoly}
    \mathcal B_{k,x}:=&\left\{z\in\mathbb Z^d: z\in R\left(x+\Big(\sum_{i=0}^{k-1}\frac{-N_i}{11},(1+1/11)N_{k}+\sum_{i=0}^{k-1}\frac{N_i}{11}\Big)\times \right.\right.\\
    \nonumber
    & \qquad \qquad \left.\left.\Big(-\sum_{i=0}^{k}\frac{N_i^3}{10}, (1+1/10)N_k^3 + \sum_{i=0}^{k-1}\frac{N_i^3}{10}\Big)^{d-1}\right) \right\},
    \end{align}
    which is straightforward to prove by induction. We further observe that
    $$
    \sum_{i=0}^{k-1}\frac{N_i^3}{10}\leq \frac{N_k^3}{10} \mbox{ and }\sum_{i=0}^{k-1}\frac{N_i}{11}\leq \frac{N_k}{11}.
    $$
\item[(3)] In a given box $B_{2,k}(z)$, with fixed $z\in \mathfrak L_k$ and $k\geq 0$, the number of other boxes in $\mathfrak B_k$ intersecting $B_{2,k}(z)$ along a given fixed direction $\pm R(e_i),\, i\in[1,d]$ is at most two. Furthermore, along a given direction we have that nonconsecutive boxes are disjoint, indeed they are separated by at least a distance $9N_k/11$ in terms of the $\ell_1$-norm. Consequently, disjoint boxes have attached boxes of the type (\ref{attachedboxpoly}) separated by at least $7N_k/11$ in $\ell_1$-norm  and thus the present choice of scales allows the use of the mixing condition on renormalization schemes.
\end{itemize}
\end{remark}

The next objective is to get doubly exponential upper bounds in $k$ for the probability of a given box in $\mathfrak B_k$ to be $N_k$-\textit{Bad}. We begin with case $k=0$.
\begin{lemma}\label{lemmaseedpolybad}
Assume $(P_J)|\ell$ to be fulfilled for given $J>0$ and $\ell \in \mathbb S^{d-1}$. Then for any $z\in\mathfrak L_0$
\begin{equation}\label{estimabad0}
\mathbb P[B_{2,0}(z) \mbox{ is } N_0-\mbox{\textit{Bad}} ]\leq \frac{1}{N_0^{J-3(d+1)}} \, .
\end{equation}
\end{lemma}
\noindent \emph{Proof.}
Let $J$, $\ell$ and $z$ be as in the statement of the lemma. From Definition \ref{defgoodboxpoly} we have
\begin{align*}
\mathbb P[B_{2,0}(z) \mbox{ is } N_0-\mbox{\textit{Bad}}]&= \mathbb P [\sup_{\substack x\in \widetilde B_{1,0}(z)}P_{x,\omega}[X_{T_{B_{2,0}(z)}}\notin \partial^+B_{2,0}(z)]\geq\frac{1}{N_0^5}]\\
&\leq N_0^5|\widetilde B_{1,0}(z)|\sup_{\substack{x\in \widetilde B_{1,0}(z)}}P_x[X_{T_{B_{2,0}(z)}}\notin \partial^+B_{2,0}(z)]\\
&\leq N_0^5|\widetilde B_{1,0}(z)|\frac{1}{N_0^J}\leq\frac{1}{N_0^{J-5-3d+2}} = \frac{1}{N_0^{J-3(d+1)}},
\end{align*}
where we have used Chebyshev's inequality in the first inequality and the hypothesis $(P_J)|\ell$ in the second one.\hfill $\square$

\medskip
We continue with the estimates for general $k\geq1$ in the following:
\begin{proposition}\label{lemmakpolybad}
Assume $(P_J)|\ell$ to be fulfilled for given $J>9d$ and $\ell\in \mathbb S^{d-1}$.  There exist positive constants $\eta_1:=\eta_1(d)$ and $\eta_2:=\eta_2(\eta_1, d)$ such that whenever $N_0\geq \eta_1$ and $N_0$ satisfies (\ref{requpolno}) we have
\begin{equation}\label{estimabadk}
\mathbb P[B_{2,k}(z) \mbox{ is } N_k-\mbox{\textit{Bad}}]\leq \exp\left(-\eta_2 2^{k}\right)
\end{equation}
for all $k\geq 0$ and $z\in \mathfrak L_k$.
\end{proposition}
\noindent \emph{Proof.}
Let $J$, $\ell$ and $z$ be as in the statement of this proposition. Consider a sequence of positive constants $(c_j)_{j\geq0}$ defined by (recall Definition \ref{def:smg} for notation)
\begin{align}
\label{ckpoly}
\qquad c_0=&(J-3(d+1))\ln(N_0), \mbox{ and for } k\geq 0\\
\nonumber
c_{k+1}=&c_{k}-\frac{\ln \left(\frac{12}{10}\left(\left\lfloor \frac{15\nu_1 N_0 \ln(1/(2\kappa))}{2\ln(N_0)}\right\rfloor +1\right)v^{k+1}\right)^{6d-4}}{2^{k+1}}-
\frac{\left(C|B_{2,k}|^2e^{-\frac{gN_k}{4}}\right)}{2^{k}}.
\end{align}

By induction, we shall prove: 
\noindent
\textit{For any} $k\geq 0$, $z\in \mathfrak L_k$ 
\begin{equation}\label{claimbadboxk}
\mathbb P[B_{2,k}(z) \mbox{ is } N_k-\mbox{\textit{Bad}}]\leq \exp\left(-c_k2^{k}\right).
\end{equation}
Later on we shall prove that $\inf_{\substack{k\geq 0}}c_k\geq1$ for some sufficiently large $N_0$, which will finish the proof.

For the case $k=0$ in (\ref{estimabadk}), we simply recall that $J>9d$ and rewrite (\ref{estimabad0}) as
\begin{equation}
\mathbb P[B_{2,0}(z) \mbox{ is } N_0-\mbox{\textit{Bad}} ]\leq \exp\left((-J+3d+3)\ln(N_0)\right).
\end{equation}
Hence, it suffices to prove the induction step $k$ to $k+1$. For this end, we now assume that inequality (\ref{claimbadboxk}) holds for some $k\geq0$, we shall prove that it also holds when $k$ is replaced by $k+1$. Indeed, let $z\in \mathfrak L_{k+1}$ be fixed and consider the environmental event: $z$-Bad:=$\{ B_{2,k+1}(z)$ is $N_{k+1}$-\textit{Bad}$\}$. Using Definition \ref{defgoodboxpoly} one sees that $z$-Bad is a subset of the event: ``there exist two disjoint boxes $B_{2,k}(t_1),B_{2,k}(t_2)$ of scale $k$ such that $\dot B_{1,k}(t_1),\dot B_{1,t_2}(t_2)$ are contained in $B_{2,k+1}(z)$ and $B_{2,k}(t_1),B_{2,k}(t_2)$ are $N_k$-\textit{Bad}". Therefore, introducing for $k\geq0$, the set
\begin{align*}
\Lambda_k:=&\{(t_1, t_2)\in \mathfrak L_k\times \mathfrak L_k:\,B_{2,k}(t_1)\cap B_{2,k}(t_2)=\varnothing,\\
&\dot B_{1,k}(t_1), \dot B_{1,k}(t_2)\subset B_{2,k+1}(z) \},
\end{align*}
we have
\begin{align}\nonumber
\mathbb P[z-\mbox{Bad}]\leq& \mathbb P[\exists (t_1, t_2)\in \Lambda_k: B_{2,k}(t_1), B_{2,k}(t_2) \mbox{ are }N_k-\mbox{\textit{Bad}} ]\\
\label{inebadz}
\leq &\sum_{\substack{(t_1,t_2)\in \Lambda_k}}\mathbb P[B_{2,k}(t_1), B_{2,k}(t_2) \mbox{ are }N_k-\mbox{Bad}].
\end{align}
Observe that by (iii) in Remark \ref{remarkgoodpoly} disjoint boxes of scale $k$ 
are separated in $\ell_1$-norm by at least $7N_k/11$. As a result, using 
\textbf{(SMG)}$_{C,g}$, the last probability inside the sum in (\ref{inebadz}) splits into a product of two factors up to a mixing correction. More precisely, we introduce for $k\geq 0$ and $(t_1,t_2) \in \Lambda_k$ the mixing factor (under notation of Remark \ref{remarkgoodpoly}):
\begin{equation*}
\Gamma_k:=\exp\left(\sum_{\substack{x\in \mathcal B_{k,t_1}\\ y\in \mathcal B_{k,t_2} }}Ce^{-g|x-y|_1}\right).
\end{equation*}
Then under the previous notation and 
\textbf{(SMG)}$_{C,g}$ we have
\begin{align*}
&\mathbb P[B_{2,k}(t_1) \mbox{ and } B_{2,k}(t_2) \mbox{ are }N_k-\mbox{Bad}]\leq\\
&\Gamma_k \times \mathbb P[B_{2,k}(t_1)\mbox{ is }N_k-\mbox{\textit{Bad}}]\times
\mathbb P[B_{2,k}(t_2)\mbox{ is }N_k-\mbox{\textit{Bad}}].
\end{align*}
For easy of notation,
we set
$$
T_{k+1}:=\frac{12}{10}\left(\left\lfloor \frac{15\nu_1 N_0 \ln(1/(2\kappa))}{2\ln(N_0)}\right\rfloor +1\right)v^{k+1} = \frac{12}{10} \frac{N_{k+1}}{N_k}.
$$
Going back to (\ref{inebadz}) together with induction hypothesis (\ref{claimbadboxk}) we have that $\mathbb P[z-\mbox{Bad}]$ is bounded above by
\begin{align}
\nonumber
& | \Lambda_k|\times \Gamma_k\times e^{-c_k2^{k+1}}\\
\nonumber
\leq& \left(T_{k+1}\right)^{6d-4}
\times e^{-\left(c_k-\frac{\ln(\Gamma_k)}{2^{k+1}}\right)2^{k+1}}\\
\label{inepoly1}
\leq& \exp\left(-2^{k+1} \left(c_k-\frac{(6d-4)\ln \left(T_{k+1}\right)}{2^{k+1}}-\frac{\left(C|B_{2,k}|^2e^{-\frac{gN_k}{4}}\right)}{2^{k}}\right)\right).
\end{align}
We have used in the previous chain of inequalities (\ref{inepoly1}) a non-sharp estimate:
$$
\Gamma_k\leq e^{2C|B_{2,k}|^2e^{-\frac{gN_k}{4}}},
$$
together with the construction prescribed by (\ref{scalespoly}). (Recall also the discussion after (\ref{unionscalespoly}) and in (\ref{requpolno}).) We now use the definitions given in (\ref{ckpoly}) to get that
\begin{equation*}
\mathbb P[z-\mbox{Bad}]\leq \exp\left(-c_{k+1}2^{k+1}\right)
\end{equation*}
provided $N_0\geq \eta$, where $\eta$ is a suitable constant.
%
The induction is finished, therefore at this point it suffices to prove the claim: \textit{There exists constant} $\eta'(d)\geq \eta$ \textit{such that the inequality:}
$$
\inf_{\substack{j\geq 0}}c_j\geq 1
$$
\textit{holds, provided that} $N_0\geq \eta'$. However a rough estimate gives that
$$
\lim_{k\rightarrow\infty} |c_k-c_0|
$$
is smaller than
\begin{align*}
&\sum_{\substack{j\geq 0}}\frac{\ln \left(T_{j+1}\right)^{6d-4}}{2^{j+1}}+
\frac{2C|B_{2,j}|^2e^{-\frac{gN_j}{4}}}{2^{j+1}} \\
\leq&(6d-4)\ln\left(\frac{12}{10}\left(\left\lfloor \frac{15\nu_1 N_0 \ln(1/(2\kappa)}{2\ln(N_0)}\right\rfloor +1\right)\right)+(6d-4)2\ln(v)+1.
\end{align*}
It is now straightforward to verify that $\inf_{k\geq0} c_k\geq 1$, provided $N_0$ is large but fixed and using  $c_0=(M-3(d+1))\ln(N_0)$, which ends the proof. \hfill $\square$

\medskip
We now take a further step into the proof of Theorem \ref{mainth1}. Similarly as the argument given in \cite{BDR14},  quenched exponential bounds for the unlikely exit event from a \textit{Good} box will be needed. Notice that we will have a weaker decay in the quenched estimate than the similar result stated in \cite[Proposition 3.9]{BDR14}. This is mostly due to the present choice of scales (\ref{scalespoly}) and definition of good boxes in Definition \ref{defgoodboxpoly}, which were introduced in that form to clearly avoid any intersection problem. 
Roughly speaking, in comparison with the construction in \cite{BDR14},  for a given box of scale $k$ we have less boxes of scale $k-1$ that intersect it under our scaling construction. The formal statement is as follows:

\begin{proposition}\label{propquenpol}
 There exist positive constants $\eta_3:=\eta_3(d)$ and $\eta_4:=\eta_4(\eta_3,d)$ such that whenever $N_0\geq \eta_3$ and $N_0$ satisfies (\ref{requpolno}) we have, for each  $k\geq 0$, each $z\in \mathfrak L_k$ and each $B_{2,k}(z)\in \mathfrak B_k$ assumed $N_k$-\textit{Good}:
\begin{equation}\label{Quenchedexpbo}
\sup_{x\in \widetilde B_{1,k}(z)}P_{x,\omega}\left[X_{T_{B_{2,k}(z)}}\notin \partial^+B_{2,k}(z)\right]\leq e^{-\frac{\eta_4 N_k}{v^{k+1}}}
\end{equation}
\end{proposition}
\noindent \emph{Proof.}
By stationarity under spatial shifts of the probability measure $\mathbb P$, it is enough to prove the proposition for general $k\geq 0$ and $z=0$. Henceforth, we assume that $z=0$ along with $B_{2,k}(0)$ is $N_k$-\textit{Good} for given $k\geq 0$. In order to simplify notation, we shall drop the dependence of $0\in \mathbb R^d$ from the boxes $\widetilde B_{1,k}:=\widetilde B_{1,k}(0)$ and $B_{2,k}:=B_{2,k}(0)$.

We will first prove the following by induction: For each $k\geq 0$ there exists a sequence of positive numbers $(c_k)_{k\geq0}$ such that whenever $N_0\geq \zeta_1$
\begin{equation}\label{quenchestgoodp}
\sup_{x\in \widetilde B_{1,k}}P_{x,\omega}\left[X_{T_{B_{2,k}}}\notin \partial^+B_{2,k}\right]\leq e^{-c_k N_k}
\end{equation}
holds. The sequence $(c_k)_{k\geq0}$ can be derived along the proof of the induction, however we keep the convention as in the previous proposition of defining them here:
\begin{align}\label{ckquepoly}
&c_0:=\frac{5\ln N_0}{N_0}, \mbox{ and }
c_k:=\frac{5\ln N_0}{(44)^k N_0},  \mbox{ for $k\geq1$.}
\end{align}

\noindent
Our method to prove inequality (\ref{quenchestgoodp}) will follow a similar analysis as the one given in Section \ref{sectionce}, whose argument was originally due to Sznitman \cite{Sz02}.

Let us start with $k=0$, in this case we observe that whenever $B_{k,2}$ is $N_k$-\textit{Good} by Definition \ref{defgoodboxpoly} we have
\begin{equation}\label{0induction}
\sup_{\substack{x\in \widetilde B_1}}P_{x}[X_{T_{B_{2,0}}}\notin\partial^+B_{2,0}]\leq \frac{1}{N_0^5}=e^{-5\frac{\ln{N_0}}{N_0}N_0}.
\end{equation}
Hence we have that (\ref{quenchestgoodp}) holds from the very definition of $c_0$ in (\ref{ckquepoly}).

We now proceed to prove the induction step for $k\geq1$. It amounts to prove the statement of (\ref{quenchestgoodp}) for $k$ assuming that (\ref{quenchestgoodp}) is satisfied for $k$ replaced by $k-1$. For this purpose, we introduce for $k\geq 0$, $u\in\mathbb R$ and $i\in [2,d]$ the $\mathcal F_n$-stopping times $\sigma_u^{+i}$ and $\sigma_u^{-i}$ defined by
\begin{align}
\sigma_u^{+i}:=&\inf\{n\geq 0 : (X_n-X_0)\cdot R(e_i)\geq u\} \mbox{ and}\\
\sigma_u^{-i}:=&\inf\{n\geq 0 : (X_n-X_0)\cdot R(e_i)\leq u\},
\end{align}
which are analogous to the random times defined in \eqref{sigmatimes}. Let $i\in [2,d]$ and consider the following notation for the stopping times
\begin{align*}
\vartheta_k^{+i} = \sigma_{\frac{N_k^3}{10}}^{+i} \mbox{  and } \vartheta_k^{-i} = \sigma_{-\frac{N_k^3}{10}}^{-i} \, .
\end{align*}
It will be useful to introduce as well, the path event $\mathfrak I_k$ for $k\geq 1$
\begin{equation}\label{=Jk}
\mathfrak I_k:=\{\exists i\in [2,d]: \vartheta_{k-1}^{+i}\leq T_{B_{2,k}(0)} \ \mbox{ or } \ \vartheta_{k-1}^{-i}\leq T_{B_{2,k}(0)}\}.
\end{equation}
\noindent
Plainly letting $x\in \widetilde B_{1,k}$ be arbitrary, the following decomposition
\begin{align}
\label{decomquenpoly}
P_{x,\omega}\left[X_{T_{B_{2,k}}}\notin \partial^+B_{2,k}\right]\leq &P_{x,\omega}[\mathfrak I_{k}]+ \\
\nonumber
&P_{x,\omega}\left[\mathfrak I_{k}^c, \ X_{T_{B_{2,k}}}\cdot \ell\leq -N_{k}/11 \right]
\end{align}
holds. We now need upper bounds for the terms on the right hand side of (\ref{decomquenpoly}), and we start decomposing the first term:
\begin{equation}\label{ordecpoly}
P_{x,\omega}\left[\mathfrak I_k\right]\leq \sum_{i=2}^{d}\left(P_{x,\omega}\left[\vartheta_{k-1}^{+i}\leq T_{B_{2,k}}\right]+ P_{x,\omega}\left[\vartheta_{k-1}^{-i}\leq T_{B_{2,k}}\right]\right).
\end{equation}
It suffices to prove an upper bound in the spirit of (\ref{quenchestgoodp}) for a single stopping time of the $\vartheta$-type, we say $\vartheta_{k-1}^{+2}$, the others are analogous and the general procedure to be displayed below will provide the same bound for them.

For easy of notation we set
$$
n_{k-1}:=\frac{1}{2}\left(\frac{13N_k}{11N_{k-1}}+1\right)
$$
and note $2n_{k-1}-1$ is the number of consecutive boxes $B_{2,k-1}(t), \, t\in\mathfrak L_{k-1}$, along direction $R(e_1)=\ell$ such that $\dot B_{1,k-1}(t) \subset B_{2,k}$ (recall comments after (\ref{unionscalespoly})). In these terms, it will be convenient to define for each fixed integer $k\geq1$, numbers $I_{k-1}$ and $\overline L_{k-1}$ by
\begin{equation}\label{scalesorthofluc}
I_{k-1}:= \left\lfloor \frac{10N_{k}^{3}}{220N_{k-1}^3n_{k-1}}\right\rfloor \ \ \mbox{ and } \ \
\overline{L}_{k-1}:= \frac{22n_{k-1}N_{k-1}^3}{10}.
\end{equation}

\noindent
Note that for any $x\in \widetilde B_{1,k}$ we have $P_{x,\omega}$-a.s.
$\left\{\vartheta_k^{+2}\geq\sigma_{I_{k-1}\overline L_{k-1}}^{+2}\right\}.$

As in the proof of Proposition \ref{proprecurce}, it will be useful to introduce cylinders $c_\perp(k,j)$ and $\overline c(k,j)$ to determine levels along direction $R(e_2)$ as follows:
\begin{align*}
c_\perp(k,j):=\left\{  z\in \mathbb R^d:\,z\cdot R(e_2)\in \overline L_{k-1}[j,j+1),\ z\in B_{2,k} \right\},
\end{align*}
together with the discrete truncated cylinder $\overline c(k,j)$ defined by
\begin{align*}
\overline c(k,j):=\left\{z\in \mathbb Z^d:\,\inf_{w\in c_\perp(k,j)}|z-w|_\perp\leq 2n_{k-1},\ z\in B_{2,k}  \right\},
\end{align*}
where the \textit{orthogonal norm} $|\,\cdot\,|_\perp$ was defined in \eqref{nperp}. Since $k$ is fixed, we will write simply $c_\perp(j)$ and $\overline c(j)$ unless it seems necessary.

Let us further introduce for technical matters to be clarified soon, the stopping time $\sigma_{u}^{2,x}$ depending on $u\in\mathbb R$ and $x\in \widetilde B_{1,k}$, defined by
$$
\sigma_{u}^{2,x}:=\inf\{n\geq0:\, (X_n-x)R(e_2)\geq u\}.
$$
Using the strong Markov property for arbitrary $x\in \widetilde B_{1,k}$, one gets
\begin{align}
\label{decompsigma}
&P_{x,\omega}\left[\vartheta_{k-1}^{+2}\leq T_{B_{2,k}}\right]\leq P_{x,\omega}\left[\sigma_{I_{k-1}\overline L_{k-1}}^{+2}\leq T_{B_{2,k}}\right]=\\
\nonumber
& E_{x,\omega}\left[\sigma_{(I_{k-1}-2)\overline L_{k-1}}^{+2}<T_{B_{2,k}},P_{X_{\sigma_{(I_{k-1}-2)\overline L_{k-1}}^{+2}},\omega}\left[\sigma_{I_{k-1} \overline L_{k-1}}^{2,x}\leq T_{B_{2,k}}\right]\right].
\end{align}
Plainly we have $P_{x,\omega}$-a.s. on the event $\{\sigma_{(I_{k-1}-2)\overline L_{k-1}}^{+2}<T_{B_{2,k}}\}$, the random variable
$X_{\sigma_{(I_{k-1}-2)\overline L_{k-1}}^{+2}}$
belongs to $c_\perp(I_{k-1}-2)$. The strategy will be to bound from below the probability of the complementary event $\{\sigma_{I_{k-1}\overline L_{k-1}}^{2,x}> T_{B_{2,k}}\}$ for paths starting from arbitrary $y\in c_\perp(I_{k-1}-2)$. Indeed, the strategy consists in pushing the walk successively $(2n_{k-1}-1)$-times to exit boxes of scale $k-1$ inside of $B_{2,k}$ by their boundary side where $\ell$ points out. This will be fulfilled with relatively high probability since $B_{2,k}$ is $N_k-\textit{Good}$. Formally we introduce a sequence of stopping times $(H^i)_{i\geq0}$ along with sequences of successive random positions $(Y_i)_{i\geq 0}$ and $(Z_i)_{i\geq0}$, which are recursively defined as follows:
\begin{align*}
&H^0=0,\,Y_0=X_0,\, Z_0 = z \in \mathfrak L_{k-1} \mbox{ such that }  Y_0 \in \widetilde B_{1,k-1}(z), \\
&H^1=T_{B_{2,k}}\wedge T_{ B_{2,k-1}(Z_0)},\ \mbox{ and for $i>1$,}  \\
&Y_{i-1}=X_{H^{i-1}}, \,\ Z_{i-1} = z\in \mathfrak L_{k-1} \mbox{ such that }  Y_{i-1} \in \widetilde B_{1,k-1}(z),\\
&H^i=H^{i-1}+H^1\circ \theta_{H^{i-1}}.
\end{align*}
Notice that at each step $i$, there might be multiple possibilities for the choice of $Z_i$. We can take any of them.
We also define the stopping time $S$ to indicate a \textit{wrong} exit from box $B_{2,k-1}$ by
\begin{equation*}
S:=\inf\left \{n\geq 0: X_n \in \partial B_{2,k-1}(Z_0)\setminus \partial^+ B_{2,k-1}(Z_0)\right\}.
\end{equation*}
For any $y\in c_\perp(I_{k-1}-2)$, since (\ref{unionscalespoly}) one has $P_{y,\omega}$-a.s. on the event
\begin{equation} \label{HS}
\bigcap_{i=0}^{2n_{k-1}-2}\theta_{H^i}^{-1}\{S> H^1\}
\end{equation}
the walk exits $B_{2,k}$ through $\partial^+ B_{2,k}$ before time $\sigma_{I_{k-1}\overline L_{k-1}}^{2,x}$.
We need to introduce some further definitions in order to get estimates under this strategy. We define for $i\in[0,2n_{k-1}-2]$ and $j\in \mathbb Z$ the random variable
$$
\psi_{i,j-2}:=\inf_{h\in \Theta_{i,j-2}}P_{h,\omega}[S>H_1],
$$
where
\begin{align*}
\Theta_{i,j-2}:=&\{z\in \mathbb Z^d:\, \exists w\in \mathfrak L_{k-1}, z\in \widetilde B_{1,k-1}(w), \\
&w\cdot \ell= iN_{k-1} - \frac{1}{11}N_k,
\dot B_{1,k-1}(w)
\subset \overline c\left(j-2\right)\}.
\end{align*}
From \eqref{HS}, by applying the strong Markov property successively, we get
\begin{align}\label{inequenortho}
P_{y,\omega}[T_{B_{2,k}}<\sigma_{I_{k-1}\overline L_{k-1}}^{2,x}]\geq \prod_{i=0}^{2n_{k-1}-2}\psi_{i,I_{k-1}-2}=:\varphi(I_{k-1}-2).
\end{align}
Hence in virtue of (\ref{decompsigma}), for any $x\in \widetilde B_{1,k}$, we have got the estimate
\begin{align}
\label{ineinothop}
&P_{x,\omega}[\sigma_{I_{k-1}\overline L_{k-1}}^{+2}\leq T_{B_{2,k}}]\leq\\
\nonumber
& P_{x,\omega}\left[\sigma_{(I_{k-1}-2)\overline L_{k-1}}^{+2}\leq T_{B_{2,k}}\right]\left(1-\varphi(I_{k-1}-2)\right).
\end{align}
In turn, let us further decompose the probability in the right hand side of (\ref{ineinothop}), in order to iterate the given process up to this point, to be precise:
\begin{align} \label{ineindorthpoly}
&P_{x,\omega} \left[\sigma_{(I_{k-1}-2)\overline L_{k-1}}^{+2}\leq T_{B_{2,k}}\right] \leq \\
\nonumber
&E_{x,\omega} \left[\sigma_{(I_{k-1}-5)\overline L_{k-1}}^{+2} \leq T_{B_{2,k}} , P_{X_{\sigma_{(I_{k-1}-5)\overline L_{k-1}}^{+2}}} \left[ \sigma_{(I_{k-1}-2) \overline L_{k-1}}^{2,x} \leq T_{B_{2,k}} \right]  \right] .
\end{align}
Now notice that for any $y$ in $c_\perp((I_{k-1}-5)\overline L_{k-1})$, $P_{y,\omega}$-a.s.
\begin{equation}\label{strategyortho}
\bigcap_{i=0}^{2n_{k-1}-2}\theta_{H^i}^{-1}\{S> H^1\}\subset
\{\sigma_{(I_{k-1}-3)\overline L_{k-1}}^{2,x}> T_{B_{2,k}}\}.
\end{equation}
Thus following the same steps leading to \eqref{ineinothop} from \eqref{decompsigma}, we obtain from \eqref{ineindorthpoly} that
\begin{align*}
&P_{x,\omega}\left[\vartheta_{k-1}^{+2}\leq T_{B_{2,k}}\right]\leq\\
\nonumber
& P_{x,\omega}\left[\sigma_{(I_{k-1}-5)\overline L_{k-1}}^{+2}\leq T_{B_{2,k}}\right]\left(1-\varphi(I_{k-1}-5)\right)\left(1-\varphi(I_{k-1}-2)\right).
\end{align*}
It is useful to stress that the lower bound $\left(1-\varphi(I_{k-1}-5)\right)$ in this case, makes use of disjoint boxes of the ones previously used to compute the random variable $\varphi(I_{k-1}-2)$.

By induction for arbitrary $x\in \widetilde B_{1,k}$ one obtains
\begin{equation}\label{inducorthpoly}
P_{x,\omega}\left[\vartheta_{k-1}^{+2}\leq T_{B_{2,k}}\right]\leq \prod_{i=0}^{\lfloor (I_{k-1}-2)/3 \rfloor}\left(1-\varphi(I_{k-1}-3i-2)\right).
\end{equation}
Moreover, 
by assumption $B_{2,k}$ is $N_k$-\textit{Good}, the disjointness of the boxes involved in the events into different product terms in (\ref{inducorthpoly}), and (1) and (3) in Remark \ref{remarkgoodpoly}
we have that there exist at most three $N_{k-1}$-\textit{Bad} boxes in at most one factor $1-\varphi(I_{k-1}-3i-2)$, with $i\in [1,\lfloor (I_{k-1}-2)/3 \rfloor] $.
Thus for any $x\in \widetilde B_{1,k}$ we will bound from above by one the eventual level $i\in [1,\lfloor (I_{k-1}-2)/3 \rfloor]$ containing the bad boxes, recalling (\ref{decompsigma}) and applying the induction hypothesis to each besides one term inside product (\ref{inducorthpoly}), we have
\begin{align} \label{bound+2}
P_{x,\omega}\left[\vartheta_{k-1}^{+2}\leq T_{B_{2,k}}\right] \leq&
\left(1-(1-e^{-c_{k-1}N_{k-1}})^{2n_{k-1}-1}\right)^{\lfloor (I_{k-1}-2)/3 \rfloor-1} \nonumber \\
\leq&\left((2n_{k-1}-1)e^{-c_{k-1}N_{k-1}}\right)^{\lfloor (I_{k-1}-2)/3 \rfloor-1} \nonumber \\
\leq&e^{-c_{k-1}N_{k-1}I_{k-1}/6+\ln(2n_{k-1}-1)I_{k-1}/6}
\end{align}
which is satisfied provided that $N_0\geq \eta$ for certain positive constant $\eta$.

Now recall the definition of $\mathfrak I_k$ from \eqref{=Jk}. As we have pointed out, the same bound as in \eqref{bound+2} holds for each term in (\ref{ordecpoly}).  Hence for every $x\in\widetilde B_{1,k}$
\begin{eqnarray}\label{orthpolyest}
P_{x,\omega}\left[\mathfrak I_k\right]
& \le & (2d-2)e^{-c_{k-1}N_{k-1}I_{k-1}/6+\ln(2n_{k-1}-1)I_{k-1}/6} \nonumber \\
&\le & (2d-2) e^{-c \, c_{k-1} N_{k-1} n_{k-1}^2},
\end{eqnarray}
for some $c>0$.
It remains to treat the second term on the right hand side of (\ref{decomquenpoly}). 
We define for $i\in \mathbb Z$, the strip $\mathcal H_i$ via
$$
\mathcal H_i:=\left\{z\in \mathbb Z^d:\  \exists z'\in \mathbb Z,\  |z-z'|_1=1, (z\cdot \ell- iN_{k-1})\cdot(z'\cdot\ell -i N_{k-1} )\leq 0\right\}.
$$
Recall that we have fixed $k\geq1$ and $x\in\widetilde B_{1,k}$. For later purposes, we further define the set $\mathcal{\widehat H}_{0}$ by
$$
\mathcal{\widehat H}_{0}:=\mathcal H_0\cap \Big\{ z\in\mathbb R^d: \ \forall j\in[2,d], \ |(z-x)\cdot R(e_j)|<\frac{N_k^3}{10} \Big\},
$$
and throughout this part of the proof we let $x_1\in \mathcal{\widehat H}_{0}$ be arbitrary.
Furthermore, we also define the function $I(\cdot)$ on $\mathbb R^{d}$ such that $I(z)=i$ whenever  $z\cdot\ell\in[iN_{k-1}-\frac{N_{k-1}}{2}, iN_{k-1}+\frac{N_{k-1}}{2})$. If $N_0\geq \eta'$ for certain constant, then $I(\mathcal H_i)=i$ for each $i\in \mathbb Z$. We introduce a sequence $(V_j)_{j\geq0}$ the successive times of visits to different strips $\mathcal H_i,\, i\in \mathbb Z$ recursively via:
\begin{align}\label{Vtimes}
  &V_0=0, \, V_1=\inf\{n\geq0:\, X_n\in \mathcal H_{I(X_0)-1}\cup \mathcal H_{I(X_0)+1}\}, \\
  \nonumber
  &\mbox{and for $n\geq 1$  }  V_n=V_{n-1}+V_1\circ \theta_{V_{i-1}}.
\end{align}
Let us also define, for $\omega \in \Omega$ and $y\in \mathbb Z^d$, random variables $q(y,\omega)$, $p(y,\omega)$ and $\rho(y,\omega)$, as follows:
\begin{align*}
q(y,\omega):=&P_{y,\omega}\left[X_{V_1}\in \mathcal H_{I(X_0)-1}\right]=:1-p(y,\omega), \\
\rho(y,\omega):=&\frac{q(y,\omega)}{p(y,\omega)}.
\end{align*}
In these terms for $i\in \mathbb Z$ we further define the random variable $\rho(i):=\rho_\omega (i)$ given by
$$
\rho(i):=\sup_{x\in \mathcal{\widetilde H}_i}\frac{q(x,\omega)}{p(x,\omega)},
$$
where in turn for $i\in \mathbb Z$, the truncated strip $\mathcal{ \widetilde H}_i$ is defined by 
$$
\mathcal{ \widetilde H}_i=\mathcal H_i\cap B_{2,k}.
$$
\noindent
It will be useful to introduce an environmental random function $f_{\omega}:\mathbb Z \rightharpoonup\mathbb R$ defined as
\begin{align*}
f_\omega(j)=0, \mbox{ for $j\geq \frac{12N_{k}}{11N_{k-1}}+1=:m_{k-1}\in \mathbb N$,} \\
f_\omega(j)=\sum_{j\leq m\leq m_{k-1}-1}\prod_{m<i\leq m_{k-1}-1}\rho^{-1}(i).
\end{align*}
For simplicity we drop the variable $\omega$ when there is not risk of confusion. We set
$$
w_{k-1}:=\frac{N_k}{11N_{k-1}}\mbox{($\in \mathbb N$ by construction, see (\ref{requpolno})),}
$$
and then we claim that (recall $x_1\in\mathcal{\widehat H}_{0}$)
\begin{equation}\label{quenchestleftpoly}
P_{x_1,\omega}\left[X_{T_{B_{2,k}}}\cdot \ell\leq \frac{-N_{k}}{11} , (\mathbb Z^{d})^{\mathbb N}\setminus \mathfrak I_k\right]\leq \frac{f(0)}{f(-w_{k-1})}
\end{equation}
holds. Indeed the proof of claim (\ref{quenchestleftpoly}) is similar to the one after (\ref{quenf}), therefore we omit it.
We continue with estimating the rightmost term in (\ref{quenchestleftpoly}) as follows
\begin{equation}\label{inefpoly}
\frac{f(0)}{f(-w_{k-1})}\leq\frac{\sum_{0\leq m\leq m_{k-1}-1}\prod_{\substack{m<j\leq m_{k-1}-1}}\rho^{-1}(j)}{\prod_{\substack{-w_{k-1}<j\leq m_{k-1}-1}}\rho^{-1}(j)}.
\end{equation}

On the other hand, recall that by Remark \ref{remarkgoodpoly} for a $N_k$-\textit{Good} box $B_{2,k}$ the maximum number of consecutive bad boxes of scale $k-1$ \textit{inside} of it and along direction $\ell$, is three.  
We also observe that when $z\in \mathcal{\widetilde{H}}_i$ for some $i\in \mathbb Z$, there exists $y:=y(z)\in \widetilde B_{1,k-1}(t)$ for some $t\in \mathfrak L_{k-1}$ such that $|z-y|_1=1$, together with $y\cdot \ell\geq iN_{k-1}$. As a result of applying uniform ellipticity (\ref{simplex}) to the precedent discussion, we have
$$
\rho(i)\leq \sup_{\substack{z\in \widetilde B_{1,k-1}(t)\\ t\in \mathfrak L_{k-1}\cap \mathcal H_i\\ \dot B_{1,k-1}(t)\subset B_{2,k}}}  \frac{\frac{1}{2\kappa}P_{z,\omega}\left[X_{T_{B_{2,k-1}}}\notin \partial^+ B_{2,k-1}\right]}{1-\frac{1}{2\kappa}P_{z,\omega} \left[X_{T_{B_{2,k-1}}}\notin \partial^+ B_{2,k-1}\right]}.
$$
Using the last estimate, uniform ellipticity again for those eventually bad boxes and applying the induction hypothesis into (\ref{inefpoly}), we get that 
\begin{align}\label{estimateleftpoly}
&P_{x_1 ,\omega}\left[X_{T_{B_{2,k}}}\cdot \ell\leq \frac{-N_{k}}{11} , (\mathbb Z^{d})^{\mathbb N}\setminus \mathfrak I_k\right]\leq\\
\nonumber
&e^{3\nu_1\ln(1/\kappa)N_{k-1}}\sum_{i=0}^{m_{k-1}-1} \left(e^{-c_{k-1} N_{k-1}+\ln(1/(2\kappa))-\ln\left(1-(1/(2\kappa))e^{-c_{k-1}N_{k-1}}\right)}\right)^{w_{k-1}+i-2}.
\end{align}
It was assumed that $x_1\in \mathcal{\widehat H}_{0}$ to obtain (\ref{estimateleftpoly}), however for $z\in \widetilde B_{1,k}$ one has $P_{z,\omega}$-a.s. on the event
$$
\{X_{T_{B_{2,k}}}\cdot \ell\leq \frac{-N_{k}}{11} , (\mathbb Z^{d})^{\mathbb N}\setminus \mathfrak I_k\},
$$
the random time $\widehat H:=\inf\{n\geq 0:\, X_n\in \mathcal{\widehat H}_{0}\}$ is finite and $\widehat H\leq T_{B_{2,k}}$. Thus applying the Markov property
\begin{align*}
&P_{z,\omega}\left[X_{T_{B_{2,k}}}\cdot \ell\leq \frac{-N_{k}}{11} , (\mathbb Z^{d})^{\mathbb N}\setminus \mathfrak I_k\right]\\
\leq&\sum_{x_1\in \mathcal{\widehat H}_{0}}P_{z,\omega}[X_{\widehat H}=x_1]P_{x_1,\omega}[X_{T_{B_{2,k}}}\cdot \ell\leq \frac{-N_{k}}{11} , (\mathbb Z^{d})^{\mathbb N}\setminus \mathfrak I_k]\\
\leq&\sup_{x_1\in \mathcal{\widehat H}_{0}}P_{x_1,\omega}[X_{T_{B_{2,k}}}\cdot \ell\leq \frac{-N_{k}}{11} , (\mathbb Z^{d})^{\mathbb N}\setminus \mathfrak I_k].
\end{align*}
which proves that estimate (\ref{estimateleftpoly}) holds for any $x\in \widetilde B_{1,k}$.

\noindent
In view of (\ref{decomquenpoly}), estimate (\ref{orthpolyest}) and inequality (\ref{estimateleftpoly}) for arbitrary $x\in \widetilde B_{1,k}$ under $B_{2,k}$ is $N_k$-\textit{Good} we have that $P_{x,\omega}[X_{T_{B_{2,k}}}\notin \partial^+B_{2,k}] $  is bounded by
\begin{align}
\label{lastpolyine}
& (2d-2) e^{-c \, c_{k-1} N_{k-1} n_{k-1}^2} + \nonumber \\
& \frac{e^{3\nu_1\ln(1/2\kappa)N_{k-1}}}{1-e^{-c_{k-1} N_{k-1}+\ln(1/(2\kappa)-\ln\left(1-(1/(2\kappa))e^{-c_{k-1}N_{k-1}}\right)}}\\
\nonumber
&\times e^{-c_{k-1} N_{k-1}(m_{k-1}-3)+\ln(1/(2\kappa))(m_{k-1}-3)-\ln(1-(1/(2\kappa))e^{-c_{k-1}N_{k-1}})(m_{k-1}-3)}.
\end{align}
Since $ c c_{k-1} N_{k-1} n_{k-1}^2 > c_k N_k$ for $N_0$ sufficiently large, the first term in  (\ref{lastpolyine}) is bounded from above by $\frac{1}{2}e^{-c_{k}N_k}$. For the second term, note that $m_{k-1}-3 >m_{k-1}/2$ for large but fixed $N_0$, then we get the  upper bound (see Definition (\ref{ckquepoly}))
\begin{align}
\label{finalquenchespoly}
\leq& e^{6\nu_1\ln(1/(2\kappa))N_{k-1}-c_{k-1} N_{k-1}(\frac{1}{2}w_{k-1})+2\ln(1/(2\kappa))(\frac{1}{2}w_{k-1})}\\
\nonumber
\leq & e^{6\nu_1\ln(1/(2\kappa))N_{k-1}-\frac{c_{k-1}N_{k-1} N_k}{33N_{k-1}}} \, \leq \, \frac{1}{2} e^{-c_{k}N_k},
\end{align}
provided that $N_0\geq \zeta$ for some constant $\zeta$. Thus (\ref{finalquenchespoly}) completes the induction and the proof as well. \hfill $\square$

\smallskip
\subsection{Equivalence between the polynomial condition and the effective criterion}\label{secpolthm}

\medskip
Lemma \ref{lemmakpolybad} along with Proposition \ref{propquenpol} enable us to prove a stronger than polynomial decay for the probability of an unlikely exit of the random walk from certain boxes. We introduce the condition representing that stronger decay.

\begin{definition}
\label{deftgamaofn}
Having fixed $\ell\in \mathbb S^{d-1}$ and a rotation $R$ of $\mathbb R^d$ such that $R(e_1)=\ell$, we define the box $\widehat{B}_N$ for $N\geq 3\sqrt d$  by
\begin{equation*}
\widehat B_N:= R\left(\left(-N,N+1 \right)\times\left(-\frac{21N^3}{10}, \frac{21N^3}{10}\right)^{d-1}\right)\cap \mathbb Z^d
\end{equation*}
together with its \textit{frontal boundary part $\partial^+\widehat B_N$ defined as in \eqref{fb+}.} 
Setting $\Gamma: [3\sqrt d , +\infty) \mapsto (0,1)$ such that $\Gamma(N) = 1/(\ln(N))^{\frac{1}{2}}$. We let $N\geq 3\sqrt d$, $\ell\in \mathbb S^{d-1}$ and $R$ be a rotation as above. We say that condition $(\widetilde{T}^{\Gamma})|\ell$ holds, if there exist $N_0\geq 3\sqrt d$ and a constant $c>0$ such that for all $N\geq N_0$ one has
\begin{equation*}
P_0\left[X_{T_{\widehat B_N}}\notin \partial ^+ \widehat B_N\right]\leq e^{-cN^{c \Gamma(N)}}.
\end{equation*}
\end{definition}

\medskip
Throughout the remaining of this section, we let $\ell\in\mathbb S^{d-1}$, $R$ a rotation of $\mathbb R^d$ such that $R(e_1)=\ell$ , $\Gamma$ be the function of Definition \ref{deftgamaofn} and  $J>9d$.
\begin{proposition}\label{proptgampoly}
Condition $(P_J)|\ell$ implies the validity of $(\widetilde{T}^{\Gamma})|\ell$.
\end{proposition}
\noindent \emph{Proof.}
Under condition $(P_J)|\ell$, we consider the construction of scales and boxes of (\ref{scalespoly}). We also use Proposition \ref{lemmakpolybad} together with Proposition \ref{propquenpol} along this proof. Thus we consider for large $N$ the natural $k$ such that $N_k\leq N<N_{k+1}$, which implies that $k$ is of order $\sqrt{\ln(N)}$.

We let $G$ be the environmental event defined by
\begin{equation*}
\{\omega\in \Omega:\, B_{2,k}(z)\in \mathfrak B_k \mbox{ is } N_{k}-\mbox{Good for all } z\in \mathfrak L_k \mbox{ with } B_{2,k}(z) \cap\widehat B_N\neq\varnothing  \}.
\end{equation*}
We next decompose the underlying probability $P_0\left[X_{T_{\widehat B_N}}\notin \partial ^+ \widehat B_N\right]$  as the sum (under notation as in Definition \ref{deftgamaofn}):
\begin{align*}
 &\mathbb E\left[\mathds 1_{G}P_{0, \omega}[X_{T_{\widehat B_N}}\notin \partial ^+ \widehat B_N]\right] 
+\mathbb E\left[\mathds 1_{\Omega\setminus G}P_{0,\omega}[X_{T_{\widehat B_N}}\notin \partial ^+ \widehat B_N]\right].
\end{align*}
Using Proposition \ref{lemmakpolybad} and a rough counting argument, we have
\begin{align}\label{ngoodtgamma}
  \mathbb E\left[\mathds 1_{\Omega\setminus G}P_{0,\omega}[X_{T_{\widehat B_N}}\notin \right.&\left. \partial ^+ \widehat B_N]\right]
  \leq  \mathbb P \left[\Omega\setminus G\right]  \\
  \nonumber
  &
  \leq \left(\frac{13N_{k+1}}{11N_k}+2\right)\times\left(\frac{12N_{k+1}^3}{10N_k^3}+2\right)^{d-1}e^{-\eta_2 2^k}\\
  \nonumber
  &\leq c(d,N_0,\kappa)e^{-\eta_2 2^{k}+(3d-2)(k+1)\ln v}.
\end{align}
Following the proof method developed to prove Proposition \ref{propquenpol} and under notation therein, we define a sequence of stopping times $(H^i)_{i\geq 0}$ as well as two sequences of random positions $(Y_i)_{i\geq0}$ and $(Z_i)_{i\geq0}$ given by:
\begin{align*}
&H^0:=0, \, Y_0=X_0, \, Z_0 = z \in \mathfrak L_{k} \mbox{ such that } Y_0\in \widetilde B_{1,k}(z), \\
&H^1=T_{\widehat B_N}\wedge T_{B_{2,k}(z)}, \mbox{  and for $i>1$}\\
&Y_{i-1}=X_{H^{i-1}},\, Z_{i-1} = z\in\mathfrak L_k \mbox{ such that } Y_{i-1}\in \widetilde B_{1,k}(z), \\
&H^i=H_{i-1}+H^1\circ \theta_{H^{i-1}}.
\end{align*}
We stress again that, in view of (\ref{scalespoly}) and (\ref{unionscalespoly}), the construction of $Z_i$ involves finitely many choices for each $i$.

We also define the stopping time
$$
S:=\inf\left \{n\geq 0: X_n \in \partial B_{2,k-1}(Z_0)\setminus \partial^+ B_{2,k-1}(Z_0)\right\} .
$$
\noindent
By a similar argument as in the proof of Proposition \ref{propquenpol}, $P_{0, \omega}$-a.s.
\begin{equation*}
\bigcap_{\substack{i=0}}^{[N/N_k]+1} \theta_{H^i}^{-1}\{S>H^1\}\subset \left\{X_{T_{B_{2,k+1}(0)}}\in \partial^+B_{2,k+1}(0)\right\}.
\end{equation*}
It is now straightforward to verify that
\begin{align}
\label{goodtgamma}
\mathbb E\left[\mathds 1_{G}P_{0,\omega}[X_{T_{\widehat B_N}}\notin \partial ^+ \widehat B_N]\right]&\leq 1-\left(1-e^{-\frac{\eta_4 N_k}{v^{k}}}\right)^{[N/N_k]+1} \nonumber \\
&\leq (N/N_k+1)e^{-\frac{\eta_4 N_k}{v^{k}}} .
\end{align}
We combine (\ref{goodtgamma}) with (\ref{ngoodtgamma}) to get the existence of a constant $\widetilde c>0 $ such that for large $N$,
\begin{align*}
P_0\left[X_{T_{\widehat B_N}}\notin \partial ^+ \widehat B_N\right]&\leq 2 c(d)e^{-2^{k}+(3d-2)(k+1)\ln v}\leq \exp\left\{-\widetilde c N^{\frac{\widetilde c}{(\ln N)^{1/2}}}\right\} .
\end{align*}
This finishes the proof. \hfill $\square$

\medskip
The next technical tool needed to prove the effective criterion under our polynomial condition, is an estimate for the $\mathbb P$-probability of a quenched atypical event. In spite of the different scales (\ref{scalespoly}) to avoid trouble with intersections in Lemma \ref{lemmakpolybad}, a variant of the argument in \cite[Proposition 2.4] {BDR14} works for the proof.

\begin{proposition}\label{propoweakaqe}
Let $(\widetilde{T}^{\Gamma})|\ell$ be fulfilled. Let, for $N$ Large, $\mathcal U_N\subset \mathbb Z^d$ be the box defined by
$$
\mathcal U_N:=R\left((N-2, N+2)\times (-N^3+1, N^3-1)^{d-1}\right)\cap\mathbb Z^d
$$
and recall the definition of $\partial^+\mathcal U_N$ from \eqref{fb+}. 
Finally let us set $\epsilon:[3\sqrt{d},\infty)\rightharpoonup [0,1]$, defined by $\epsilon(N)=(\ln N)^{-\frac{3}{4}}$. Then, for each function $\beta:[3\sqrt d, \infty)\rightharpoonup[0,1]$ such that $\lim_{M\rightarrow\infty}\epsilon(M)/\beta(M)<1$, we have for large $N$,
\begin{equation}\label{weakaqe}
\mathbb P\left[ P_{0,\omega}\left[X_{T_{\mathcal U_N}}\in \partial^+\mathcal U_N \right]\leq \frac{e^{-\nu_1\ln(1/2\kappa) N^{\epsilon(M) + \beta(N)}}}{2}\right]
\leq \frac{(4\times 6^{d-1})e}{\left[\frac{N^{\beta(N)}}{4\times 6^{d-1}}\right]!}.
\end{equation}
\end{proposition}
\begin{remark}
Observe that $\mathcal U_N$ is within the class of boxes entering in the infimum of (\ref{effectivecriterion}).
\end{remark}
\noindent \emph{Proof.}
Under the notation introduced in statement of this proposition, let us assume that $(\widetilde{T}^\Gamma)|\ell$ holds, and let $\beta$ be any function with the prescribed assumptions. We need to establish a one-step renormalization scheme which guarantees the prescribed decay by (\ref{weakaqe}). To this end, we consider $N_0:=N^{\epsilon(N)}$ for large $N\geq 3\sqrt d$ and, recalling  (\ref{setposboxpoly}),
we introduce, for $z\in \mathfrak L_{0}$, boxes $\widetilde B_1(z)$, $B_2(z)$ 
defined by: 
\begin{align*}
  \widetilde B_1(z) &:=R\left(z+[0,N_0]\times [0, N_0^3]^{d-1} \right)\cap \mathbb Z^d, \\
  B_2(z)&:= R\left(z+(-N_0,N_0+1)\times (-\frac{21N_0}{10}, \frac{31N_0^3}{10}) \right)\cap \mathbb Z^d \mbox{ and }
\end{align*}
with the frontal boundary part of $B_2(z)$ defined as in \eqref{fb+}.
It is then convenient to define \textit{Good/Bad boxes}: for $z\in \mathfrak L_{0}$, $B_2(z)$ is $N_0$-\textit{Good} if
\begin{equation*}
  \inf_{\substack{y\in\widetilde B_1(z)}}P_{y,\omega}\left[X_{T_{B_2(z)}}\in \partial^+ B_2(z)\right]>1-\frac{1}{N^{\epsilon^{-1}}}.
\end{equation*}
\noindent
Otherwise we say that $B_2(z)$ is $N_0$-\textit{Bad}.

\vspace{1ex}
\noindent
Observe that for arbitrary $z\in \mathfrak L_{0}$, using condition $(\widetilde T^\Gamma)|\ell$, we have that
\begin{align}
\nonumber
\mathbb P\left[B_{2}(z)\mbox{ is }N_0-\mbox{\textit{Bad}}\right]&=\mathbb P\left[\sup_{\substack{y\in \widetilde
B_1(z)}}P_{y,\omega}[X_{T_{B_2(z)}}\notin \partial^+B_2(z)]\geq \frac{1}{N^{\epsilon^{-1}}}\right]\\
\nonumber
&\leq N^{\epsilon^{-1}} |\widetilde B_1(z)| e^{-\widetilde c N_{0}^{\frac{\widetilde c}{(\ln N_0)^{1/2}}}}\\
\nonumber
&\leq e^{-\widetilde c N^{\frac{\widetilde c(N)}{(\ln N)^{7/8}}}+(\epsilon^{-1}(N)+3d\epsilon(N))\ln(N)}\\
\label{estimabadwaqe}
&\leq e^{-\widehat c N^{\frac{\widehat c}{(\ln N)^{7/8}}}} ,
\end{align}
for certain positive dimensional constants $\widetilde c$ and $\widehat c$.

\noindent
We need to apply the well-known strategy to exit from box $\mathcal U_N$, starting at $0\in \mathbb R^d$. Further definitions will be required so as to describe the environmental event where the strategy fails. We define the set of boxes involved in the forthcoming strategy $\mathcal G_{N}$ by
\begin{align*}
\mathcal G_N:=&\{B_{2}(z),\, z\in \mathfrak L_{0},\, z\cdot e_1=j N_0,\, j\in[0,[N/N_0]+1],\,\mbox{and for $i\in[2,d]$ }\\
 &z\cdot e_i=jN_0^3, j\in [-[N/N_0]-1, ([N/N_0]+1)]\}.
\end{align*}
Define as well, the environmental event
$$
\mathfrak G_{\beta(N)}:=\left\{\sum_{B_{2}(z)\in \mathcal G_N, z\in \mathfrak L_{0}}\mathds{1}_{\{ B_2(z) \mbox{ is } N_0-\mbox{\textit{Bad} } \}}\leq N^{\beta(N)}\right\}.
$$
As in Proposition \ref{propquenpol} and Proposition \ref{proptgampoly}, we introduce the strategy of successive exits from boxes of type $B_2(z),$  $z\in\mathfrak L_{0}$ through their frontal boundary part. Afterwards we shall see that on event $\mathfrak G_{\beta(N)}$, the strategy shall imply the complementary event involved in the probability (\ref{weakaqe}), and then prove $\mathfrak G_{\beta(N)}$ has high $\mathbb P$-probability. Despite the heavy notation, and for sake of completeness, we introduce again sequences of random positions $(Y_i)_{i\geq0}$, $(Z_i)_{i\geq0}$ together with stopping times $(H^i)_{i\geq 0}$ and $S$ defined as follows:
\begin{align}
\nonumber
&H^0=0,\,Y_0=X_0, \, Z_0 = z \in \mathfrak L_{0} \mbox{ such that } Y_0 \in\widetilde B_1(z) ,\\
\nonumber
&H^1=T_{\mathcal U_N}\wedge T_{B_2(Z_0)}, \,\mbox{  and for $i>1$} \\
\nonumber
&Y_{i-1}=X_{H^{i-1}},\, Z_{i-1} = z \in \mathfrak L_{0} \mbox{ such that } Y_{i-1}\in \widetilde B_1(z), \\
\nonumber
&H^i=H^{i-1}+H^1\circ \theta_{H^{i-1}}, \\
\nonumber
&S=\inf\{n\geq 0: \, X_n\in \partial B_2(Z_0)\setminus \partial^+B_2(z)\}.
\end{align}
It is routine, using uniform ellipticity (\ref{simplex}) to see that for large $N$, on the event $\mathfrak G_{\beta(N)}$ one has $\mathbb P$-a.s.
\begin{align}
\nonumber
P_{0,\omega}\left[X_{T_{\mathcal U_N}}\in \partial^+\mathcal U_N \right]&>P_{0,\omega}\left[\bigcap_{0\leq i\leq \frac{N}{N_0}}\theta_{H^i}^{-1}\left\{S>H^1\right\}\right]\\
\nonumber
&\geq (2\kappa)^{\nu_1 N^{\epsilon(N)+\beta(N)}}\left(1-\frac{1}{N^{\epsilon^{-1}(N)}}\right)^{[N/N_0]+1}\\
\label{goodestpoly}
&\geq \frac{1}{2}e^{\nu_1 \ln(2\kappa)N^{\epsilon(N) + \beta (N)}}.
\end{align}

We now decompose the set $\mathcal G_N$ into smaller sets which shall have boxes elements with appropriate disjointness, so as to apply mixing conditions. Keeping this in mind, we introduce for integer $i\in [0,3]$ subsets $\mathcal G_{N,i, e_1}$ of $\mathcal G_N$ defined by
$$
\mathcal G_{N,i, e_1}:=\{B_{2}(z)\in \mathcal G_N\colon \, z\in \mathfrak L_{0},\, \frac{z\cdot e_1}{N_0} =i  \mod 4\}.
$$
Next, we further decompose each one of the four $\mathcal G_{N,i, e_1}$, $i\in [0,3]$ along the \textit{orthogonal space to} $\ell$. For this end, define for $j_2,\ldots, j_{d}\in [0,5]$ subsets $\mathcal G_{N, i,j_2,\ldots,j_{d}}$ of $\mathcal G_{N, i,e_1}$ given by:
$$
\mathcal G_{N, i,j_2,\ldots,j_{d}}:=\{B_{2}(z)\in \mathcal G_{N, i,e_1}\colon \, z\in \mathfrak L_{0},\, \frac{z\cdot e_k}{N_0^3}=j_k\mod 6\, \forall k\in[2,d]  \}.
$$
Thus, we have constructed $4\times 6^{d-1}$ subsets $\mathcal G_{N, i,j_1,\ldots,j_{d-1}}$ of $\mathcal G_N$, where $i\in [0,3]$ and $j_2,\ldots, j_d\in [0,5]$ such that each pair of boxes $B_2(z_1), B_2(z_2)\in \mathcal G_{N, i,j_2,\ldots,j_{d}}$ with $z_1,z_2\in\mathfrak L_{0} $ are at least $2N_0=2N^{\epsilon(N)}$ separated in terms of $\ell^1$-distance. We plainly have as well,
$$
\bigcup_{\substack{i,j_2,\ldots,j_d\in[0,3]\times[0,5]^{d-1}}} \mathcal G_{N, i,j_2,\ldots,j_{d}}=\mathcal G_N.
$$
Under these terms, it is clear that $\mathbb P$- a.s. the environmental event $\Omega \setminus \mathfrak G_{\beta(N)}$ is contained in the set:
\begin{align}\label{Gcom}
& \bigcup_{\substack{i\in[0,3],\\ j_2,\ldots,j_d[0,5]^{d-1}}} \left\{\sum_{\substack{B_2(z)\in  \mathcal G_{N, i,j_2,\ldots,j_{d}}\\z\in \mathfrak L_{0} }}\mathds 1_{\{B_2(z) \mbox{ is }N_0-\mbox{\textit{Bad}} \}}>\frac{N^{\beta(N)}}{4\times 6^{d-1}}\right\}.
\end{align}

On the other hand, using the remark above, for integers $i, j_2,\ldots, j_d$ and any integer
$$
k\in \left[0, \frac{([N/N_0]+1)\times (2([N/N_0]+1))^{d-1}}{4\times 6^{d-1}}\right]=:[0, N']
$$
using 
\textbf{(SMG)}$_{C,g}$ and a counting argument we have for large $N$,
\begin{align*}
&\mathbb P\left[B_2(z_1), B_2(z_2),\ldots,B_2(z_k)\in \mathcal G_{N, i,j_2,\ldots,j_{d}} \mbox{ are (and not more) } N_0-\mbox{\textit{Bad}} \right]\stackrel{(\ref{estimabadwaqe})}\leq\\
&\qquad \qquad \qquad \exp\left( C(d) k^2 N_0^{3d-2} e^{-2gN_0}\right)\exp\left(- \widehat c \, k  \, e^{\widehat c(\ln N)^{1/8}}\right).
\end{align*}
Hence, for large $N$ there exists constant $\breve{c} = \breve{c}(d) >0$ such that
\begin{align*}
&\mathbb P\left[\sum_{\substack{B_2(z)\in \mathcal G_{N, i,j_2,\ldots,j_{d}}\\z\in \mathfrak L_{0}}}\mathds 1_{\{B_2(z)\mbox{ is }N_0-\mbox{\textit{Bad}}\}}=k\right] \leq \binom{N'}{k}\exp\left(-\breve c \, k\, e^{ \breve c (\ln N)^{1/8}}\right).
\end{align*}
Combining this last estimate with (\ref{Gcom}) and some basic estimates one sees
\begin{align}
\nonumber
&\mathbb P[\Omega \setminus \mathfrak G_{\beta(N)}]\\
\nonumber
&\leq \sum_{\substack{i,j_2,\ldots,j_d\in[0,3]\times[0,5]^{d-1}}}\mathbb P\left[\mathds 1_{\{B_{2}(z)\mathcal G_{N, i, i,j_2,\ldots,j_d}, B_2(z)\mbox{ is }N_0-\mbox{\textit{Bad}}\}}>\frac{N^{\beta(N)}}{4\times 6^{d-1}}\right]\\
\nonumber
&\leq 4\times 6^{d-1}\sum_{\frac{N^{\beta(N)}}{4\times 6^{d-1}} \le j \le N'}\binom{N'}{j}\exp\left(-\breve c \, j\, e^{ \breve c (\ln N)^{1/8}}\right)\\
\label{lastestquenchpoly}
&\leq 4\times 6^{d-1}\sum_{\frac{N^{\beta(N)}}{4\times 6^{d-1}} \le j \le N'}\binom{N'}{j}\frac{1}{(N')^{j}}\leq 4\times 6^{d-1}\frac{e}{\left[\frac{N^{\beta(N)}}{4\times 6^{d-1}}\right]!}.
\end{align}
With (\ref{goodestpoly}) and (\ref{lastestquenchpoly}) we conclude the proof. \hfill $\square$

\medskip
The last required step in order to get the effective criterion will be to integrate the random variable $\rho$ of (\ref{enrique}). The formal statement comes in the next:

\begin{theorem} \label{PgoesEC}
Assume $(P_J)|\ell$ be fulfilled for $J>9d$ and $\ell\in \mathbb S^{d-1}$ (for some $N_0$ large enough but fixed, see Definition \ref{defpoly}). Then $(EC)|\ell$ holds (see Definition \ref{defpolasympandec}).
\end{theorem}
\noindent \emph{Proof.}
Assume condition $(P_J)|\ell$ be fulfilled for some $N_0$ large enough. Consider $N \ge 1$ and recall the notation introduced for a box specification $\mathcal B:= (R,N-2, N+2, 4N^3 )$ in Section \ref{sectionce}, where the rotation $R$ satisfies $R(e_1)=\ell$ and we set $B$ for the box attached to $\mathcal B$. In virtue of Proposition \ref{proptgampoly} there exists constant $\widetilde c>0$  such that for large $N$ (see Definition \ref{defpolasympandec} for notation)
\begin{equation}\label{est1poly}
  \mathbb E\left[q_{\mathcal B}(\omega)\right]\leq e^{-\widetilde cN^{\frac{\widetilde c}{(\ln N)^{1/2}}}}.
\end{equation}
In order to apply inequality (\ref{est1poly}) and Proposition \ref{propoweakaqe}, we define parameters:
\begin{equation}
\label{parampolyfinal}
\beta_1:=\frac{\widetilde c}{2(\ln N)^{1/2}},
\quad
\alpha:=\frac{\widetilde c}{3(\ln N)^{1/2}}
\quad \mbox{and } \ \
a:=\frac{1}{N^{\alpha}}.
\end{equation}
Following the proof argument of \cite[Section 2.2]{BDR14} we split the expectation $\mathbb E\left[\rho_{\mathcal B}^{a}\right]$ into $L$ terms, where
$$
L:=\left[\frac{2(1-\beta_1)}{\beta_1}\right]+1.
$$
Indeed, the decomposition splits terms according the underlying denominator of the random variable $\rho_{\mathcal B}$. The first term is given by
\begin{align*}
\mathcal E_0&:=\mathbb E\left[\rho_{\mathcal B}^a, P_{0,\omega}\left[X_{T_B}\in \partial^+B\right]>\frac{e^{-2\nu_1\ln(1/(2\kappa))N^{\beta_1}}}{2}\right]. \\
\end{align*}
For $j\in\{1,\ldots, L-1\}$, terms $\mathcal E_j$ are defined by
\begin{align*}
\mathcal E_j&:=\mathbb E\left[\rho_{\mathcal B}^a,\frac{e^{-2\nu_1\ln(1/(2\kappa))N^{\beta_{j+1}}}}{2} <P_{0,\omega}\left[X_{T_B}\in \partial^+B\right]\leq\frac{e^{-2\nu_1\ln(1/(2\kappa))N^{\beta_j}}}{2}\right], \\
\end{align*}
and the last term is,
\begin{align*}
\mathcal E_L&:=\mathbb E\left[\rho_{\mathcal B}^a, P_{0,\omega}\left[X_{T_{B}}\in \partial^+B\right]\leq \frac{e^{-2\nu_1\ln(1/(2\kappa))N^{\beta_L}}}{2}\right].
\end{align*}
Where the  numbers $\beta_j$ with $j\in \{1,\ldots,L\}$ involved in the definition above, are prescribed by
$$
\beta_j:=\beta_1+(j-1)\frac{\beta_1}{2}.
$$

Notice that using uniform ellipticity (\ref{simplex}) one sees that $\mathcal E_L=0$ since $\beta_L\geq 1$ and $\mathbb P$-a.s.
$$
P_{0,\omega}\left[X_{T_B}\in \partial^+B\right]>e^{-2\nu_1\ln(1/(2\kappa))N}.
$$
On the other hand, by an application of Jensen's inequality and (\ref{est1poly}), we have
\begin{align}
\nonumber
\mathcal E_0 &\leq 2^{a}e^{2\nu_1\ln(1/(2\kappa))N^{\beta_1-\alpha}}e^{-\widetilde c N^{-\alpha+\widetilde c(\ln N)^{-1/2}}}\\
\label{estterjply}
&\leq 2^{a}e^{2\nu_1\ln(1/(2\kappa))N^{\frac{\widetilde c}{6(\ln N)^{1/2}}}}e^{-\widetilde cN^{\frac{2\widetilde c}{3(\ln N)^{1/2}}}}.
\end{align}
Furthermore, we use the atypical quenched estimate provided by Proposition \ref{propoweakaqe} to get that for $j\in\{1,\ldots,L-1\}$ there exist positive constants $\widehat c, \breve c$ such that
\begin{align}
\nonumber
\mathcal E_j\leq&2\exp(N^{\beta_{j+1}-\alpha})\mathbb P\left[P_{0,\omega}[X_{T_B}\in\partial^+B]\leq \frac{e^{-2\nu_1\ln(1/(2\kappa))N^{\beta_j}}}{2}\right]\\
\nonumber
\leq&\widehat c\,  2\exp\left(2\nu_1\, \ln(1/(2\kappa))N^{\beta_{j+1}-\alpha}\right)\exp\left(-\breve c N^{\beta_j - \epsilon(N)}\ln(N^{\beta_j - \epsilon(N)})\right)\\
\label{esttermj2poly}
\leq&\widehat c\,  2\exp\left(2\nu_1\, \ln(1/(2\kappa))N^{\beta_j+(\beta_1/2)-2\beta_1/3}
-\breve c N^{\beta_j - \epsilon(N)}\ln(N^{\beta_j - \epsilon(N)})
\right)
\end{align}
Considering the orders of $\beta_1$ and $\epsilon(N)$ and combining inequalities (\ref{estterjply}) and (\ref{esttermj2poly}) we see that $\mathbb E[\rho_{\mathcal B}^a]$ is less than any polynomial function in $N$, therefore the proof is complete. \hfill $\square$

\medskip
It is now straightforward to prove Theorem \ref{mainth1}.

\medskip

\noindent \emph{Proof.} [Proof of Theorem \ref{mainth1}]
The implication $(P^*_J)|\ell \,\Rightarrow \, (P_J)|\ell$ follows by a geometric fact. We give the proof for sake of completeness (see \cite{Li16} for an alternative proof). Let boxes $\widetilde B_{1,0}$ and $B_{2,0}$ be defined as in (\ref{boxespolynomial}) for large $N_0$. Consider Definition \ref{defpolasympandec} for the polynomial \textit{asymptotic} condition; then setting $b=1/13$ we can find a neighborhood $\mathcal U_\ell\subset \mathbb  S^{d-1}$ such that $\ell\in \mathcal U_\ell$ and
$$
\lim_{N\rightarrow\infty}N^J P_0\left[\widetilde T_{-bN}^{\ell'}<T_N^{\ell'}\right]=0
$$
for all $\ell'\in \mathcal U_\ell$. Thus in particular taking $\epsilon=\frac{1}{2(d-1)}$ there exist large $L_0$ together with a neighborhood $\mathcal U_\ell\subset \mathbb S^{d-1}$ of $\ell$ such that
\begin{equation}\label{polydecayasymp}
P_0\left[\widetilde T_{-bL_0}^{\ell'}<T_{L_0}^{\ell'}\right]<\frac{1}{2(d-1)L_0^J},
\end{equation}
for all $\ell'\in \mathcal U_\ell$. We take
$$
N_0=\frac{11L_0}{12.5},
$$
and we keep in mind the box involved in Definition \ref{defpoly}. Since $\mathcal U_\ell$ is open containing $\ell$ there exists strictly positive number $\alpha$ such that the following requirements are satisfied:

\begin{itemize}
\item the vectors $\ell_i^{\pm}:=\ell\pm \alpha R(e_i)/|\ell\pm \alpha R(e_i)|_2 \,\in \mathcal U_\ell$ for each $i\in [2,d]$. We stress that there are $2(d-1)$ vectors $\ell_i^{\pm}$, i.e. the symbol $\pm$ is an abbreviation.

\item $\alpha$ is small enough such that
$$
\cos(\arctan(\alpha))\leq \min\left\{\frac{12N_0}{11L_0},\frac{11L_0}{13N_0}\right\} = \frac{12}{12.5} .
$$
\end{itemize}
Consider now the set $\mathcal D\subset\mathbb Z^d$ defined by:
$$
\mathcal D:=\Big\{x\in \mathbb Z^d:\, x\cdot \ell\in \Big(-\frac{N_0}{11}, \frac{12N_0}{11}\Big), \,\forall i\in[2,d], \ x\cdot \ell_i^{\pm}>-\frac{L_0}{13} \Big\}.
$$
Observe that for any $x\in \widetilde B_{1,0}(0)$ we have
$x+\mathcal D  \subset B_{2,0}(0).$
We then introduce the frontal part of the boundary for set $\mathcal D$, denoted by $\partial^+\mathcal D$ and defined as follows:
$$
\partial^+\mathcal D:= \Big\{z\in\mathbb Z^d:\, z\in \partial \mathcal D,\,z\cdot \forall i\in[2,d], \ x\cdot \ell_i^{\pm}>-\frac{L_0}{13} \Big\},
$$
and therefore we note that for any $x\in\widetilde B_{1,0}(0)$ we have $P_x$-a.s. one has
$$
T_{ X_0+\partial\mathcal D\setminus\partial^+\mathcal D}\leq T_{B_{2,0}(0)}.
$$
By stationarity of the probability space $(\Omega, \mathfrak F_{\Omega}, \mathbb P)$ and using (\ref{polydecayasymp}), we get for arbitrary $x\in\widetilde B_{1,0}(0)$
\begin{align*}
P_x[X_{T_{B_{2,0}(0)}}\notin \partial^+B_{2,0}(0)]&\leq P_{x}[X_{H_{x+\partial \mathcal D}}\notin x+\partial^+\mathcal D]\\
&\leq P_0[X_{H_{\partial \mathcal D}}\notin \partial^+\mathcal D]\\
&\leq \sum_{i=2}^{d} \Big( P_0[\widetilde T_{-bL_0}^{\ell_i^+}<T_{L_0}^{\ell_i^+}]+P_0[\widetilde T_{-bL_0}^{\ell_i^-}<T_{L_0}^{\ell_i^-}] \Big)\\
&\stackrel{(\ref{polydecayasymp})}<\frac{1}{N_0^{J}},
\end{align*}
which ends the proof.

\noindent
On the other hand, $(P_J)|\ell$ implies $(EC)|\ell$ by Theorem \ref{PgoesEC} and then it also implies $(T')|\ell$ by Theorem \ref{theocrite-tprime}. Finally, it is clear using any of the equivalent definitions contained in statement of Lemma \ref{lemmaTgamma}, the implication $(T^\gamma)|\ell\, \Rightarrow\,(P_J)|\ell$ holds for any $\gamma>0$. \hfill $\square$

\medskip
\section{Estimates for tails of regeneration times: Proof of Theorem \ref{mainth2}}
\label{secproofmainth}

In this section we prove Theorem \ref{mainth2}. Due to Theorem \ref{mainth1}, we may assume condition $(T')|_\ell$ holds. Under this condition we obtain an upper bound for the tail of the approximate regeneration time $\tau_1$, which implies that $\tau_1$ has finite moments of all orders. We end the proof by applying the annealed central limit theorem of Comets and Zeitouni \cite{CZ02}.

Henceforward, we will assume $(T')|\ell$ with $h\ell=l\in \mathbb Z^d$ for some $h>0$ (see (\ref{rul})), 
\textbf{(SMG)}$_{C,g}$ (see Definition \ref{def:smg}), where $C>0$, $g>2\ln(1/\kappa)$ and $\kappa \in (0,1/(4d))$ is the uniform elliptic constant in (\ref{simplex}). Furthermore, we introduce for $L\in |l|_1\mathbb N$ the approximate regeneration time $\tau_1^{(L)}$ as in Section \ref{section2}, as well as the approximate asymptotic direction  
\begin{equation}
\label{veloc}
\hat v_L:= \frac{E_0 \big[X_{\tau_1^{(L)}} \big| D'=\infty \big]}
{\big| E_0 \big[X_{\tau_1^{(L)}} \big| D'=\infty \big] \big|_2}.
\end{equation}

It is convenient to consider the orthogonal projector $\Pi_{\hat v_L}: \,\mathbb Z^d\,\rightarrow\,\mathbb Z^d$, which projects vectors onto the subspace orthogonal to direction $\hat v_L$:
\begin{equation}\label{projvl}
\Pi_{\hat v_L}(z)=z-(z\cdot v_L)v_L.
\end{equation}
We introduce, for $M>0$, the time $\mathcal L_M$ of the last visit to the half space $\mathcal H_M=\{z\in \mathbb Z^d: z\cdot l\leq M\}$, i.e.
\begin{equation}
\label{lastvistM}
\mathcal L_M := \sup\{n\geq0: X_n\cdot l\leq M \}.
\end{equation}
Furthermore, we let
\begin{equation}
t=\frac{1}{2} \Big(\frac{2\ln\left(\frac{1}{\kappa}\right)}{g}+1 \Big).
\label{defintransfluct}
\end{equation}
The next theorem is our main result in this section. It establishes a faster than polynomial decay for the tail of $\tau^{(L)}_1$.

\begin{theorem}\label{thboundtails}
Let
$$
\sigma=\min\left\{\frac{gt-2\ln(1/\kappa)}{3gt-2\ln(1/\kappa)},\, \frac{(d-1)gt+2\ln(1/\kappa)}{(3d+1)gt-2\ln(1/\kappa)}\right\},
$$
and notice that $\sigma>0$. There exist positive constants $\hat c_{1}$, $\hat c_{2}$ and $L_0$ such that for each $L\geq L_0$ with $L\in |l|_1\,\mathbb N$ and $\alpha\in(1,1+\sigma)$,
\begin{equation}\label{tailtau1}
\overline P_0\left[\tau_1^{(L)}>u\right]\leq e^{- \hat c_{1}\kappa^L\left(\ln(u)\right)^\alpha}+e^{- \hat c_{2}\left(\ln(u)\right)^\alpha}.
\end{equation}
\end{theorem}

\medskip
Using Theorem \ref{thboundtails} and standard analysis (see \cite[Chapter 8, Theorem 8.16]{Ru87}), we obtain that for each $L\geq L_0$ with $L\in |l|_1\mathbb N$ there exists a constant $C=C(L)$ such that 
\begin{equation}\label{tdmomenttau1}
\mathbb P\left[\frac{\overline E_0\left[\left(\kappa^L\tau_1(L)\right)^3, \,D'=\infty|\,\mathfrak F_{0,L}\right]}{P_0\left[D'=\infty|\,\mathfrak F_{0,L}\right]}>C\right]=0,
\end{equation}
where
\begin{equation*}\label{sigmafrak}
 \mathfrak{F}_{x,L}:=\sigma(\omega(y,\cdot): (y-x)\cdot l\leq -(L|l|_2)/|l|_1 ).
\end{equation*}
Theorem \ref{mainth2} follows from equation (\ref{tdmomenttau1}) as in \cite[Theorem 2]{CZ02}.
\medskip


Aiming at the proof of Theorem \ref{thboundtails}, we now state and prove a few technical estimates. The next proposition will be fundamental to apply renormalization schemes. 

\begin{proposition}\label{protransfluctuation}
Let $\eta>0$  and
\begin{equation}
\varsigma \in \left(\frac{gt+2\ln\left(\frac{1}{\kappa}\right)}{2gt},1\right),
\end{equation}
where $t$ was defined in \eqref{defintransfluct}. (Note that the interval above is nonempty under assumption $g>2\ln(1/\kappa)$.) For large $M>0$, we let
\begin{equation}
\label{LM}
L = L(M) = \min \left\{\overline L\in |l|_1\mathbb N: \, e^{-gt\overline L}\leq M^{2\varsigma -2-\epsilon}\right\},
\end{equation} 
where
\begin{equation}\label{epsilonchoise}
\epsilon:=\frac{(2\ln(1/\kappa)(2-2\varsigma))}{(gt-2\ln(1/\kappa))}>0.
\end{equation}
\vspace{2ex}
\noindent
Assuming that $(T')|_\ell$ holds, there exists $c_8>0$ such that for large $M$ 
\begin{equation}
\label{claimproptranv}
P_0\left[\sup_{0\leq n\leq \mathcal L_M}\,|\Pi_{\hat{v}_L}(X_n)|_2>\eta \,M^\varsigma\right]\leq
\exp\left(-c_{8}M^{(2\varsigma-1)-\epsilon}\right)
\end{equation}
 with $\hat{v}_L$ defined in \eqref{veloc}.
\end{proposition}

\begin{remark}
Notice that $(2\varsigma-1)-\epsilon>0$. Taking $g\rightarrow\infty$, one can see that we recover Sznitman's result for i.i.d. random environments \cite[Theorem A.2]{Sz02}. The proof is based on a combination of \cite[Theorems A.1-A.2]{Sz02} and \cite[Proposition  4.5]{Gue17}.
\end{remark}
\noindent \emph{Proof.}[Proof of Proposition \ref{protransfluctuation}]
We let $\varsigma$ and $\eta>0$ be as in the statement. Notice that $2\varsigma-1<\varsigma$ and choose $\gamma\in (0,1)$ such that $2\varsigma-1<\gamma \varsigma$ and $(T^\gamma)|\ell$ holds, which is possible since we assume that $(T')|\ell$ holds.
We take $M>0$ large enough so that $L\ge L_0$ with $L_0$ given in the statement of Corollary \ref{corexp}.  

Observe now that in order to prove claim (\ref{claimproptranv}), we can replace
$$
\sup_{0\leq n\leq \mathcal L_M}|\Pi_{\hat v_L} (X_n)|_2
$$
by $\sup_{0\leq n\leq\mathcal L_M}X_n\cdot w$, where $w\in \mathbb S^{d-1}$ and $w\cdot \hat v_L=0$. We consider the approximate regeneration time $\tau_1^{(L)}$ constructed as in Section \ref{Apre} along the direction $l$ and define for $n\geq0$ (recall convention $\tau_0^{(L)}:=0$):
\begin{equation*}
K_n=\sup\{k\geq 0:\,\tau_k^{(L)}\leq n \}.
\end{equation*} 
Setting $\nu_l:=|l|_1/|l|_2$, we see by the very definition of $\tau_1^{(L)}$ that $\overline P_0$-a.s. on the event $0\leq n \leq \mathcal L_M$ one has $K_n\leq (\nu_lM)/L$. As a result,  writing $X_n\cdot w=(X_n-X_{\tau_{K_n}})\cdot w + X_{\tau_{K_n}}\cdot w$ and recalling notation (\ref{supt1}), we get
\begin{gather}
\nonumber
P_0\left[\sup_{0 \leq n \leq \mathcal L_M}\,X_n\cdot w>\eta M^{\varsigma}\right]\leq \sum_{0\leq k\leq\frac{\nu_l M}{L} }\overline P_0\left[X^* \circ \theta_{\tau_k}>\frac{\eta M^\varsigma}{3}\right]+\\
\label{decomptranflucv}
\overline P_0\left[X_{\tau_1}\cdot w> \frac{\eta M^\varsigma}{3}\right]+\sum_{2\leq k\leq \frac{\nu_l M}{L}}\overline P_0\left[(X_{\tau_k}-X_{\tau_1})\cdot w>\frac{\eta M^\varsigma}{3}\right].
\end{gather}
The first two terms on the right hand side of inequality (\ref{decomptranflucv}) can be bounded from above using condition $(T^\gamma)|\ell$. More precisely, we first introduce constants $b$ and $\phi$ via:
\begin{gather}
\label{constantphiandb}
b=\frac{2\ln(1/\kappa)}{gt-2\ln(1/\kappa)} \,\, \mbox{ and } \,\,
\phi=e^{-\frac{bgt|l|_1}{2(b+1)}},
\end{gather}
then applying, for integer $k\in [0, \nu_lM/L]$, Proposition \ref{propare} in combination with Chebyshev's inequality to get
\begin{gather*}
\overline P_0\left[X^* \circ \theta_{\tau_k}>\frac{\eta M^\varsigma}{3}\right]=\overline P_0\left[c_6(\phi M^{-\epsilon}X^* \circ \theta_{\tau_k})^\gamma>\frac{c_{6}(\eta \phi)^\gamma M^{\gamma(\varsigma-\epsilon)}}{3}\right]\\
\leq \exp \left(-\frac{c_{6}(\eta \phi)^\gamma M^{\gamma (\varsigma-\epsilon)}}{3^\gamma}\right)\overline E_0\left[\exp\left(c_{6}(\phi M^{-\epsilon}X^* \circ \theta_{\tau_k})^\gamma\right)\right]\\
\leq e^{\exp\left(-gtL\right)}\exp \left(-\frac{c_6(\eta\phi)^\gamma M^{\gamma (\varsigma-\epsilon)}}{3^\gamma}\right)\overline E_0\left[\exp\left(c_6(\phi M^{-\epsilon}X^* )^\gamma\right)|D'=\infty\right].
\end{gather*}
Plainly, the same upper bound holds for the second term on the right side of inequality (\ref{decomptranflucv}). Thus
\begin{gather}
\nonumber
\sum_{0\leq k\leq\frac{\nu_l M}{L} }\overline P_0\left[X^* \circ \theta_{\tau_k}>\frac{\eta M^\varsigma}{3}\right]+\overline P_0\left[X_{\tau_1}\cdot w> \frac{\eta M^\varsigma}{3}\right]\\
\label{boundfirsttwo}
\leq \left(\frac{4M\nu_l}{L}\right)\exp \left(-\frac{c_6(\eta\phi)^\gamma M^{\gamma (\varsigma-\epsilon)}}{3^\gamma}\right)\overline E_0\left[\exp\left(c_6(\phi M^{-\epsilon}X^* )^\gamma\right)|D'=\infty\right].
\end{gather} 

Let us now proceed by examining the order of the last term on the right-most sum in (\ref{decomptranflucv}). For this end, we let $2\leq k\leq (\nu_lM)/L$ and notice that defining $\mathcal H_{k,M}:=\{\exists j\in [2, k]:\, |(X_{\tau_{j}}-X_{\tau_{j-1}})\cdot w| \geq \delta M^{\frac{2\varsigma-1}{\gamma}}\}$ for $\delta>0$ to be chosen later on, we have
\begin{gather}
\nonumber
\overline P_0\left[(X_{\tau_k}-X_{\tau_1})\cdot w>\frac{\eta M^\varsigma}{3}\right]\leq \\
\label{decotran2}
\overline P_0\left[ \mathcal H_{k,M}\right]+\overline P_0\left[(X_{\tau_k}-X_{\tau_1})\cdot w>\frac{\eta M^\varsigma}{3}, (\mathcal H_{k,M})^c \right].
\end{gather}
We first bound from above the first term on the right hand side of (\ref{decotran2}).
Observe that an argument close to that leading to (\ref{boundfirsttwo}) gives us
\begin{gather}
\nonumber
\overline P_0\left[ \mathcal H_{k,M}\right]\leq \sum_{2\leq j\leq k}\overline P_0\left[|(X_{\tau_j}-X_{\tau_{j-1}})\cdot w|\geq\delta M^{\frac{M^{2\varsigma-1}}{\gamma}}\right]\\
\nonumber
\leq \sum_{2\leq j\leq k}2\exp\left(-(\delta\phi)^\gamma c_6M^{2\varsigma-1-\gamma\epsilon}\right)\overline E_0\left[ \exp\left(c_6(\phi M^{-\epsilon}X^*)^\gamma\right)|D'=\infty\right]\\
\label{boundintwhofirst}
\leq 2k\exp\left(-c_6(\delta\phi)^\gamma M^{2\varsigma-1-\gamma\epsilon}\right)\overline E_0\left[ \exp\left(c_6(\phi M^{-\epsilon}X^*)^\gamma\right)|D'=\infty\right].
\end{gather}
We now turn to bound the remaining (and hardest to handle) term. However, we must keep (\ref{boundintwhofirst}) in mind and continue fixing $k\in[2, (\nu_lM)/L]$. For that purpose, we fix some $\zeta\in (1/2,1]$ and apply a Chernoff bound to the second term on right most hand of inequality (\ref{decotran2}), getting
\begin{gather}
\label{decomhard}
\overline P_0\left[(X_{\tau_k}-X_{\tau_1})\cdot w>\frac{\eta M^\varsigma}{3}, (\mathcal H_{k,M})^c \right]\leq
e^{-\zeta\frac{\eta M^\varsigma}{M^{(1-\varsigma)+\epsilon}}}\times\\
\nonumber
\left(\overline E_0\left[e^{\zeta\frac{\sum_{2\leq j\leq k}(X_{\tau_j}-X_{\tau_{j-1}})\cdot w}{M^{(1-\varsigma)+\epsilon}}}, \forall i\in[2,k],\,|(X_{\tau_i}-X_{\tau_{i-1}})\cdot w|<\delta M^{\frac{2\varsigma-1}{\gamma}} \right]\right).
\end{gather}
On the other hand, 
\textbf{(SMG)}$_{C,g}$ along with
Proposition \ref{propare} and successive conditioning lead us to:
\begin{gather*}
\left(\overline E_0\left[e^{\zeta\frac{\sum_{2\leq j\leq k}(X_{\tau_j}-X_{\tau_{j-1}})\cdot w}{M^{(1-\varsigma)+\epsilon}}}, \forall i\in[2,k],\,|(X_{\tau_i}-X_{\tau_{i-1}})\cdot w|<\delta M^{\frac{2\varsigma-1}{\gamma}} \right]\right)\\
\leq \left(\exp\left(e^{-gtL}\right)\overline E_0\left[\exp\left(\frac{\zeta X_{\tau_1}\cdot w}{M^{1-\varsigma+\epsilon}}\right), |X_{\tau_1}\cdot w| <\delta M^{\frac{2\varsigma-1}{\gamma}} \Big|D'=\infty\right]\right)^{\frac{\nu_l M}{L}},
\end{gather*}
for $2 \leq k\leq (\nu_lM)/L$. Using that for $|u|\leq 1$ one can find some $\nu\in (1/2,1]$ such that
\begin{equation*}
|e^{u}-1-u|<\nu u^2,
\end{equation*}
and choosing $\alpha>\epsilon/2$, with $\alpha<(2\varsigma-1)/\gamma\wedge (1-\varsigma+\epsilon)$,  we have
\begin{gather}
\nonumber
 \overline E_0\left[\exp\left(\frac{\zeta X_{\tau_1}\cdot w}{M^{1-\varsigma+\epsilon}}\right), |X_{\tau_1}\cdot w| <\delta M^{\frac{2\varsigma-1}{\gamma}} \Big|D'=\infty\right]\\
 \nonumber
\leq \overline E_0\left[1 +\frac{\zeta X_{\tau_1}\cdot w}{M^{(1-\varsigma)+\epsilon}}+\nu \frac{\zeta^2 (X_{\tau_1}\cdot w)^2}{M^{(2-2\varsigma)+2\epsilon}},\,|X_{\tau_1}\cdot w|\leq M^{\alpha} \Big|D'=\infty\right]\\
\nonumber
+\overline E_0\left[\exp\left(\frac{\zeta X_{\tau_1}\cdot w}{M^{1-\varsigma +\epsilon}}\right), M^\alpha<|X_{\tau_1}\cdot w|<\delta M^{\frac{2\varsigma-1}{\gamma}} \Big| D'=\infty\right]\\
\nonumber
\leq 1+ \nu \frac{\zeta^2}{M^{2-2\varsigma+\epsilon}}\overline E_0\left[\left(M^{-\frac{\epsilon}{2}}X_{\tau_1}\cdot w\right)^2 \Big| D'=\infty\right]-\frac{\zeta}{M^{1-\varsigma+(\epsilon/2)}}\times\\
\nonumber
\overline E_0\left[M^{-\frac{\epsilon}{2}}X_{\tau_1}\cdot w,\, |X_{\tau_1}\cdot w|> M^{\alpha} \Big| D'=\infty\right]\\
\label{hh}
+\overline E_0\left[\exp\left(\frac{\zeta X_{\tau_1}\cdot w}{M^{1-\varsigma +\epsilon}}\right), M^\alpha<|X_{\tau_1}\cdot w|<\delta M^{\frac{2\varsigma-1}{\gamma}}\Big| D'=\infty\right].
\end{gather}
Applying Cauchy-Schwarz inequality to the third term on the right hand side of  (\ref{hh}), we have
\begin{gather*}
\overline E_0\left[M^{-\frac{\epsilon}{2}}X_{\tau_1}\cdot w, |X_{\tau_1}\cdot w|> M^\alpha|D'=\infty\right]\\
\leq \overline E_0\left[\left(M^{-\frac{\epsilon}{2}}X_{\tau_1}\cdot w\right)^2|D'=\infty\right]^{\frac{1}{2}}\left(\overline P_0\left[|X_{\tau_1}\cdot w|> M^\alpha|D'=\infty\right]^{\frac{1}{2}}\right).
\end{gather*}
We then use condition $(T')|\ell$ on the probability term above, and we finally get that the third term on the right side of (\ref{hh}) is bounded by
\begin{gather}
\nonumber
\overline E_0\left[\left(M^{-\frac{\epsilon}{2}}X_{\tau_1}\cdot w\right)^2|D'=\infty\right]^{\frac{1}{2}}\\
\label{bound3}
\times \exp(-\frac{c_{6}}{2}\phi^\gamma M^{\gamma(\alpha-\epsilon/2)})\overline E_0\left[\exp\left(c_{6}(\phi M^{-\frac{\epsilon}{2}}X^*)^\gamma\right)|D'=\infty\right]^{\frac{1}{2}},
\end{gather}
where we have proceeded as in (\ref{boundfirsttwo})-(\ref{boundintwhofirst}) in the last inequality.
At this point we are going to leave the main proof subject to prove that the expectations entering at inequalities (\ref{boundfirsttwo})-(\ref{boundintwhofirst}) are bounded. From the very definition of $L$ at the beginning of the proof, one has 
\begin{equation*}
e^{-\frac{bgt|l|_1}{2(b+1)}}M^{-\frac{\epsilon}{2}}=e^{-\frac{bgt|l|_1}{2(b+1)}}M^{-\frac{b(2-2\varsigma)}{2}}<e^{-\frac{bgtL}{2(b+1)}}=\kappa^L.
\end{equation*}

As a result, using first the bound in (\ref{boundfirsttwo}) and with the help of Corollary \ref{corexp}, we obtain
\begin{gather}
\nonumber
\sum_{0\leq k\leq\frac{\nu_l M}{L} }\overline P_0\left[X^*\circ \theta_{\tau_k}>\frac{\eta M^\varsigma}{3}\right]+\overline P_0\left[X_{\tau_1}\cdot w> \frac{\eta M^\varsigma}{3}\right]\\
\label{recastbounds1}
\leq \left(\frac{4M\nu_l}{L}\right)\exp \left(-\frac{c_{6}(\eta\phi)^\gamma M^{\gamma (\varsigma-\epsilon)}}{3^\gamma}\right)c_{7}.
\end{gather}
Analogously, using bounds in (\ref{boundintwhofirst}) and then the bound in (\ref{bound3}), Corollary \ref{corexp} provides the bounds:
\begin{align}
\label{recastbounds2}
\overline P_0\left[ \mathcal H_{k,M}\right]\leq 2k\exp\left(-c_{6}(\delta\phi)^\gamma M^{2\varsigma-1-\gamma\epsilon}\right)c_{7},
\end{align}
along with:
\begin{align}
\label{recastbounds3}
&\overline E_0\left[M^{-\frac{\epsilon}{2}}X_{\tau_1}\cdot w, |X_{\tau_1}\cdot w|> M^\alpha|D'=\infty\right]\\
\nonumber
&\leq\exp\left(-\frac{c_6}{2}\phi^\gamma M^{\gamma(\alpha-\epsilon/2)}\right)(c_{7}/c_{6}\phi).
\end{align}
We now go back to the last term on the right side of (\ref{hh}). Writing
\begin{gather}
\label{powerint}
\overline E_0\left[\exp\left(\frac{\zeta X_{\tau_1}\cdot w}{M^{1-\varsigma +\epsilon}}\right), M^\alpha<|X_{\tau_1}\cdot w|<\delta M^{\frac{2\varsigma-1}{\gamma}}|D'=\infty\right]
\end{gather}
as a sum of expectations over two sets: $\{-\delta M^{\frac{2\varsigma-1}{\gamma}}<X_{\tau_1}\cdot w<-M^\alpha\}$ and $\{M^\alpha<X_{\tau_1}\cdot w<\delta M^{\frac{2\varsigma-1}{\gamma}}\}$, using a similar analysis as the one led to (\ref{bound3}), the expectation on the first set is bounded by 
\begin{equation*}
\exp\left(-c_6 \phi^\gamma M^{\gamma(\alpha-\epsilon/2)}\right)c_7.
\end{equation*}
While the expectation of the integrand in (\ref{powerint}) over the another set: $\{M^\alpha<X_{\tau_1}\cdot w<\delta M^{\frac{2\varsigma-1}{\gamma}}\}$, equals
\begin{gather*}
\frac{\zeta}{M^{1-\varsigma+(\epsilon/2)}}\times\\
\int_{M^{\alpha-(\epsilon/2)}}^{\delta M^{\frac{2\varsigma-1}{\gamma}-(\epsilon /2)}}\exp\left(\zeta\frac{u}{M^{1-\varsigma+(\epsilon/2)}}\right)\overline P_0[M^{-\frac{\epsilon}{2}}X_{\tau_1}\cdot w>u|D'=\infty]du.
\end{gather*}
Thus, as a result of a further application of $(T^\gamma)|\ell$ on the probability inside of the integrand above, we get that (\ref{powerint}) is bounded by
\begin{gather}
\label{layercake1}
\exp\left(-c_6 \phi^\gamma M^{\gamma(\alpha-\epsilon/2)}\right)c_7+\frac{\zeta'}{M^{1-\varsigma+(\epsilon/2)}}\\
\times\int_{M^{\alpha-(\epsilon/2)}}^{\delta M^{\frac{2\varsigma-1}{\gamma}-(\epsilon /2)}}\exp\left(\zeta\frac{u}{M^{1-\varsigma+(\epsilon/2)}}-\xi u^\gamma\right)du,
\end{gather}
where $\xi=c_6\phi^\gamma$ and $\zeta'=\zeta \overline E_0\left[e^{c_6(\phi M^{-\frac{\epsilon}{2}}X_{\tau_1}\cdot w)^\gamma}\right]$.
\noindent
Observe now that $(2\varsigma-1)/\gamma+(\varsigma-1)-\epsilon<\gamma \varsigma$, since $2\varsigma-1<\gamma \varsigma$. Using this, we find that the second term in (\ref{layercake1}), for large $M$, is bounded by
$$
\frac{\zeta'}{M^{1-\varsigma+(\epsilon/2)}}\int_{M^{\alpha-(\epsilon/2)}}^{\infty}\exp\left(-\frac{\xi u^\gamma}{2}\right)du\leq \exp\left(-\frac{\xi M^{\gamma(\alpha-(\epsilon/2))}}{3}\right).
$$
In view of this last inequality and (\ref{bound3}), we go back to (\ref{hh}) and find that for large $M$,
\begin{gather*}
\overline E_0\left[\exp\left(\frac{\zeta X_{\tau_1}\cdot w}{M^{1-\varsigma+\epsilon}}\right), |X_{\tau_1}\cdot w| <\delta M^{\frac{2\varsigma-1}{\gamma}}|D'=\infty\right]\\
\leq \exp(2\nu \zeta^2 (c_{7}/c_{6}^{2})M^{2\varsigma-2-\epsilon}).
\end{gather*}
Moreover, we apply this last inequality to (\ref{decomhard}) and we get that for each integer $k\in[2,(\nu_lM)/L]$, 
\begin{equation}\label{desinterlac1}
\overline P_0\left[(X_{\tau_k}-X_{\tau_1})\cdot w,\ (\mathcal H_{k,M})^c\right]\leq e^{-\Big(\zeta\eta - \frac{2 \nu \zeta^2 (c_{7}/c_{6}^{2})}{L} \Big) M^{2\varsigma-1-\epsilon}}.
\end{equation}
On the other hand, taking $L\in |l|_1\mathbb N$ sufficiently large, one has that
$$
-\zeta \eta +\left(2\nu \zeta^2 c_{7}/c_{6}^2 +1\right)\frac{\nu_l}{L}<-(\zeta\eta)/4.
$$
Applying the inequality above to (\ref{desinterlac1}), combining the resulting inequality and (\ref{boundintwhofirst}) in (\ref{decotran2}), after summing over $k\in[2,(\nu_lM)/L]$, we have that for large $M$
\begin{equation}\label{ineq1tran}
\sum_{2\leq k\leq \nu_lM/L}\overline P_0\left[(X_{\tau_k}-X_{\tau_{1}})\cdot w>\frac{\eta M^{\varsigma}}{3}\right]\leq \exp\left(-\frac{\zeta\eta M^{2\varsigma-1-\epsilon}}{8}\right).
\end{equation}
Thus, we have finished the proof provided that we combine (\ref{recastbounds1}), (\ref{recastbounds2}), (\ref{recastbounds3}) and (\ref{ineq1tran}) with (\ref{decomptranflucv}). \hfill $\square$



\medskip
For $L\in |l|_1\mathbb N$, in order to bound tails of the random variable $\tau_1^{(L)}$ under assumption $(T')|\ell$, we  follow \cite[Section 5]{Gue17} (see also \cite[Section III]{Sz00} for the original argument for i.i.d. environments).

\smallskip
The next lemma is a first connection between condition $(T')|\ell$ and tails of regeneration times. We fix a rotation $R$ of $\mathbb R^d$ with $R(e_1)=\ell$ and as in Section \ref{secPrel}, we choose $\mathfrak r$ such that (\ref{Tgammasquare}) is satisfied.

\smallskip

We define for $M>0$, the hypercube (see (\ref{generalboxes}) for notation):
\begin{equation}\label{hypercube}
C_M:=B_{M,\mathfrak rM,\ell}(0),
\end{equation}
we then have the next result:

\begin{lemma}\label{decomtail1}
There exist constants $c_9:=c_9(d,\kappa,g)>0$ and $L_0>0$ such that for any function $M:\mathbb R^+ \mapsto  \mathbb R^+$, with $\lim_{\substack{u\rightarrow\infty}}M(u)=\infty$, we have that for large $u$ and $\gamma\in(0,1)$,
\begin{equation}\label{integratingtau1}
\overline P_0\left[\tau_1>u\right]\leq P_0\left[T_{C_M}=T_{M(u)}^\ell>u\right]+e^{-c_{9}\left(\kappa^LM(u)\right)^\gamma}.
\end{equation}
\end{lemma}
\noindent \emph{Proof.}
For large $u$, using Proposition \ref{expmpr} we have
\begin{gather*}
\overline P_0\left[\tau_1>u\right]\leq\\
\overline P_0\left[\tau_1>u,\, X_{\tau_1}\cdot l <|l|_2 M(u)\right]+\overline P_0\left[X_{\tau_1}\cdot l\geq |l|_2M(u)\right]\leq\\
\overline P_0\left[\tau_1>u,\, X_{\tau_1}\cdot l <|l|_2 M(u)\right]+e^{-\frac{c_2|l|_2^\gamma\left(\kappa^L M(u)\right)^\gamma}{2}}.
\end{gather*}
On the other hand, using that $\tau_1=T_{X_{\tau_1}\cdot l}^l$ we find that:
\begin{gather*}
\overline P_0\left[\tau_1>u,\, X_{\tau_1}\cdot l <|l|_2 M(u)\right]\leq P_0\left[T_{M(u)}^\ell>u\right] \\
\leq P_0\left[T_{M(u)}^\ell=T_{C_{M(u)}}>u \right]+ P_{0}\left[X_{T_{C_{M(u)}}}\notin \partial^+ C_{M(u)}\right].
\end{gather*}
Using $(T')|\ell$ we obtain
\begin{gather*}
P_{0}\left[X_{T_{C_{M(u)}}}\notin \partial^+ C_{M(u)}\right]\leq e^{-\widetilde c (M(u))^\gamma}
\end{gather*}
where $\widetilde c>0$ is certain constant. The inequality (\ref{integratingtau1}) follows. \hfill $\square$

\medskip
We continue with a version of the so-called atypical quenched estimate. Its construction strongly depends on Proposition \ref{protransfluctuation}. Furthermore, the next result shall be the cornerstone to get the final estimate. 
We first introduce the set
\begin{equation}\label{slabUl}
U_{M}=\left\{x\in \mathbb Z^d:\, |x\cdot \ell|<M\right\}
\end{equation}
for $M>0$. We have: 
\begin{proposition}\label{propatyquenest}
For $c>0$ and
\begin{equation}
\label{betbig}
\beta \in \Big( \frac{2gt}{3gt-2\ln(1/ \kappa)} , 1 \Big)
\end{equation}
 with $t$ as in \eqref{defintransfluct}, we have that
\begin{equation}\label{atyquenest}
\limsup_{\substack{M\rightarrow\infty}}M^{-\chi}\ln \mathbb P\left[P_{0,\omega}\left[X_{T_{U_M}}\cdot \ell>0\right]\leq e^{-cM^\beta}\right]<0,
\end{equation}
where 
$\chi < \frac{dgt(3\beta - 2)}{gt-2\ln(1/\kappa)}$.
\end{proposition}

We postpone the proof of Proposition \ref{propatyquenest} to Section \ref{secproofaqe} which is entirely devoted to it.

\medskip
We finally finish the procedure to bound tails of $\tau_1^{(L)}$ for $L\in |l|_1\mathbb N$. Roughly speaking, it is a Markov chain argument which links Proposition \ref{propatyquenest} and Lemma \ref{decomtail1} with the upper bound in the statement of Theorem \ref{thboundtails}. 

\medskip

\noindent \emph{Proof.} [Proof of Theorem \ref{thboundtails}]
Let $\alpha\in (1,1+\sigma)$ and for large $u$ define
\begin{gather}
\nonumber
  \Delta(u)=\frac{1}{10\sqrt{d}}\, \frac{\ln(u)}{\ln\left(\frac{1}{\kappa}\right)}\,\, \mbox{ and }\,\,M(u)= N(u)\Delta(u), \\
  \label{scalestail}
  \hspace{15ex}\mbox{where } N(u)=\left[(\ln(u))^{\alpha-1}\right].
\end{gather}
Henceforward, we shall omit the dependence on $u$ of $M, \,\Delta$ and $N$. We apply Lemma \ref{decomtail1} using the function $M$ defined above, to get
\begin{equation*}
\overline P_0\left[\tau_1^{(L)}>u\right]\leq e^{-\hat c_{1}\kappa^L\left(\ln(u)\right)^{\gamma \alpha}}+P_0\left[T_{C_M}>u\right].
\end{equation*}
Since $\gamma$ can be taken close to $1$ and the upper bound for $\alpha$ is not attained, it suffices to prove that
\begin{equation}\label{claimlemmatail}
\limsup_{\substack{M\rightarrow\infty}}\ln(u)^{-\alpha}\ln\left(P_0\left[T_{C_M}>u\right]\right)<0.
\end{equation}
We decompose $P_0\left[T_{C_M}>u\right]$ as follows:
\begin{gather}
\nonumber
P_0\left[T_{C_M}>u\right]\leq \mathbb E\left[P_{0,\omega}\left[T_{C_M}>u\right],\, \forall x\in C_M\, P_{x,\omega}\left[T_{C_M}> \frac{u}{(\ln(u))^\alpha}\right]\leq \frac{1}{2}\right]+\\
\label{decomtauitail}
\mathbb E\left[P_{0,\omega}\left[T_{C_M}>u\right],\, \exists x_1\in C_M\, P_{x_1,\omega}\left[T_{C_M}>\frac{u}{(\ln(u))^\alpha}\right]>\frac{1}{2}\right].
\end{gather}
For the first term on the right hand of the inequality (\ref{decomtauitail}), we apply the strong Markov property and we see that
$$
\mathbb E\left[P_{0,\omega}\left[T_{C_M}>u\right],\, \forall x\in C_M\, P_{x,\omega}\left[T_{C_M}> \frac{u}{(\ln(u))^\alpha}\right]\leq \frac{1}{2}\right]\leq \left(\frac{1}{2}\right)^{[\ln(u)^\alpha]}.
$$
We turn to bound the second term on the right most hand in (\ref{decomtauitail}). Notice that, by the argument in \cite[pages 127 and 128]{Sz00}, for large $u$
\begin{gather}
\nonumber
\mathcal R:=\left\{\omega\in \Omega: \ \exists x_1\in C_M, \  P_{x_1,\omega}\left[T_{C_M}>\frac{u}{(\ln(u))^\alpha}\right]>\frac{1}{2} \right\}\\
\label{contained2}
\subset \{ \omega\in \Omega: \ \exists x, i: \ i\in \mathbb [-N+2,N-1], \ x\in C_M\cap \mathcal G_i, \\
\nonumber
 P_{x,\omega}[\widetilde T_{(i-1)\Delta}^\ell>T_{C_M}]\leq 1/\sqrt u\}.
\end{gather}
For $i\in \mathbb Z$, we define random variables
\begin{equation*}
X_i:=\left\{
\begin{array}{ll}
-\ln\left( \inf_{x\in C_M \cap \mathcal G_i }
P_{x,\omega}\left[\widetilde{T}_{(i-1)\Delta}^\ell>T_{(i+1)\Delta}^\ell\right]\right)
&\quad{\rm if}\quad C_M\cap \mathcal G_i\neq\emptyset,\\ 0&\quad{\rm
if}\quad C_M\cap \mathcal G_i=\emptyset.
\end{array}
\right.
\end{equation*}
For $i\in [-N+1, N]$ and $x\in \mathcal G_i$, a standard application of the strong Markov property  yields
\begin{equation*}
P_{x,\omega}\left[\widetilde T_{(i-1)\Delta}>T_{C_M}\right]\geq \exp\left(-\sum_{j=i}^N X_i\right).
\end{equation*}
Therefore,
\begin{equation}
\label{probrarevent}
\mathbb P\left[\mathcal R\right]\leq \mathbb P\left[\sum_{i=-N+1}^N X_i\geq \frac{\ln (u)}{2}\right]\leq 2N \sup_{\substack{i\in [-N+1, N]}}\mathbb P\left[X_i\geq \frac{\ln(u)}{4N}\right].
\end{equation}
On the other hand, observe that for arbitrary $i\in \mathbb Z$ and $\vartheta>0$ we have (see \ref{slabUl} for notation)
\begin{equation}\label{atypquen}
\mathbb P\left[X_i>\vartheta\right]\leq |C_M|\mathbb P\left[P_{0,\omega}\left[X_{T_{U_\delta}}\cdot\ell \geq \Delta\right]\leq e^{-\vartheta}\right].
\end{equation}
Hence, we apply (\ref{atypquen}) into inequalities (\ref{probrarevent}) and then we use Proposition \ref{propatyquenest} to conclude that for large $u$ there exists a suitable constant $\widetilde c$ such that (see (\ref{defintransfluct}) for notation)
\begin{equation}
\mathbb P\left[\mathcal R\right]\leq \exp\left(-\widetilde c \ln(u)^\chi\right),
\end{equation}
where $\chi<  \frac{dgt(3(2-\alpha) - 2)}{gt-2\ln(1/\kappa)}$ with $(2-\alpha)$, that plays the role of $\beta$ in the statement of Proposition \ref{propatyquenest}, satisfying 
\begin{equation}\label{req1}
2-\alpha>\frac{2gt}{3gt-2\ln(1/\kappa)}\Leftrightarrow \alpha<1+\frac{gt-2\ln(1/\kappa)}{3gt-2\ln(1/\kappa)}.
\end{equation}
In turn, we require that
\begin{gather}
\label{req2}
\alpha< \frac{dgt(3(2-\alpha) - 2)}{gt-2\ln(1/\kappa)} 
\Leftrightarrow
\alpha<1+\frac{(d-1)gt +2\ln(1/\kappa)}{(3d+1)gt-2\ln(1/\kappa)}.
\end{gather}
The theorem now follows from (\ref{req1})-(\ref{req2}). \hfill $\square$

\subsection{Proof of the atypical quenched estimate}
\label{secproofaqe}

This section is entirely devoted to the proof of Proposition \ref{propatyquenest}.

\medskip
Notice that if $\beta$ satisfies (\ref{betbig}), then $\beta > \frac{(2d+1)gt-2\ln(1/\kappa)}{3dgt}$ and therefore
$$
 \frac{dgt(3\beta - 2)}{gt-2\ln(1/\kappa)} > 1.
$$
With the help of Proposition \ref{protransfluctuation}, we now follow closely the proof of  \cite[Proposition 5.2]{Gue17} which is an extension of arguments used in \cite{Sz00} for i.i.d. environments.
For large $M$, we shall construct strategies which will happen with high probability on the environment law ensuring that the random walk escapes from slab $U_M$ by the boundary side
$$
\partial^+U_M:=\partial U_M\cap\left\{z:\,z\cdot \ell \geq M\right\}
$$
with quenched probability bigger than $e^{-cM^\beta}$. Recall definitions (\ref{defintransfluct})-(\ref{epsilonchoise}) and let $\gamma$ be a real number in the interval:
$$
\left(\frac{gt+2\ln(1/\kappa)}{2gt},1\right).
$$
We pick $M_0>3\sqrt d$ large enough such that the number $L$ defined by \eqref{LM} verifies $L\geq L_0$ ($L_0$ being the maximum between the ones of Proposition \ref{propare}, Corollary \ref{corexp} and Proposition \ref{protransfluctuation}).
For that given $L$ one then chooses a rotation $\widehat R$ of $\mathbb R^d$ with
$$
\widehat R(e_1)=\widehat v_L,
$$
with $v_L$ given by \eqref{veloc}, introduced just before Theorem \ref{thboundtails}.
We now introduce for $z\in M_0\,\mathbb Z^d$ the following blocks:
\begin{gather}
\widetilde{B}_1(z):=\widehat R\left(z+(0,M_0)^d\right)\cap \mathbb Z^d \hspace{2ex}\mbox{and}\\
\label{blocksaqe}
\widetilde{B}_2(z):=\widehat R\left(z+(-M_0^\gamma, M_0+M_0^\gamma)^d\right)\cap \mathbb Z^d.
\end{gather}
with their frontal boundary defined as in \eqref{fb+}.
The aforementioned strategy involves a definition of \textit{good} and \textit{bad} blocks. However, differently than the process in Section \ref{secpoly}, now we have a one-step renormalization.
\smallskip
We say that a site $z\in M_0\,\mathbb Z^d$ is $M_0$-good if
\begin{equation}\label{goodbox}
\sup_{\substack{x\in \widetilde B_1(z)}}P_{x,\omega}\left[X_{T_{B_2(z)}}\in \partial^+B_2(z)\right]\geq \frac{1}{2}
\end{equation}
and $M_0$-bad otherwise.

\begin{lemma}\label{lemmagoodbox}
Let $\gamma\in \left(\frac{gt+2\ln(1/\kappa)}{2gt},1\right)$. Then one has that
\begin{equation}\label{estbadbox}
\limsup_{\substack{M_0\rightarrow\infty}}M_0^{-(2\gamma-1-\epsilon)}\,\sup_{\substack{z\in M_0\, \mathbb Z^d}} \ln \mathbb P\left[z\,\mbox{ is }\,M_0-\mbox{bad}\right]<0,
\end{equation}
where $\epsilon$ was defined in \eqref{epsilonchoise}.
\end{lemma}
\noindent \emph{Proof.}
For $z\in M_0\, \mathbb Z^d$,
\begin{equation*}
\mathbb P\left[z\,\mbox{ is }\,M_0-\mbox{bad}\right]\leq 2|\widetilde B_1(z)|\sup_{\substack{x\in \widetilde B_1(z)}}P_x\left[X_{T_{\widetilde B_2(z)}}\notin \partial^+\widetilde B_2(z)\right].
\end{equation*}
Observe now that for arbitrary $x\in \widetilde B_1(z)$, the block $\widetilde B_{2}(z)$ is included in a ball of radius $3\sqrt dM_0$ centered at $x$. Thus,  $P_x$-a.s we have
$$
T_{\widetilde B_2(z)}\leq T_{x\cdot \ell +3\sqrt d M_0}^\ell.
$$
Also, notice that on the event $\left\{X_{T_{\widetilde B_2(z)}}\notin \partial^+\widetilde B_2(z)\right\}$ we have $P_x$-a.s. either
\begin{equation*}
\left(X_{T_{\widetilde B_2(z)}}-x\right)\leq -\frac{M_0^\gamma}{2}\hspace{2ex}\mbox{or}\hspace{2ex}\left|\Pi_{\widehat v}\left(X_{T_{\widetilde B_2(z)}}-x\right)\right|_2\geq \frac{M_0^\gamma}{2},
\end{equation*}
where the notations as at the beginning of this section, see (\ref{projvl}). Consequently as in the proof of \cite[Lemma 5.3]{Gue17} for a suitable constant $c(d)$ we find
\begin{gather}\label{lastinbadboxest}
\mathbb P\left[z\,\mbox{ is }\,M_0-\mbox{bad}\right]\leq \\
\nonumber
c(d)M_0^d \left(P_0\left[\sup_{\substack{0\leq n\leq T_{3\sqrt d M_0}^\ell}}\Pi_{\widehat v}\left(X_n\right)\geq \frac{\widehat v\cdot \ell M_0^\gamma }{4}\right]+P_0\left[\widetilde T_{-\frac{\widehat v\cdot \ell M_0^\gamma }{4}}^\ell<\infty\right]\right).
\end{gather}
We construct the random variable $\tau_1:=\tau_1^{(L_0)}$ along direction $\ell$ change by $l$ with $L_0$ as in Proposition \ref{expmompos}, and using that
$$
P_0\left[\widetilde T_{-\frac{\widehat v\cdot \ell M_0^\gamma }{4}}^\ell<\infty\right]\leq \overline P_0\left[X_{\tau_1}\cdot \ell \geq \frac{\widehat v\cdot \ell M_0^\gamma }{4}\right],
$$
it is clear that a Chernoff bound with the help of Proposition \ref{expmpr}, bounds $\gamma$-stretched exponential the second term in the right most hand of inequality \ref{lastinbadboxest}. The assertion (\ref{estbadbox}) follows since the existence of non-dependent on $L$ constants $k_1>0$ and $k_2>0$ (see \cite[Proposition 7.2]{GR17}) such that
$$
k_1\leq \widehat v_L\cdot \ell \leq k_2
$$
together with applying Proposition \ref{protransfluctuation}. \hfill $\square$


\medskip
It is convenient to introduce some further terminology. Consider $M>0$, $M_0$ as above, and attach to each site $z\in M_0\, \mathbb Z^d$ the column
\begin{equation}\label{columnz}
Col(z):=\left\{z' \in M_0\,\mathbb Z^d: \exists j \in [0,J] \,  \mbox{ with }\,  z'=z + j M_0 e_1 \right\},
\end{equation}
where $J$ is the smallest integer satisfying $JM_0\widehat v_L\cdot \ell\geq 3M $. We choose $M_1>0$ an integer multiple of $M_0$ and we also define the tube attached to site $z\in M_0\,\mathbb Z^d$,
\begin{equation}\label{tubez}
Tube(z):=\left\{z'\in M_0\,\mathbb Z^d: \exists j_1,\ldots, j_d \in \left[0,\frac{M_1}{M_0}\right],\, z'=z+\sum_{i=1}^d j_iM_0e_i\right\}.
\end{equation}
The crucial point of the strategy is that one way for the walk starting from $0$ to escape from $U_M$ by $\partial^+U_M$ is to get to one of the bottom block in $Tube(0)$ containing the largest amount of good boxes, and then moving along this column up to its top. Thus, for $z\in M_0\mathbb Z^d$ we define
\begin{equation}\label{Topz}
Top(z):=\bigcup_{\substack{z'\in Tube(z)}}\partial^+\widetilde B_2(z'+JM_0e_1),
\end{equation}
together with the neighbourhood of a tube,
\begin{equation}
\label{neighbouz}
V(z):=\left\{x\in \mathbb Z^d:\exists y\in \bigcup_{\substack{z'\in Tube(z)\\ 0\leq j \leq J}}\widetilde B_1(z'+jM_0e_1),  |x-y|_1\leq 3dM_1\right\}.
\end{equation}
For a better understanding of notation introduced in \eqref{columnz}-\eqref{neighbouz} and the strategy to escape from $U_M$ through $\partial^+ U_M$, we recommend taking a look at \cite[Figure 2 in page 3037]{Gue17}. 
One has the following:
We define the first entrance to set $A\subset \mathbb Z^d$ by $H_A:=T_{\mathbb Z^d\setminus A}$. We have the following result:
\begin{lemma}
\label{lemmacompaqe}
For $z\in M_0\,\mathbb Z^d$, we let $n(z,w)$ be the random variable:
\begin{equation}\label{minbadblocks}
\min_{\substack{z'\in Tube(z)}}\left\{\sum_{j=0}^J \mathds{1}_{z'+jM_0e_1\hspace{0.5ex}\mbox{is }M_0-\mbox{bad}}\right\}.
\end{equation}
There exists $c_{10}>0$ such that for any $z\in M_0\, \mathbb Z^d$ and any $x\in D(z)$ where
$$
D(z):=\bigcup_{\substack{z'\in Tube(z)\\ 0\leq j\leq J}}\widetilde B_1(z'+jM_0e_1),
$$
we have
\begin{equation}\label{eqlemcomaqe}
P_{x,\omega}\left[H_{Top(z)}<T_{V(z)}\right]> (2\kappa)^{c_{10}\left(M_1+JM_0^\gamma+n(z,\omega)M_0\right)}\frac{1}{2^{J+1}}.
\end{equation}
\end{lemma}
It is easy to see that replacing $\kappa$ by $2\kappa$ in virtue of (\ref{simplex}) the proof of \cite[Lemma 3.3]{Sz00} provides the claim (\ref{eqlemcomaqe}).

\medskip
Keeping in mind Lemma \ref{lemmagoodbox}, we choose $\gamma\in ((gt+2\ln(1/\kappa))/2gt,1)$ such that
\begin{equation}\label{xidef}
\xi:=\frac{1-\beta}{1-\gamma}<\beta<1.
\end{equation}
Let us note that the choice of $\gamma$ is possible under assumption (\ref{betbig}), since:
\begin{equation*}
\beta>\frac{2gt}{3gt-2\ln(1/\kappa)}\,\Leftrightarrow\,\frac{2\beta-1}{\beta}>\frac{gt+2\ln(1/\kappa)}{2gt}.
\end{equation*}
We then choose $\nu>1-\gamma$
and, for large values of $M>0$, we put
\begin{equation*}
M_0=r_1 M^\xi,\hspace{2ex}M_1=\left\lfloor r_2M^{\beta-\xi}\right\rfloor M_0,\hspace{2ex}
N_0=\left\lfloor r_3M^{\beta-\xi}\right\rfloor,
\end{equation*}
where the constants $r_1,\,r_2,\,r_3$ possibly depend on the constant of the model and $c$ in \eqref{atyquenest}.
They are chosen so that for large $M$: 
\begin{gather}\label{constrain1aqe}
   (2\kappa)^{c_{10}JM_0^\gamma}, \, (2\kappa)^{c_{10}M_1}, \, (2\kappa)^{c_{10}N_0M_0}, \,\left(\frac{1}{2}\right)^{J+1}>\exp\left(-\frac{c}{5}M^\beta\right).\\
  \label{constrain2aqe}
  \frac{N_0}{3}>(J+1)\frac{(e^2-1)}{M_0^\nu}, \, \, \mbox{  and}\\
  \label{constrain3aqe}
  \mbox{any nearest neighbor path within $V(0)$, starting from the origin} \\
  \nonumber
  \mbox{and ending at $Top(0)$, first exits $U_M$ through $\partial^+ U_M$.}
\end{gather}
That choice is possible because in order to satisfy \eqref{constrain1aqe} and \eqref{constrain2aqe} it is sufficient to take $r_1$ large enough and then $r_2=r_3= \big(c_{10}r_1\ln\frac{1}{2\kappa}\big)^{-1}$. We also note that (\ref{constrain3aqe}) is fulfilled since $\beta<1$ implies $\beta -\xi>1-(1+\nu)\xi$.

We now borrow the last arguments in proof of  \cite[Proposition 5.2]{Gue17} to conclude that
\begin{equation*}
 \limsup_{M}\, M^{-d(\beta-\xi)}\,\ln \mathbb P\left[P_{0,\omega}\left[X_{T_{U_M}}\cdot \frac{l}{|l|_2}\geq M\right]\leq e^{-cL^\beta}\right]<0,
\end{equation*}
which ends the proof of claim (\ref{atyquenest}) provided we vary $\gamma$ according to (\ref{xidef}).

\appendix
\section{Proof of Lemma \ref{lemmaTgamma}} \label{apx}
\normalsize

This proof is very similar to that of \cite[Lemma  2.2]{Gue17}. However the last reference missed a proof for the equivalence of $(T^\gamma)|_{l_0}$ in Definition \ref{deftgammaandtprime} with conditions (i), (ii) and (iii). Thus for the sake of completeness we provide a proof here for the case $\gamma<1$.

\medskip

\noindent \emph{Proof.} [Proof of Lemma \ref{lemmaTgamma}]
The proof of $(i)\Rightarrow(ii)$ can be found in \cite[pages 516-517]{Sz02}. In sequence we show that $(ii)\Rightarrow(iii)$, $(iii)\Rightarrow(i)$, $(iii)\Rightarrow (T^\gamma)|_{l_0}$ and $(T^\gamma)|_{l_0}\Rightarrow(i)$.

We turn to prove $(ii)\Rightarrow(iii)$. By $(ii)$, there exist $b, \ \hat{r}, \ \widetilde{c}>0$ such that for large $L$ there are finite subsets $\Delta_L$ with $0\in\Delta _L \subset \{x\in \mathbb Z^d: x\cdot l_0\geq -bL\}\cap \{x\in \mathbb R^d: |x|_2\leq \hat{r}L\}$ and
$$
   P_0\left[X_{T_{\Delta_L}}\notin \partial^+\Delta_L\right]<e^{-\widetilde{c}L^\gamma}.
$$
Furthermore, since we are considering the outer boundary of our discrete sets as $\Delta_L$, by a simple geometric argument we can replace $\Delta_L$ by its intersection with $\{x\in \mathbb Z^d: \,x\cdot l_0<L\}$, so without loss of generality we can assume that $\Delta_L\subset \{x\in \mathbb Z^d: \, x\cdot l_0<L \}$.  Consider the box $\widetilde B_{L, \hat{r},b, l_0}(x)$ defined by
$$
\widetilde B_{L, \hat r, b, l_0}(x):= x + \widetilde{R}\left((-bL,L)\times(-\hat{r}L, \hat{r}L)^{d-1}\right),
$$
where $\widetilde{R}$ is a rotation on $\mathbb R^d$ with $\widetilde {R}(l_0)=e_1$. Its frontal boundary is defined as in \eqref{fb+}. 
We then have that $\Delta_L \subset \widetilde B_{L, \hat{r}, b, l_0}(0)$, and consequently for large $L$,
$$
P_0\left[X_{T_{\widetilde B_{L, \hat{r}, b, l_0}(0)}}\in \partial^+ \widetilde B_{L, \hat{r}, b, l_0}(0)\right]\geq P_0\left[X_{T_{\Delta_L}}\in \partial^+\Delta_L\right]> 1-e^{-\widetilde{c}L^\gamma}.
$$
Notice that if $b\leq1$, we choose $c$ in $(iii)$ as $\hat r$, and we finish the proof.  If $b>1$, we proceed as follows: let $N=bL$ and consider the box
$$
B_{N, \hat{r}(\lfloor b \rfloor + 1) L, l_0}(0)
$$
defined according (\ref{generalboxes}). We introduce for each integer $i\in [1,\lfloor b \rfloor]$ a sequence $(T_i)_{1 \leq i \leq \lfloor b \rfloor}$ of $(\mathcal F_n)_{n\geq0}$-stopping times via
\begin{align}
\label{stoppinti}
T_1=T_{\widetilde{B}_{L, \hat{r}, b, l_0}(X_0)}, \,\, \mbox{ and }  \
T_i=T_1\circ\theta_{T_{i-1}}+T_{i-1}, \mbox{ for $i>1$.}
\end{align}
We also introduce the stopping time $S$ which codifies the unlikely walk exit from box $\widetilde{B}_{L, \hat{r}, b, l_0}(X_0)$ and is defined by:
$$
S:=\inf_{n\geq 0}\{n\geq0:\, X_n\in\partial \widetilde{B}_{L, \hat{r}, b, l_0}(X_0)\setminus \partial^+\widetilde{B}_{L, \hat{r}, b, l_0}(X_0)\}.
$$
The previous definitions imply
\begin{gather}
\nonumber
  P_0\left[X_{T_{B_{N, \hat{r}(\lfloor b \rfloor+1)L, l_0}(0)}}\in \partial^+B_{N, \hat{r}(\lfloor b \rfloor+1)L, l_0}(0)\right]\geq \\
 \label{caixas}
  P_0\left[T_1<S, (T_1<S)\circ\theta_{T_1},\ldots, (T_1<S)\circ \theta_{T_{\lfloor b \rfloor}} \right].
\end{gather}
It is convenient at this point to introduce \textit{boundary sets} $F_i$, $i\in [1, \lfloor b \rfloor]$ as follows
\begin{align*}
  F_1=\partial^+ \widetilde B_{L, \hat{r}, b, l_0}(0)\,  \mbox{ and }  F_i=\bigcup_{y\in F_{i-1}}\partial^+ \widetilde B_{L, \hat{r}, b, l_0}(y), \mbox{ for $i>1$},
\end{align*}
We also define for $i\in [1,\lfloor b \rfloor]$, \textit{good environmental events} $G_i$ via
\begin{equation*}
  G_i=\left\{\omega\in \Omega:\, P_{y,\omega}\left[T_1<S\right]\geq 1 -e^{-\frac{\widetilde{c}L^\gamma}{2}},\,\forall y\in F_i  \right\}.
\end{equation*}
Observe that the left hand of inequality (\ref{caixas}) is greater than
\begin{gather*}
   \mathbb E\left[P_{0,\omega} \left[T_1<S, (T_1<S)\circ\theta_{T_1},\ldots, (T_1<S)\circ \theta_{T_{\lfloor b \rfloor}}\right] \mathds{1}_{G_{\lfloor b \rfloor}}\right].
\end{gather*}
In turn, writing for simplicity $\widetilde B_{L, \hat{r}, b, l_0}(y)$ as $B(y)$ for $y\in \mathbb Z^d$, the last expression equals
\begin{align*}
  &\sum_{y\in F_{\lfloor b \rfloor}}\mathbb E\left[P_{0,\omega}\left[T_1<S,\ldots ,X_{T_{\lfloor b \rfloor}}=y\right]P_{y,\omega}\left[X_{T_{B(y)}}\in \partial^+B(y)\right]\mathds{1}_{G_{\lfloor b \rfloor}}\right]\\
  & \geq \left(1-e^{-\frac{\widetilde{c}}{2}L^{\gamma}}\right)\left(P_0\left[T_1<S,\ldots,(T_1<S)\circ \theta_{T_{\lfloor b \rfloor-1}} \right]-\mathbb P[(G_{\lfloor b \rfloor})^c]\right),
\end{align*}
where we have made use of the Markov property to obtain the last inequality. We iterate this process to see
\begin{align}
\nonumber
&P_0\left[X_{T_{B_{N, \hat{r}(\lfloor b \rfloor+1)L, l_0}(0)}}\in \partial^+B_{N, \hat{r}(\lfloor b \rfloor+1)L, l_0}(0)\right]\\
\label{finalcaixa}
&\geq \left(1-e^{-\frac{\widetilde{c}L^\gamma}{2}}\right)^{\lfloor b \rfloor+1}-\sum_{i=1}^{\lfloor b \rfloor}\left(1-e^{-\frac{\widetilde{c}L^\gamma}{2}}\right)^{\lfloor b \rfloor-i}\mathbb P[G_i^c].
\end{align}
Notice that using Chebyshev's inequality, we have for $i\in [1,\lfloor b \rfloor]$ and large $L$
\begin{equation}\label{Gi}
\mathbb P[G_i^c]\leq\sum_{y\in F_i}\mathbb P\left[P_{y,\omega}\left[X_{T_{B(y)}}\notin \partial^+B(y)\right]>e^{-\frac{\widetilde{c}L^\gamma}{2}}\right]\leq e^{-\frac{\widetilde{c}L^\gamma}{4}}.
\end{equation}
From (\ref{finalcaixa}), the fact that $b$ is finite and independent of $L$ and the estimate (\ref{Gi}), there exists $\eta>0$ such that for large $N$
\begin{equation*}
P_0\left[X_{T_{B_{N, \hat{r}(\lfloor b \rfloor+1)L, l_0}(0)}}\in \partial^+B_{N, \hat{r}(\lfloor b \rfloor+1)L, l_0}(0)\right]\geq 1-e^{-\eta N^\gamma}
\end{equation*}
and this ends the proof of the required implication.

\smallskip

To prove the implication $(iii)\Rightarrow(i)$, we fix a rotation $R$ on $\mathbb R^d$, with $R(e_1)=l_0$ and such that $R$ is the underlying rotation of hypothesis in $(iii)$. For small $\alpha$ we define $2(d-1)$-\textit{directions} $l_{+i}$ and $l_{-i}$, $i\in [2,d]$
\begin{gather*}
  l_{+i}=\frac{l_0+\alpha R(e_i)}{|l_0+\alpha R(e_i)|_2} \,\, \mbox{ and }
  l_{-i}=\frac{l_0-\alpha R(e_i)}{|l_0-\alpha R(e_i)|_2}.
\end{gather*}
Following a similar argument as the one in \cite[Proposition 4.2]{GR17}; 
 but using $\gamma$-stretched exponential decay instead of polynomial one; we conclude that there exists $\alpha > 0$ sufficiently small such that for each $2 \le i \le d$ there are some $r_i>0$, with
\begin{equation}\label{piexp}
\limsup_{\substack{L\rightarrow}\infty}\,L^{-\gamma} \ln P_0\left[X_{T_{B_{L, r_i L, l_{\pm i}}(0)}}\notin \partial^+B_{L, r_i L, l_{\pm i}}(0)\right]<0.
\end{equation}
Observe that from (\ref{piexp}), we finish the proof that $(iii)$ implies $(i)$ by taking
\begin{gather*}
a_0=1,\,a_1=\ldots=a_{2(d-1)}=\frac{1}{2}, \
b_0=\ldots=b_{2(d-1)}=1 ,\\
l_0,\, l_1=l_{+1},l_2=l_{-1},\ldots, l_{2(d-1)-1}=l_{+(d-1)}, l_{2(d-1)}=l_{-(d-1)},
\end{gather*}
and then observing that for $i\in[0,2(d-1)]$ we have
\begin{equation*}
 P_0\left[\widetilde{T}_{-b_iL}^{l_i}<T_{a_iL}^{l_i} \right]\leq P_0\left[X_{T_{B_{L, r_i L,l_{i} }(0)}}\notin \partial^+B_{L, r_i L,l_{i} }(0)\right].
\end{equation*}

Finally, it is straightforward by the previous argument to see that $(iii)$ implies $(T^\gamma)|_{l_0}$, indeed for $\alpha$ as above, the set
$$
\mathcal A_{\alpha}:= \Big\{\ell\in\mathbb S^{d-1}:\ \ell=\frac{l_0+\alpha' R(e_i)}{|l_0+\alpha' R(e_i)|_2},\ \alpha'\in(-\alpha,\alpha),\ i\in[2,d] \Big\}
$$
is open in $\mathbb S^{d-1}$ and contains $l_0$. Furthermore, for any $b>0$ we have
$$
\limsup_{\substack{L\rightarrow\infty}}L^{-\gamma}\ln P_0[\widetilde T_{-bL}^\ell<T_{L}^\ell]<0,
$$
for each direction $\ell\in \mathcal A_{\alpha}$, an thus we have that $(T^\gamma)|_{l_0}$ holds. Conversely, assume $(T^\gamma)|_{l_0}$, then there exits $\varepsilon > 0$ such that
$$
\limsup_{\substack{L\rightarrow\infty}}L^{-\gamma}\ln P_0[\widetilde T_{-bL}^\ell<T_{L}^\ell]<0,
$$
for each direction $\ell\in \left(l_0+\{x\in\mathbb R^d: |x|_2<\varepsilon\}\right)\cap\mathbb S^{d-1}=: \Theta(\varepsilon)$. Since the set $\Theta(\varepsilon)$ is not plane and in fact has a non-zero curvature, any $d$-different elements $l_1$, $l_2$, $\ldots,$ $l_d$ of $\Theta(\varepsilon)$ are linearly independents and thus they span $\mathbb R^d$. Take data: $l_0,$ $l_1,$
$\ldots,$ $l_d$, $a_0=1,$ $a_1=1,$ $\ldots,$ $a_d=0$ and arbitrary positive numbers $b_0,$ $b_1,$ $\ldots,$ $b_d$ and it is clear that they generate an $l_0$-directed system of slabs of order $\gamma$, which ends the proof. \hfill $\square$

\end{document}